\documentclass[11pt]{article}

\usepackage[latin1]{inputenc}
\usepackage[T1]{fontenc}
\usepackage{amsmath,amssymb,amsthm}

\DeclareMathOperator{\depth}{depth} \DeclareMathOperator{\im}{im}
\DeclareMathOperator{\coker}{coker} 
\DeclareMathOperator{\HH}{H} \DeclareMathOperator{\Proj}{Proj}
\DeclareMathOperator{\Spec}{Spec} \DeclareMathOperator{\Hom}{Hom}
\DeclareMathOperator{\Tor}{Tor} \DeclareMathOperator{\Ext}{Ext}
 
\DeclareMathOperator{\rank}{rank} \DeclareMathOperator{\Hilb}{Hilb}
\DeclareMathOperator{\GradAlg}{GradAlg} \DeclareMathOperator{\PGor}{PGor}
\DeclareMathOperator{\Gor}{Gor}
\DeclareMathOperator{\ext}{ext}
\DeclareMathOperator{\reg}{reg}
\DeclareMathOperator{\DD}{D} 
\newtheorem{theorem}{Theorem}
\newtheorem{lemma}[theorem]{Lemma}
\newtheorem{corollary}[theorem]{Corollary}
\newtheorem{definition}[theorem]{Definition}
\newtheorem{example}[theorem]{Example}
\newtheorem{proposition}[theorem]{Proposition}

\newtheorem{remark}[theorem]{Remark}

\usepackage{latexsym}            % for the qed symbol
\usepackage{amssymb}

\newcommand{\pp}{{\mathbb P}}

\newcommand\sE{{\mathcal E}}

\newcommand\sH{{\mathcal H}}
\newcommand\sI{{\mathcal I}}

\newcommand\sM{{\mathcal M}}
\newcommand\sN{{\mathcal N}}
\newcommand\sO{{\mathcal O}}

\newcommand{\proj}[1]
{ \mathchoice
           { {\mathbb P}^{#1} }
           { {\mathbb P}^{#1} }
           { {\mathbb P}^{#1} }
           { {\mathbb P}^{#1} }
         }

\begin{document}
\bibliographystyle{plain}

\title{Unobstructedness and dimension of families of Gorenstein algebras}

\author{Jan O. Kleppe}

\date{\ }

\maketitle

\vspace*{-0.75in}
\begin{abstract}
  \noindent The goal of this paper is to develop tools to study maximal
  families of Gorenstein quotients $A$ of a polynomial ring $R$. We prove a
  very general Theorem on deformations of the homogeneous coordinate ring of a
  scheme $\Proj(A)$ which is defined as the degeneracy locus of a regular
  section of the dual of some sheaf $\widetilde M$ of rank $r$ supported on
  say an arithmetically Cohen-Macaulay subscheme $\Proj(B)$ of $\Proj(R)$.
  Under certain conditions (notably; $M$ maximally Cohen-Macaulay and
  $\wedge^r{\widetilde { M}} \simeq \widetilde {{ K}}_{ B}(t)$ a twist of the
  canonical sheaf), then $A$ is Gorenstein, and under additional assumptions,
  we show the unobstructedness of $A$ and we give an explicit formula the
  dimension of any maximal family of Gorenstein quotients of $R$ with fixed
  Hilbert function obtained by a regular section as above. The theorem also
  applies to Artinian quotients $A$.

  The case where $M$ itself is a twist of the canonical module ($r=1$) was
  studied in a previous paper, while this paper concentrates on other low rank
  cases, notably $r=2$ and $3$. In these cases regular sections of the first
  Koszul homology module and of normal sheaves to licci schemes (of say
  codimension 2) lead to Gorenstein quotients (of e.g. codimension 4) whose
  parameter spaces we examine. Our main applications are for Gorenstein
  quotients of codimension 4 of $R$ since our assumptions are almost always
  satisfied in this case. Special attention are paid to arithmetically
  Gorenstein curves in $\pp^5$.

% 14C05 Param (Chow and HIlb schemes),
% 13D10 Deform and inf methods
% 13D03 (Co)homology of comm.ringds and alg.
% 14M05 Varieties defined by ring conditions
% 14F05 Vector bundles, sheaves, related constructions
% 14J40 n-folds
% 3D02 Syzygies and resolutions
% 13C14 Cohen-Macaulay modules

 \end{abstract}
 \vspace*{-0.25in}
{\small
 \tableofcontents}
\thispagestyle{empty}

\section{Introduction}
In this paper we study deformations of arithmetically Gorenstein (AG)
subschemes $X= \Proj(A)$ of a projective space $\pp^{N}=\Proj(R)$. In
particular we consider deformations of schemes $X$ defined as the degeneracy
locus of a regular section of some maximal Cohen-Macaulay sheaf $\widetilde M$
supported on an arithmetically Cohen-Macaulay (ACM) subscheme $Y= \Proj(B)$ of
$\pp^{N}$. Our results allow us to study well three natural parameter spaces
classifying AG subschemes of $\pp^{N}$ (resp. graded Gorenstein quotients of
$R$) with Hilbert polynomial $p$ (resp. Hilbert function $H$). The three
schemes we have in mind are the open subscheme of Grothendieck's Hilbert
scheme $\Hilb^{p}(\pp^{N})$ consisting of AG subschemes of positive dimension,
the open part of the postulation Hilbert scheme $\Hilb^H(\pp^{c})$ consisting
of AG zeroschemes with Hilbert function $H$ and Iarrobino-Kanev's
determinantal scheme $\PGor(H)$ classifying graded Artinian Gorenstein
quotients of $R$ via the Macaulay correspondence. A common denominator of
these parameter spaces is the $k$-scheme $ \GradAlg^{H}(R)$ of \cite{K98}
which parametrizes graded quotients $B =R/I_B$ such that $\ \depth B \geq
\min(1,\dim B)$ and $H_B=H$, $H_B$ the Hilbert function of $B$. Its open
subscheme consisting of Gorenstein quotients $B$ of codimension $c$ in $R$,
$\Gor_c^H(R)$, is naturally {\it isomorphic} to the three schemes in the
mentioned order provided the degree of $p$ is greater or equal $1$, equal $0$
and equal $-1$ (i.e. $p(t)=0$) respectively. The isomorphism is even
scheme-theoretic except possibly for $\PGor(H)$ where the isomorphism is at
least infinitesimal and topological, making questions of dimension and
smoothness of $\PGor(H)$ equivalent to the corresponding questions for
$\GradAlg(H):=\GradAlg^{H}(R)$, cf. \cite{K03} for details and \cite{HS} for a
generalization of $\GradAlg(H)$.

If the codimension $c \leq 3$, the smoothness, the irreducibility and the
dimension of $\Gor_c^H(R)$ are known (\cite{K98}, \cite{GPS}, \cite{D96},
\cite{KR98}). Hence we focus on Gorenstein quotients of codimension $4$ of $R$
which, nowadays, get quite a lot of attention (\cite{GKS}, \cite{IS},
\cite{K03}, \cite{IK}, \cite{B99}). Moreover since some of our results require
that the degree of $p$ is greater or equal to $1$, we mostly consider AG curves
in $\pp^{5}$, having also the case $\PGor(H)$ in mind since we have tried very
much, and partially succeeded, to generalize our results to the Artinian case.

Let $B \rightarrow A$ be a morphism of graded quotients of $R$. In the
background section we recall a basic result of how we can deduce the
smoothness and the dimension of $ \GradAlg^{H_A}(R)$ from the corresponding
properties of $ \GradAlg^{H_B}(R)$ (Theorem~\ref{varunobstr} and
Corollary~\ref{corvarunobstr}). The assumptions of Theorem~\ref{varunobstr}
are satisfied for a complete intersection (c.i.) $B$ and for some other
classes of quotients with nice properties (licci is often enough) except for
the vanishing ``$_0\!\HH^2(B,A,A)=0$'' of the second algebra cohomology group.
Theorem~\ref{varunobstr} and Corollary~\ref{corvarunobstr} are enough to treat
the rank $r= 1$ case of $M$ satisfactorily (\cite{K03}), while, for higher
rank cases, they are insufficient. A main result, Theorem~\ref{gorgenth}, of
this paper replaces the assumption ``$_0\!\HH^2(B,A,A)=0$'' with assumptions
which are more natural to the application we have in mind; deforming schemes
$X$ defined as the degeneracy locus of a regular section. Note, however, that
Theorem~\ref{varunobstr} and Corollary~\ref{corvarunobstr} apply also to
quotients %in which $B \rightarrow A$
which not necessarily are given by such a regular section, and we take the
opportunity to include a result which generalizes \cite{IK}, Thm.\! 4.17
(Corollary~\ref{corvarunobstr2}).

To state the assumptions of Corollary~\ref{corvarunobstr} and
Theorem~\ref{gorgenth} and its dimension formulas in a more computable form we
consider low rank cases on a licci scheme $Y=\Proj(B)$. Let $K_B$ be the
canonical module of $B$. If $M$ is a graded maximal Cohen-Macaulay $B$-module
of rank $r=2$ such that $\widetilde {M}\arrowvert_U $ is locally free and
$\wedge^2{\widetilde { M}}\arrowvert_U \simeq \widetilde {K_B}(t)\arrowvert_U$
in a large enough open set $U=Y-Z$ of $Y$, then a regular section $\sigma$ of
$\widetilde{M}^*(s)\arrowvert_U$, $M^* = \Hom_B(M,B)$, defines a graded
Gorenstein quotient $A$ given by the exact sequence
\begin{equation} \label{introMs}
 {\footnotesize
0 \rightarrow K_B(t-2s)\rightarrow M(-s)  \stackrel{\sigma}{\rightarrow}  B
\rightarrow A \rightarrow 0 \ ,
} 
\end{equation}
similar to what happens in the usual Hartshorne-Serre correspondence
(\cite{Ha}, Thm.\! 4.1). Moreover $M \simeq \Hom_B(M,K_B(t))$
(Theorem~\ref{mainth}). Let $S_2(M)$ be the second symmetric power of $M$.
Using small letters for the $k$-dimension of the $v$-graded piece
$_v\!\Ext_B^i(-,-)$ of $\Ext_B^i(-,-)$, we define $\gamma(S_2M)_v$ and
$\delta(Q)_v$ by
$$
\delta(Q)_{v} =\ _{v}\!\hom_B(I_B/I_B^2,Q)-\ _{v}\!\ext_B^1(I_B/I_B^2,Q) \ ,
\ \ {\rm and} $$
$$\gamma(S_2M)_v =\ _{v}\!\hom_B(S_2(M),K_B(t))-\ 
_{v}\!\ext_B^1(S_2(M),K_B(t)) \ . 
$$
Note that if $\depth_{I(Z)}B \geq 3$ (resp. $\depth_{I(Z)}B \geq 4$) and
$char(k) \neq 2$, one may show that $\gamma(S_2M)_0 =\ _{0}\!\hom_B(M,M)-\
_{0}\!\ext_B^1(M,M)-1$ (resp. $_{0}\!\Ext_B^i(S_2(M),K_B(t)) \simeq \
_{0}\!\Ext_B^i(M,M)$ for $i=1,2$) since $r=2$, cf. \eqref{M,M}. A main result
of Section 2 is the following Theorem (Theorem~\ref{Mmainth}) which applies
also to an Artinian $A$.

\begin{theorem} \label{introMmainth} Let $B = R/I_B$ be a graded licci
  quotient of $R$, %of codimension $c$
  let $M$ be a graded maximal Cohen-Macaulay B-module, and suppose
  $\widetilde{M}$ is locally free of rank $2$ in $U:=\Proj(B) -Z$, that $\dim
  { B}-\dim {B}/{I}(Z)\geq 2$ and $\wedge^2{\widetilde { M}}\arrowvert_U
  \simeq \widetilde {K_B}(t)\arrowvert_U$. Let $A$ be defined by a regular
  section $\sigma$ of $\widetilde{M}^*(s)$ on $U$, i.e. given by
  \eqref{introMs}, and suppose  $\  _s\!\Ext^1_B(M,B)=0$. \\[2mm]
  \indent {\rm A)} \ If $\ \ \ \ _t\!\Ext_B^2(S_2(M),K_B)=0 \ $ and $\
  _{-s}\!\Ext^2_B(I_B/I_B^2,M)= 0 \ $,   \\[2mm]
  then $A$ is unobstructed as a graded $R$-algebra (indeed
  $_0\!\HH^2(R,A,A)=0$), $A$ is Gorenstein, and {\small $$ \ \dim_{(A)}
    \GradAlg(H_A)= \ \dim (N_B)_0 + \ \dim (M^*)_s - 1 - \gamma(S_2M)_0 + \
    \dim(K_B)_{t-2s} +\delta(K_B)_{t-2s} - \delta(M)_{-s}.$$ } Moreover if
  $char(k)=0$ and $(B \rightarrow A)$ is general with respect to $\
  _{0}\!\hom_R(I_B,I_{A/B}) $, then the codimension of the stratum of
  quotients given by \eqref{introMs} around $(A)$ is $$
  _{-s}\!\ext^1_B(I_B/I_B^2,M) - \dim (\im \beta)$$ where $\beta$ is the
  homomorphism $\ _{-2s}\!\Ext_B^1(I_B/I_B^2,K_B(t)) \rightarrow \
  _{-s}\!\Ext_B^1(I_B/I_B^2,M)$ induced by \eqref{introMs}.
  \\[2mm]
  \indent {\rm B)} \ If \ $(M,B)$ is unobstructed along any graded deformation
  of $B$ and $ _{-s}\!\Ext^1_B(I_B/I_B^2,M)= 0$, then $A$ is Gorenstein and
  $H_B$-generic. Moreover $A$ is unobstructed as a graded $R$-algebra and the
  dimension formula for $\ \dim_{(A)} \GradAlg^{H_A}(R)$ of part {\rm A)}
  holds.
\end{theorem}

Here $N_B = \Hom_B(I_B/I_B^2,B)$, $I_{A/B}:= \ker (B \to A)$ and ``$(M,B)$
unobstructed along any deformation of $B$'' means that for every deformation
$(M_S,B_S)$ of $(M,B)$, $S$ local and Artinian, there is a deformation of
$M_S$ to any deformation of $B_S$ (cf. Definition~\ref{unobstrdef}). Moreover
``$(B \rightarrow A)$ is general with respect to $\gamma:=\
_{0}\!\hom_R(I_B,I_{A/B}) $'' if and only if $ \
_{0}\!\hom_R(I_{B'},I_{A'/B'}) \ge \gamma $  for every $(B' \rightarrow A')$
in an {\it open neighbourhood} of $(B \rightarrow A)$ in $\GradAlg(H_B,H_A)$
(see Remark~\ref{mre}$(ii)$). Furthermore ``$A$ is $H_B$-generic'' essentially
means that the codimension of the stratum mentioned in
Theorem~\ref{introMmainth}A) is zero, see the text before
Theorem~\ref{varunobstr} and Theorem~\ref{Mmainth} for precise definitions.

\begin{remark} \label{introrem} If $ \max n_{2,j}$ is the largest degree of a
  minimal relation of $I_B$ and $a$ is an integer such that $M_v = 0$ for $v
  \leq a$, then $$\ _{-s}\!\Ext^1_B(I_B/I_B^2,M)=0 \ \ {\rm provided} \ \ s
  \geq \max n_{2,j} - a.$$ In this case $ \dim(K_B)_{t-2s} +\delta(K_B)_{t-2s}
  - \delta(M)_{-s} = 0$. This significantly simplifies the dimension formula
  of Theorem~\ref{introMmainth}B) provided we use the theorem for such $s$.
\end{remark}
 
On a Cohen-Macaulay (CM) quotient $B=R/I_B$ of codimension 2 in $R$ with
minimal resolution
\begin{equation*}
 0 \rightarrow  \oplus_{j=1}^{\mu -1} R(-n_{2,j}) \rightarrow
  \oplus_{i=1}^{\mu} R(-n_{1,i}) \rightarrow R \rightarrow B \rightarrow 0  \ ,
\end{equation*}
Theorem~\ref{introMmainth} applies to $M=N_B$, the normal module, and to
$M=H_1$, the 1. Koszul homology of $I_B$. In Section 3 we study these
applications in detail. To see that the assumptions of
Theorem~\ref{introMmainth} are satisfied, let $\Proj(B) -Z \hookrightarrow
\pp^N$ be a local complete intersection (l.c.i.). Then we prove
 
\begin{proposition} \label{introprop}
  Let $\ B=R/I_B \ $ be a codimension two CM quotient of $R$ and suppose
  $\  \dim { B}-\dim {B}/{I}(Z) \geq 2 $. Then 
  
  (i) $\ \Ext_B^2(I_B/I_B^2,H_1)= \ \Ext_B^1(H_1,B)=0$ and $S_2(H_1)$ is a
  maximal CM $B$-module.
 
  (ii) $\ \Ext_B^i(I_B/I_B^2,I_B/I_B^2)=0 \ \ {\rm for} \ 1 \leq i \leq
   \dim { B}-\dim {B}/{I}(Z) . $
\end{proposition}

Using Proposition~\ref{introprop}$(i)$ we get the following result, with a
sufficiently weak assumption on $\depth_{I(Z)}B = \dim { B}-\dim {B}/{I}(Z)$
so that the result applies also to an Artinian $A$.

\begin{corollary} \label{introcorH} Let $B=R/I_B$ be a graded codimension two
  CM quotient of $R$, let $U= \Proj(B)-Z \hookrightarrow \pp^N$ be an l.c.i.
  and suppose the number of minimal generators of $I_B$ is $\mu =4$ and
  $\depth_{I(Z)}B \geq 2$. If $A$ is defined by a regular section of
  $\widetilde{H_1}^*(s)$ on $U$, then $A$ is unobstructed as a graded
  $R$-algebra (indeed $\ \HH^2(R,A,A)=0$), $A$ is Gorenstein of codimension 4
  in $R$, and {\small
  $$     \ \dim_{(A)} \GradAlg(H_A)= \ \dim (N_B)_0 + \ \dim (H_1^*)_s - 1 -
  \gamma(S_2H_1)_0 + \ \dim(K_B)_{t-2s} +\delta(K_B)_{t-2s} - \delta(H_1)_{-s}
  $$ } where $t = \dim R - \sum n_{1,i}$. Moreover if $\ s > \max n_{2,j}-\min
n_{2,j} $, then $A$ is $H_B$-generic, and $$\ \ \dim_{(A)} \GradAlg^{H_A}(R)=
\ \dim (N_B)_0 + \ \dim (H_1^*)_s - 1 -\ _0\!\hom_B( S_2(H_1), K_B(t)).$$
\end{corollary}

Furthermore in the situation of Corollary~\ref{introcorH} we show how to
compute every term of the dimension formulas and we give several examples. 

Applying Theorem~\ref{introMmainth} to the normal module $N_B$ we need the
vanishing of the $\Ext$-groups of Proposition~\ref{introprop}$(ii)$ to prove
that the assumptions of Theorem~\ref{introMmainth} hold. To do so we have to
increase $\depth_{I(Z)} B$. Hence the following result does not automatically
apply to Gorenstein quotients $A$ of dimension less or equal to $1$ unless we
in a given example are able to verify the assumptions of
Theorem~\ref{introMmainth}. In Section 3 we give examples of the result.
% and we show how to compute the numbers of the dimension formula.
   
\begin{corollary} \label{introcorNB} Let $B=R/I_B$ be a graded codimension two
  CM quotient of $R$, let $U= \Proj(B)-Z \hookrightarrow \pp^N$ be an l.c.i.
  and suppose $\depth_{I(Z)} B \geq 4$. If $A$ is given by a regular section
  of $\widetilde{N_B}^*(s)$ on $U$, and if $s > 2 \max n_{2,j}- \min n_{1,i}$
  and $char(k) \neq 2$, then $A$ is $H_B$-generic and unobstructed as a graded
  $R$-algebra. Moreover $A$ is Gorenstein of codimension 4 in $R$, and letting
  $X=\Proj(A)$ and $\eta(v):= \dim (I_B/I_B^2)_v$, we have
  $$
  \ \dim_{(A)} \GradAlg^{H_A}(R)= \ \dim_{(X)} \Hilb^{p}(\pp^{N})=
  \eta(s)+ \sum_{j=1}^{\mu -1} \eta( n_{2,j})-\sum_{i=1}^{\mu} \eta( n_{1,i})
  \ . $$
\end{corollary}

Finally we prove a Theorem, similar to Theorem~\ref{introMmainth}, in the rank
$r = 3$ case which also applies to an Artinian $A$ (Theorem~\ref{Mmainth3}).
It admits the following Corollary,  

\begin{corollary} \label{introcorH3} Let $B=R/I$ and $U= \Proj(B)-Z$ be as in
  Corollary~\ref{introcorH} and suppose that $\mu=5$, $char(k) \neq 2$ and
  $\depth_{I(Z)}B \geq 3$. If $A$ is defined by a regular section of
  $\widetilde{H_1}^*(s)$ on $U$, then $A$ is unobstructed as a graded
  $R$-algebra (indeed $\ \HH^2(R,A,A)=0$), $A$ is Gorenstein of codimension 5
  in $R$, and $ \ \dim_{(A)} \GradAlg^{H_A}(R)= $
  $$
  \ \dim (N_B)_0 + \ \dim (H_1^*)_s + \ _{-s}\!\hom_B(S_2(H_1), K_B(t)) - \ 
  _0\! \hom_B(H_1, H_1)- \dim(K_B)_{t-3s} - \delta \ ,$$
  where $ \delta:=
  \delta(H_1)_{-s} + \delta(K_B)_{t-3s} - \delta(H_2)_{-2s}$ and $t = \dim R -
  \sum n_{1,i}$.  If in addition $\ s > \max n_{2,j}-\min n_{2,j} $, then
  $A$ is $H_B$-generic, $\ \delta = 0$ and $ \dim(K_B)_{t-3s}=0$.
\end{corollary}

We thank Chris Peterson and Johannes Kleppe for useful discussions. Indeed we
are in this work inspired by the joint works \cite{KP1} and \cite{KP2} where
we in various ways construct Gorenstein quotients of $R$. This paper uses
deformation theory to vary every object and morphism of these constructions,
to see how large the corresponding stratum of its parameter space will be. As
we see by the results above, the stratum is often a component. In some cases,
however, it is a proper stratum of a codimension which we make explicit.
Moreover we thank Rosa M. Mir\'o-Roig for interesting discussions on
$\PGor(H)$. In particular Corollary~\ref{corvarunobstr2} was included in the
paper after a discussion with her.
  
\subsection{Background}
We largely keep the notations of \cite{K03}, and we recommend the section of
Preliminaries of that paper for an overview of cohomology groups and
deformation theory. Let $B$ be an n-dimensional graded quotient of a
polynomial ring $R$ in $n+c$ variables (of degree $1$) over an
algebraically closed field $k$, and let $M$ and $N$ be finitely generated
graded $B$-modules. Let $\depth_{J}{M}$ denote the length of a maximal
$M$-sequence in a homogeneous ideal $J$ and let $\depth {M} =
\depth_{\mathfrak m}{ M}$ where $\mathfrak m$ is the irrelevant maximal ideal.
Let $\HH^i_{J}(-)$ be the right derived functor of the functor,
$\Gamma_{J}(-)$, of sections with support in $ \Spec (B/J)$. Then
%. Recall that $\depth_{J}M = \inf \{\depth M_{\wp}\ \arrowvert \ \wp
%\supseteq J,\ \wp \ {\rm graded} \}$ 
\begin{equation} \label{ONE} \depth_{J}{M} \geq r {\ \ \rm if \ and \ only \
    if \ \ (iff) \ \ } \HH^i_{J}({M})=0 {\ \rm for \ } i<r.
\end{equation} 
(cf. \cite{HARLOC}). Let $Y=\Proj ({B})$ and let $Z$ be closed in $Y$ and
$U=Y-Z$. Put $\HH^0_*(U,\widetilde {M}):= \oplus_v \HH^0(U, \widetilde
{M}(v))$. Then we have an exact sequence $0 \rightarrow \HH^0_{I(Z)}(M)
\rightarrow M \rightarrow \HH^0_*(U,\widetilde{M}) \rightarrow \HH^1_{I(Z)}(M)
\rightarrow 0$ and isomorphisms $\HH^i_{{I}(Z)}({M}) \simeq
\HH^{i-1}_*(U,\widetilde {M})\ {\rm for} \ i \geq 2$. If $\depth_{I(Z)}N \geq
i+1$, then the graded group $\Ext_B^i(M,N)$ injects into the corresponding
global $\Ext_{\sO_U}^i$-group of sheaves. Indeed we have by \cite{SGA2}, exp.
VI, an exact sequence
\begin{equation} \label{twoo}
 \small
 _0\!\Ext_B^i(M,N) \hookrightarrow \Ext_{\sO_U}^i(\widetilde M
 \arrowvert_U,\widetilde N \arrowvert_U)  \rightarrow \ _0\!\Hom_B(M,
 \HH^{i+1}_{I(Z)}(N)) \rightarrow \  _0\!\Ext_B^{i+1}(M,N)  
  \rightarrow \Ext_{\sO_U}^{i+1}(\widetilde M \arrowvert_U,\widetilde
  N\arrowvert_U) 
\end{equation}   
where the form of the middle term comes from a spectral sequence discussed in
\cite{SGA2}.

A Cohen-Macaulay (resp. maximal Cohen-Macaulay) $B$-module satisfies $\depth
{M}= \dim {M}$ (resp. $\depth {M}= \dim {B}$) by definition, or equivalently,
$\HH_{\mathfrak m}^i({ M})=0\ {\rm for}\ i<\dim {M}$ (resp. $i<\dim {B}$). If
$B$ is Cohen-Macaulay, then the $v$-graded piece of $\HH_{\mathfrak m}^i({
  M})$ is via Gorenstein duality given by $_{v}\!\HH_{\mathfrak{m}}^{i}(M)
\simeq \ _{-v}\!\Ext_B^{n-i}(M,K_B)^{\vee}$ where $K_{B}=
\Ext^{c}_{R}(B,R(-n-c))$ is the canonical module of $B$. In this case
$B=R/I_B$ has a minimal $R$-free resolution of the following form % \cite{EIS}
\begin{equation} \label{complex2K}
  0 \rightarrow {G}_c \rightarrow . . . 
  \rightarrow {G}_1 \rightarrow R \rightarrow B  \rightarrow 0  \hspace{0.2 cm}
  ,  \hspace{1 cm} G_j = \oplus_{i=1}^{r_j} R(-n_{j,i}) 
\end{equation}
and the $R$-dual sequence $ 0 \rightarrow R \rightarrow {G}_1^* \rightarrow
... \rightarrow {G}_c^* \rightarrow K_B(n+c) \rightarrow 0$ is exact. Moreover
the Castelnuovo-Mumford regularity of $I_B$ is $\reg(I_B)= \max\{ n_{c,i} \}
-c+1$. A coherent non-trivial $\mathcal{O}_Y$-module, $\sM$, is called a
maximal Cohen-Macaulay sheaf on $U$ if $\sE xt^i_{\mathcal
  O_Y}(\sM,\omega_Y)\arrowvert_U = 0$ for $i>0$ where $\omega_Y =
\widetilde{K}_{B}$.

\begin{lemma} \label{PRELEM}
  Let $B$ be Cohen-Macaulay. Let $r$ and $t$ be
  integers. Let $J \subseteq B$ be an ideal satisfying
  $\depth_{J}B \geq r$ and let $M$ be a finitely generated
  $B$-module satisfying $\depth_{\mathfrak m} M \geq \dim B -t$.
  Then $\depth_{J}M \geq r-t$.
\end{lemma}
For a proof, see \cite{KP2}, Lem.\!\ 5.  (What there looks like sheaf-$\sE
xt^j$ is actually $ \Ext^j$ as $B$-modules).

Now we recall a main result of \cite{KP1} in which $B$ is Cohen-Macaulay.
Indeed a main idea of this paper is to use deformation theory to see how
general the construction of Gorenstein quotients, given by the Theorem below,
is. Note that in this case $\depth_{ J}B= \dim{ B}-\dim{B}/{ J}$ for any
ideal $ J$ of $B$, cf. \cite{EIS}.

\begin{theorem}\label{mainth}
  Let $R$ be a polynomial ring, let $B=R/I_B$ be a codimension $c$ graded CM
  quotient of $R$ and let $M$ be a finitely generated graded maximal CM
  $B$-module. Let $Y=\Proj(B)$, let $Z$ be a closed scheme such that ${\rm
    dim}({B})-{\rm dim}({B}/{I}(Z)) \geq {\rm max}(r,2)$ and let $U=Y-Z$. Let
  ${ M}_i = \HH^0_*(U,\wedge^i \widetilde { M})$ for $i\geq 0$, and suppose
  that ${\widetilde { M}}\arrowvert_U$ is locally free (of rank $r$) and that
 $$\wedge^r{\widetilde { M}}\arrowvert_U \simeq \widetilde {{ K}}_{
   B}(t)\arrowvert_U \ \ {\rm for \ some \ integer \ } t \ .$$
 Moreover suppose
 ${ M}_i$ is a maximal CM ${B}$-module for $2 \leq i \leq r/2$.
 Then any regular section, $\sigma \in \HH^0(U,\widetilde{ M}^*(s))$ ($s$ an
 integer), defines a Gorenstein quotient ${ R} \twoheadrightarrow { A}$ of
 codimension $r+c$ given by the exact sequence
\begin{equation} \label{LESTHM}
   0\rightarrow {M}_r(-rs) \rightarrow {M}_{r-1}((1-r)s) \rightarrow  \dots
   \rightarrow  {M}_2(-2s) 
   \rightarrow {M}(-s) \stackrel{\sigma}{\rightarrow} 
   {B} \rightarrow {A} \rightarrow 0.
\end{equation}
Indeed all $M_i, 0 \leq i \leq r$, are maximal CM $B$-modules, and
$M_{r-i} \simeq \Hom_B(M_i,K_B)(t)$ for $0 \leq i \leq r$.
\end{theorem}

\begin{remark} \label{remmainth} By replacing $M$ by $M_1$ in \eqref{LESTHM} we
  see that the maximal CM assumption on $M$ is really superfluous in the rank
  2 Artinian case ($r = \dim(B) = 2$). In \cite{KP2} we succeeded generalizing
  Theorem~\ref{mainth} by weakening the maximal CM assumptions on $M$ and
  ${M_i}$ in several ways.
 \end{remark} 
 
 In \cite{KP1} we applied Theorem~\ref{mainth} to different modules $M$ on
 licci schemes (i.e. schemes whose homogeneous coordinate ring is in the
 linkage class of a complete intersection, cf.  \cite{MIG} for a survey),
 notably the normal module, $N_B$, and the 1. Koszul homology,
 $\HH_1=\HH_1(I_B)$, built on a set of minimal generators of $I_B$. Note that
 the module $\HH_1$ is essentially given by an exact sequence
\begin{equation} \label{es1}
   0 \rightarrow \HH_2(R,B,B) \rightarrow  {\HH}_1 \rightarrow
  {G}_{1} \otimes_{R} {B} \rightarrow {I_B}/{I_B}^2 \rightarrow 0.
\end{equation}
in which $\HH_2(R,B,B)$ is the 2. algebra homology \cite{VAS}. An ideal $I_B$
of $R$ is called {\it syzygetic} if $\HH_2(R,R/I_B,R/I_B) = 0$. Using
\eqref{es1} one shows that if $R \rightarrow B$ is generically a complete
intersection and $B$ is licci, then $\HH_2(R,B,B)=0$ because $\HH_1$ is a
maximal CM module in the licci case \cite{H82}.  If $B$ is licci one also
knows that $N_B$ is a maximal CM $B$-module and that $\Ext_B^1( {I_B}/{I_B}^2,
B)=0$ (\cite{BU83} and \cite{KMMNP}). Note that if $\mu({I_B})$ is the number
of minimal generators of $I_B$, we see from \eqref{es1} that $\rank{ H}_1 =
\mu({ I_B})-c$ since the rank of $\sN_Y = \widetilde N_B = \widetilde {{ I_B}
  / { I_B}^2}^*$ is the same as the codimension, $c$, of $Y$ in $\Proj(R)$. We
proved
 
\begin{proposition} \label{propnormal}
  Let $B= R/I_B$ be licci of codimension $c \geq 2$, and let $Y =
  \Proj({B}) \hookrightarrow \pp^{n+c-1}=\Proj({R})$ be a local complete
  intersection (l.c.i.) in some open $U = Y-Z$.
 \begin{enumerate}
 \item If ${\rm dim} B -{\rm dim} B/I(Z) \geq c$, then
   $\HH^0_{*}(U,\wedge^i \sN_Y)$ are maximal CM 
 $B$-modules for every  $1 \leq i \leq c$, and any regular section of 
  $\HH^0(U,{\sN}_Y^*(s))$ defines a Gorenstein quotient ${  R}
  \twoheadrightarrow {  A}$ of codimension $2c$.
\item If ${\rm dim} B -{\rm dim} B/I(Z) \geq {\rm max}(2, \mu({I_B})-c)$, then
  $\HH^0_{*}(U,\wedge^i \widetilde{H}_1)$ are maximal CM $
  B$-modules for every $1 \leq i \leq \mu({I_B})-c$ and any regular section of
  $\HH^0(U,\widetilde{H}_1^*(s))$ defines a Gorenstein quotient ${R}
  \twoheadrightarrow {A}$ of codimension $\mu({I_B})$.
   \end{enumerate}
 \end{proposition}
 For the graded group $\HH^2(R,B,B)$ we just remark that there is an exact
 sequence
\begin{equation} \label{three}
    0 \rightarrow \ _0\!\Ext_B^1(I_B/I_B^2,N)  \rightarrow \ _0\!\HH^2(R,B,N)
 \rightarrow \ _0\!\Hom_B(\HH_2(R,B,B),N) \rightarrow  
 \end{equation}   
 induced from some well known spectral sequence, cf.  \cite{AND}, Prop.\!\ 
 16.1.
 
 Let $p \in \mathbb Q [t]$ be a non-trivial polynomial and let $ \pp =
 \pp^{n+c-1}$. Then the scheme $\GradAlg(H) = \GradAlg^H (R)$ which we study
 in this paper is the stratum of Grothendieck's Hilbert scheme $ \Hilb^{ p}
 (\pp)$ (cf. \cite{SB}) consisting of points $(Y \subset \pp)$ with Hilbert
 function $ H_Y = H$ (i.e. its corresponding functor deforms the {\em
   homogeneous coordinate ring}, $B$, of $Y$ flatly as a graded $R$-algebra),
 cf. \cite{K98}. Note that we define the Hilbert function of $Y$, or $B$, by
 $H_Y(v)=H_B(v):= \dim B_v$. $\GradAlg^H(R)$ has a natural scheme structure
 whose tangent (resp.\ ``obstruction'') space at $(Y \subset \pp)$ is $\
 _0\!\Hom_B(I_{B}/I_B^2,B) \simeq \ _0\!\Hom_R(I_{B},B) \ $ (resp.\ $_0\!\HH^2
 (R,B,B)$) \cite{K79}. Since $H(v)$ does not vanish for large $v$ (i.e. $B$ is
 non-Artinian), we may look upon $\GradAlg^H(R)$ as parametrizing graded
 $R$-quotients, $R \rightarrow B$, satisfying $\depth_{\mathfrak m}B \geq 1$
 and with Hilbert function $H_B = H$. If $B$ is Artinian, then $\GradAlg^H(R)$
 still represents a functor parametrizing graded $R$-quotients with Hilbert
 function $H_B = H$. It contains an open subscheme of Gorenstein quotients
 which, at least topologically and infinitesimally, coincides with $\PGor(H)$,
 the corresponding scheme of forms with ``catalecticant structure'' (see
 \cite{K98} or \cite{K03} for details). $B$ is called {\em unobstructed} as a
 graded $R$-algebra if $\GradAlg^H(R)$ is smooth at $(R \rightarrow B)$.
 Similarly a closed subscheme $Y$ of $\pp$ is called {\em unobstructed} if
 $\!\ \Hilb^{p}(\pp)$ is smooth at $(Y \subset \pp)$. By \cite{K79}, Thm. 3.6
 and Rem. 3.7,
\begin{equation} \label{Grad}
\GradAlg^H(R)
\simeq \Hilb^{ p}(\pp) \ \ \ {\rm at} \ \ \  (Y \subset \pp) 
\end{equation}
provided $\ _0\!\Hom_R (I_B,\HH_{\mathfrak m}^1(B)) = 0$ (e.g. provided
$\depth_{\mathfrak m}B \geq 2$, in which case \eqref{Grad} also follows from
\cite{El}). In particular if the degree of $p$ is positive, then the open
subschemes of $\GradAlg^H(R)$ (resp. of $\!\ \Hilb^{p}(\pp)$) of Gorenstein
quotients (resp. of AG subschemes) are isomorphic as {\it schemes}.

Similarly we let $\GradAlg(H_B,H_A)$ be the representing object of the functor
deforming flags (surjections) $B \rightarrow A$ of graded quotients of $R$ of
positive depth (if $B$ and/or $A$ are non-Artinian) and with Hilbert functions
$H_B$ and $H_A$ of $B$ and $A$ respectively.  Let $p$ be the second projection
$$p:\GradAlg(H_B,H_A) \rightarrow \GradAlg^{H_A}(R)$$ induced by sending $(B'
\rightarrow A')$ onto $(A')$, and let $q:\GradAlg(H_B,H_A) \rightarrow
\GradAlg^{H_B}(R)$ be the first projection.

\begin{definition}\label{unobstrdef}
  Let $R \rightarrow B$ be a graded quotient, let $M$ be a graded $B$-module
  and let $\varphi:M \rightarrow B$ be a homogeneous $B$-linear map.  Let
  $(T,{m}_T) \rightarrow (S, {m}_S)$ be a small Artin
  surjection (i.e. of local Artinian $k$-algebras with residue fields $k$ whose
  kernel ${\mathfrak a}$ satisfies ${\mathfrak a} \cdot {m}_T=0$). A
  graded deformation $B_S$ of $B$ to $S$ is a graded $S$-flat quotient of $R
  \otimes_k S$ satisfying $B_S \otimes_S k \simeq B$.
  
  \ (i) \ $(M,\varphi)$ is said to be unobstructed along any graded
  deformation of $B$ if for every small Artin $T \twoheadrightarrow S$ and for
  {\rm every} graded deformation $M_S \rightarrow B_S$ (of $S$-flat
  $B_S$-modules) of $M \stackrel{\varphi}{\rightarrow}B$ to $S$, there exists,
  for {\rm every} graded deformation $B_T$ of $B_S$ to $T$, a graded
  deformation $M_T \rightarrow B_T$ over $M_S \rightarrow B_S$ (i.e. a
  morphism of T-flat $B_T$-modules reducing to $M_S \rightarrow B_S$ via $(-)
  \otimes_T S$).
  
  \ (ii) \ One correspondingly defines $(M,B)$ to be unobstructed along any
  graded deformation of $B$ by forgetting $\varphi$, i.e. by considering pairs
  $(M_S, B_S)$ instead of $M_S \rightarrow B_S$.  The unobstructedness of
  $(M,B)$ and $(M(-s),B)$, $s$ an integer, are obviously equivalent. Moreover
  a surjection of graded quotients of $R$, $B \rightarrow A$, is unobstructed
  if every graded deformation $B_S \rightarrow A_S$ deforms further to $T$ for
  every small Artin $T \twoheadrightarrow S$.  Similarly a quotient $A$ of $R$
  is unobstructed if every graded deformation $A_S$ deforms further to $T$.
\end{definition}

Note that we have defined unobstructedness of graded objects by considering
graded deformations only.  Moreover note that, in deforming quotients $B$ and
$A$ of $R$ to $S$ the corresponding deformation of $R$ is always the trivial
one (i.e. $R \otimes_k S$) while the deformations $B_S$ and $A_S$ need not be
trivial.
 
\begin{remark}\label{unobstrrem}
  Let $T \rightarrow S$ be a small Artin surjection with kernel ${\mathfrak
    a}$. The obstruction of deforming a graded $B_S$-module $M_S$ (resp. a
  graded morphism $M_S \rightarrow N_S$) to $T$ sits in $_0\!\Ext_B^2(M,M)
  \otimes {\mathfrak a}$ (resp. $_0\!\Ext_B^1(M,N) \otimes {\mathfrak a}$) and
  $_0\!\Ext_B^1(M,M) \otimes {\mathfrak a} $ (resp. $_0\!\Hom_B(M,N) \otimes
  {\mathfrak a}$) corresponds to their set of graded deformations
  respectively. For a generically complete intersection, $R \twoheadrightarrow
  B$, its obstructions sit in $_0\!\Ext_B^1(I_B/I_B^2,B) \otimes {\mathfrak
    a}$. Note that, by definition, the vanishing of the obstruction is
  equivalent to the existence of corresponding desired deformation.
\end{remark}

Since the whole deformation theory fits functorially together (by e.g. Laudal's
work on deformations of categories in \cite{L}), we get in particular

\begin{proposition}\label{unobstrp1}
  Let $M$ be a finitely generated graded $B$-module.
  
  (i) If $\varphi:M \rightarrow B$ is a $B$-module homomorphism and $\ 
  _0\!\Ext_B^1(M,B)=0$, then $(M,\varphi)$ is unobstructed along any graded
  deformation of $B$ if and only if $(M,B)$ is unobstructed along any graded
  deformation of $B$.
 
  (ii) If $\ _0\!\Ext_B^2(M,M)= 0$, then $(M,B)$ is unobstructed along any
  graded deformation of $B$.
    
  (iii) If $\ _0\!\Ext_B^1(M,M)=0$ and if for every local Artinian $k$-algebra
  $T$ with residue field $k$ and for {\rm every} graded deformation $B_T$ of
  $B$ to $T$, there exists a graded deformation $M_T$ of $M$ to $B_T$, then
  $(M,B)$ is unobstructed along any graded deformation of $B$.
\end{proposition}

\begin{proof}
  (i) If $T \rightarrow S$ is a small Artin surjection
  (Definition~\ref{unobstrdef}), it follows from the exact sequence
  $$
  0 \rightarrow B \otimes_k {\mathfrak a} \rightarrow B_{T} \rightarrow
  B_{S} \rightarrow 0 $$
  that $_0\!\Hom_{B_T}(M_T,B_T)\rightarrow \ 
  _0\!\Hom_{B_S}(M_S,B_S)$ is surjective and hence that
  $_0\!\Hom_{B_S}(M_S,B_S)\rightarrow \ _0\!\Hom_B(M,B)$ is surjective by
  induction. This implies that $\varphi$ lifts to deformations of $M$ and $B$,
  and we easily get $(i)$. (Remark~\ref{unobstrrem} may provide a quicker
  proof.)
  
  Finally $(ii)$ and $(iii)$ follow from Remark~\ref{unobstrrem}. Here we
  leave a few details to the reader, remarking only that for $(iii)$ the
  assumption $\ _0\!\Ext_B^1(M,M)=0$ implies that an isomorphism $M_S
  \rightarrow M'_S$ of deformations of $M$ to $B_S$ lifts further to an
  isomorphism of {\rm given} deformations of $M_S$ and $M'_S$ to $B_T$, i.e.
  all deformations of $M$ to $B_S$ are isomorphic.
\end{proof}

\begin{example} Let $B = R/I_B$ be a graded $n$-dimensional licci quotient of
  $R$ of codimension $c$ and suppose $\depth_{I(Z)}B \geq 2$ where $Y-Z$ is an
  l.c.i. in $ \Proj(R)$.
  
  i) Then $(K_B,B)$ is unobstructed along any graded deformation of $B$ by
  Proposition~\ref{unobstrp1}$(ii)$. Indeed Proposition~\ref{unobstrp1}$(iii)$
  also applies because $\ K_{B_T} := \!  \Ext_{B_T}^c(B_T, R_T(-n-c))$ is a
  graded deformation of $K_B$ to $B_T$ by \cite{JS}, Prop.\! A1.
  
  ii) If $c=2$, then $(I_B/I_B^2,B)$ is unobstructed along any graded
  deformation of $B$ by Proposition~\ref{unobstrp1}$(ii)$ and
  Proposition~\ref{introprop}. Moreover if $\depth_{I(Z)}B \geq 4$ we will see
  in Section 3 that $(N_B,B)$ is unobstructed along any graded deformation of
  $B$ because the assumptions of Proposition~\ref{unobstrp1}$(iii)$ holds
  (here Proposition~\ref{unobstrp1}$(ii)$ may not apply, cf.
  Remark~\ref{remcorNB}).
\end{example}

We will need the following result on how deformations of $R \rightarrow A$ is
related to deformations of $R \rightarrow B$ provided $ A \simeq B/I_{A/B}$ is
a graded quotient of $B$ (\cite{K03}, Thm. 5 and Rem.\! 6). To state it, let
$U \subset \GradAlg^{H_A}(R)$ be a sufficiently small open subset containing
$(A)$ and let $p:\GradAlg(H_B,H_A) \rightarrow \GradAlg^{H_A}(R)$ be the
second projection, i.e. given by $p((B' \rightarrow A'))= (A')$. Then the
$k$-points of the subset $p(p^{-1}(U))$ of $U$ correspond to quotients $R
\rightarrow A'$ with Hilbert function $H_A$ for which there exist some
factorization $B' \rightarrow A'$ such that $B'$ has Hilbert function $H_B$.
We will call it a {\it stratum of $H_B$-factorizations around $(A)$}, and
$\dim U - \dim p(p^{-1}(U))$ the {\it codimension of the $H_B$-stratum of
  $A$}. At least if $U$ is smooth, it is the ordinary codimension of
$p(p^{-1}(U))$ in $U$. $A$ is called {\it $H_B$-generic} if there is an {\em
  open} subset $U_A$ of $\GradAlg^{H_A}(R)$ such that $(A) \in U_A \subset
p(p^{-1}(U))$. The {\it $H_B$-stratum of $A$ at $(B \rightarrow A)$, (resp.
  its codimension)}, is defined to be $p(U')$ (resp. $\dim U - \dim p(U')$)
where $U' \subset p^{-1}(U)$ is an open subset in the union of the
irreducible components of $ p^{-1}(U)$ which contain $(B \rightarrow A)$, 
and $(B \rightarrow A) \in U'$. If $(I_{A/B})_v =0$ for $v \leq $ the largest
degree of the minimal generators of $I_B$, then $p \arrowvert_{U'}$ will be
unramified and universally injective and the two concepts of codimension above
coincide (\cite{K03}, Lem.\! 7).

\begin{theorem} \label{varunobstr} Let $R$ be a graded polynomial $k$-algebra,
  let $B = R/I_B \twoheadrightarrow A \simeq B/I_{A/B}$ be a graded
  morphism of quotients of $R$ and suppose $_0\!\HH^2(B,A,A) =0$, $\
  \depth_{\mathfrak{m}}A \geq
  \min(1,\dim A)$ and $\ \depth_{\mathfrak{m}}B \geq \min(1,\dim B)$. \\[2mm]
  \indent {\rm A)} \ If $ \ \ \ \ _0\!\Ext_B^1(I_B/I_B^2,A) =0$ and
  $(I_B)_{\wp}$ is
  syzygetic  for any graded $\wp$ of $Ass(A)$, \\[2mm]
  then $A$ is unobstructed as a graded $R$-algebra (indeed
  $_0\!\HH^2(R,A,A)=0$), and $$ \ \dim_{(A)} \GradAlg^{H_A}(R)= \
  _0\!\hom_R(I_B,B) +\ _0\!\hom_B(I_{A/B},A)- \ _0\!\hom_R(I_B,I_{A/B}) +$$
  $$ _0\!\ext_B^1(I_B/I_B^2,I_{A/B})-\ _0\!\ext_B^1(I_B/I_B^2,B).$$ Moreover
  let $B$ be unobstructed as a graded $R$-algebra, let $k$ be of
  characteristic zero and suppose $(B \rightarrow A)$ is general with respect
  to $\ _{0}\!\hom_R(I_B,I_{A/B}) $. Then the codimension of the $H_B$-stratum
  of $A$ at $(B \rightarrow A)$ is
$$
\ _0\!\ext_B^1(I_B/I_B^2,I_{A/B})-\ _0\!\ext_B^1(I_B/I_B^2,B).$$

B)  If
\hspace{0.1 cm} $_0\!\Ext_B^1(I_B/I_B^2,I_{A/B}) = 0$ and $(I_B)_{\wp}$ is
syzygetic for any graded prime $\wp$ of $Ass(I_{A/B})$, \\[2mm]
then $A$ is $H_B$-generic. % for $p(\GradAlg(H_B,H_A))$. 
Moreover $A$ is unobstructed as a graded $R$-algebra if and only if $B$ is
unobstructed as a graded $R$-algebra. Indeed \\[2mm]
 \centerline{
$\ _0\!\hom_R(I_A,A) - \dim_{(A)} \GradAlg^{H_A}(R) = \ _0\!\hom_R(I_B,B)-
\dim_{(B)} \GradAlg^{H_B}(R) $, \  and}
   $$
   \ \dim_{(A)} \GradAlg^{H_A}(R)=
 \dim_{(B)} \GradAlg^{H_B}(R) %-\ _0ext_B^1(I_{A/B},A) 
+ \ _0\!\hom_B(I_{A/B},A)-\ _0\!\hom_R(I_B,I_{A/B}) \ .$$
\end{theorem} 

\begin{remark}\label{mre} (i) Theorem~\ref{varunobstr} follows directly from
  Thm.\! 5 and Rem.\! 6 of \cite{K03}. Moreover one may replace the assumption
  ``$\ _0\!\Ext_B^1(I_B/I_B^2,I_{A/B}) = 0$ and $(I_B)_{\wp}$ syzygetic for
  any graded prime $\wp$ of $Ass(I_{A/B})$'' of Theorem~\ref{varunobstr}B) by
  ``$\ _0\!\Ext_R^1(I_B,I_{A/B}) = 0$'' and conclude exactly as in
  Theorem~\ref{varunobstr}B) because, what's needed to prove part B) is $\
  _0\!\HH^2(R,B,I_{A/B})=0$ (cf. the proof of Thm.\! 5 of \cite{K03}). Since
  we in general have an injection $\ _0\!\HH^2(R,B,I_{A/B}) \hookrightarrow \
  _0\!\Ext_R^1(I_B,I_{A/B}),$ we get the claim.

  (ii) To find the {\rm codimension of $H_B$-stratum of $A$ at $(B \rightarrow
    A)$}, the proof of \cite{K03}, Thm.\! 5 uses generic smoothness. Indeed it
  computes the dimension of the image $p(U')$, see the text before Theorem 15
  above, in terms of the dimension of $ \GradAlg (H_B,H_A)$ at $(B \rightarrow
  A)$, and the dimension, $\ _{0}\!\hom_R(I_B,I_{A/B})$, of a {\rm general}
  fiber. Unfortunately the word ``general'', i.e. the assumption ``$(B
  \rightarrow A)$ is general with respect to $\ _{0}\!\hom_R(I_B,I_{A/B}) $''
  is missing in \cite{K03}, Thm.\! 5 as well as in Theorem 15 of the published
  version, [{\em Collect. Math.} 58, 2 (2007), 199-238], of this paper. Indeed
  it is easily seen from the proof in \cite{K03} (cf.\! \cite{K03}, Prop\!
  4(i)) that we may suppose $U'$ above is a smooth irreducible scheme
 % containing $(B \rightarrow A)$ 
such that the restriction of the {\rm first}
  projection $q:\GradAlg(H_B,H_A) \rightarrow \GradAlg^{H_B}(R)$, $q((B'
  \rightarrow A'))= (B')$, to $U'$ is a smooth morphism. Hence if we let ``$(B
  \rightarrow A)$ is general with respect to $\ _{0}\!\hom_R(I_B,I_{A/B}) $''
  mean that, for a given $(B \rightarrow A)$, $\ _{0}\!\hom_R(I_B,I_{A/B}) $
  obtains its least possible value in $U'$, or equivalently that $ \
  _{0}\!\hom_R(I_B,I_{A/B}) = \ _{0}\!\hom_R(I_{B'},I_{A'/B'}) $ for every
  $(B' \rightarrow A')$ in an {\it open neighbourhood} of $(B \rightarrow A)$
  in $U'$, the proof is complete. Note that this assumption of generality is
  always satisfied if $\ _{0}\!\hom_R(I_B,I_{A/B})=0$, in which case we do not
  need to suppose ``$char(k) = 0$'' either because then $p \arrowvert_{U'}$ is
  unramified. In particular in \cite{K03}; Theorem 1, Proposition 13, Theorem
  16 (and hence in Theorem 23 of [{\em Collect. Math.} 58, 2 (2007),
  199-238]), for the {\bf codimension statement} in the {\rm A)}-part of the
  results we need to assume that ``$(B)$ is general with respect to $\
  _{-s}\!\hom_R(I_B,K_{B})$'', i.e. that $\ _{-s}\!\hom_R(I_B,K_{B})$ obtains
  its least possible value in the open subset $q(U')$ of $\GradAlg^{H_B}(R)$,
  since we in this application use \cite{K03}, Thm.\! 5 with
  $I_{A/B}:=K_B(-s)$. Several results of the published version of this paper
  in {\em Collect. Math.} 58, 2 (2007), 199-238, as well as the first version
  of this paper on the arXiv, i.e. Thm.\! 1, Cor.\! 18, Cor.\! 20, Thm.\! 25,
  Thm.\! 30, Cor.\! 37 and Cor.\! 44, use Theorem 15 and lack the generality
  assumption on $(B \rightarrow A)$ above which is needed for the validity of
  the {\bf codimension statements} of the strata. In this version all these
  results are corrected, and in the examples where the codimensions of the
  strata are stated and computed, one may check that $\
  _{0}\!\hom_R(I_B,I_{A/B})=0$, so no corrections are needed there.
\end{remark}

The following lemma is a slight improvement of \cite{K03}, Lem.\!\ 15 which we
will need to transform Theorem~\ref{varunobstr} into a result on a Gorenstein
quotient of a CM algebra.

\begin{lemma} \label{lemvarunobstr}
  Let $B$ be a graded CM quotient of $R$ and let $A \simeq
  B/I_{A/B}$ be Gorenstein such that the canonical module $K_A \simeq A(j)$
  and $\dim B -\dim A = r$. We have
  
  (i) \ If $\ r>0$ (\! resp. $ r = 0 $\! ), then $\ \ _0\!\HH^2(B,A,A)=0 \ $
  if and only if $ \ _{-j}\Ext_B^{r}(S_2(I_{A/B}),K_B)=0$ (resp. iff the
  natural map $ \! _{-j}\Hom_B(I_{A/B},K_B) \rightarrow \! 
  _{-j}\Hom_B(S_2(I_{A/B}),K_B)$ is surjective).
  
  (ii) \ $ _0\!\hom_B(I_{A/B},A) = \ _{-j}\ext_B^r(B/I_{A/B}^2,K_B) -1 $ for
  any $r \geq 0$. Moreover if $r \geq 2$, then $ _0\!\hom_B(I_{A/B},A) = \ 
  _{-j}\ext_B^{r-1}(S_2(I_{A/B}),K_B) -1 \ . $

  (iii) \ Let $\dim A = 0$. If the surjection $S_2(I_{A/B})_j \rightarrow
  (I_{A/B}^2)_j$ is an isomorphism, then  $_0\!\HH^2(B,A,A)=0$
  and  $ _0\!\hom_B(I_{A/B},A) = \dim B_{j} -\dim S_2(I_{A/B})_j -1$. In
  particular if  $S_2(I_{A/B})_j=0$, then  $_0\!\HH^2(B,A,A)=0$
  and  $ _0\!\hom_B(I_{A/B},A) = \dim B_{j} -1$.      
 \end{lemma}

\begin{proof}
  (i) \ By \cite{K03}, Lemma 15 (whose proof is quite close to the proof
  below), we get the vanishing of $_0\!\HH^2(B,A,A)$ from $ \!
  _0\!\Ext_B^{r}(S_2(I_{A/B}),K_B(-j))=0$, and conversely if $r>0$. By paying
  a little extra attention to the case $r=0$ we get $(i)$. 
  
  (ii) \ We have by $A \simeq K_A(-j)$ and Gorenstein
  duality (applied to both $A$ and $B$) that
  $$
  _0\!\Hom_A(I_{A/B}/I_{A/B}^2,A) \simeq \ _j\!\HH^{\dim
    A}_\mathfrak{m}(I_{A/B}/I_{A/B}^2)^{\vee} \simeq \ _{-j}\!\Ext_B^{r}(
  I_{A/B}/I_{A/B}^2,K_B).$$
  Moreover since $ \ \Ext_B^{r}(A,K_B) \simeq
  \Hom_A(A,K_A) \simeq A(j)$ and $ \Ext_B^{r-i}(D,K_B)=0$ for $i>0$ provided
  $D$ is an $A$-module, we get $ _0\!\hom_B(I_{A/B},A) = \! 
  _{-j}\ext_B^r(B/I_{A/B}^2,K_B) -1 $ by the long exact sequence of
  $\Hom(-,K_B)$ applied to $ 0 \rightarrow I_{A/B}/I_{A/B}^2 \rightarrow \  
  B/I_{A/B}^2 \rightarrow A \rightarrow 0$.
  Moreover applying  $\Hom(-,K_B)$ onto $ 0 \rightarrow I_{A/B}^2 \rightarrow
  \   B \rightarrow  B/I_{A/B}^2  \rightarrow 0$, we get $ \
  _{-j}\Ext_B^r(B/I_{A/B}^2,K_B) \simeq \ _{-j}\Ext_B^{r-1}(I_{A/B}^2,K_B)$
  provided $r>1$.  Since the
  kernel of the surjection $S_2(I_{A/B}) \rightarrow I_{A/B}^2$ is an A-module
  (namely $\HH_2(B,A,A)$, cf.  \cite{VAS}, Section 2.1), we get $
  \!\Ext_B^{r-1}(I_{A/B}^2,K_B) \simeq \Ext_B^{r-1}(S_2(I_{A/B}),K_B)$ and we
  have (ii).
  
  (iii) \ If $\dim A = 0$ and the surjection $S_2(I_{A/B})_{j} \rightarrow
  (I_{A/B}^2)_{j}$ is an isomorphism, we get $\dim B = r$ and $
  _{j}\!\HH_2(B,A,A)=0$ and hence $$ _{-j}\Ext_B^{r}(S_2(I_{A/B}),K_B)^{\vee}
  \simeq \! _{j}\!\HH_{\mathfrak{m}}^{0}(S_2(I_{A/B})) \simeq
  \!_{j}\!\HH_{\mathfrak{m}}^{0}(I_{A/B}^2) $$ by Gorenstein duality. If $r>0$
  (resp. $r=0$), the last group vanishes (resp. injects into
  $\!_{j}\HH_{\mathfrak{m}}^{0}(I_{A/B})$). Hence we get $_0\!\HH^2(B,A,A)=0$
  by $(i)$ and since $ _{-j}\Ext_B^r(B/I_{A/B}^2,K_B)^{\vee} \simeq \!
  _{j}\!\HH_{\mathfrak{m}}^{0}(B/I_{A/B}^2) = (B/I_{A/B}^2)_j$, we conclude by
  $(ii)$ and by assumption.
\end{proof}

\begin{corollary} \label{corvarunobstr}
  Let $B = R/I_B$ be Cohen-Macaulay and let $A \simeq B/I_{A/B}$ be a
  graded Gorenstein quotient such that $K_A \simeq A(j)$ and $\dim B -\dim A
  = r$. If 
  \begin{itemize}
  \item [(i)] $_0\!\Ext_B^1(I_B/I_B^2,B) =0$ and $(I_B)_{\wp}$ is syzygetic
    for any graded $\wp$ of $Ass(A) \cup Ass(B)$,
   \item [(ii)]  $ _{-j}\!\Ext_B^{r}(S_2(I_{A/B}),K_B)=0$, and
   \item [(iii)] $_0\!\Ext_B^2(I_B/I_B^2,I_{A/B})=0$, 
  \end{itemize} 
  then $A$ is unobstructed as a graded $R$-algebra (indeed
  $_0\!\HH^2(R,A,A)=0$).  Moreover, if $r \geq 2$, then
$$ \ \dim_{(A)} \GradAlg^{H_A}(R)=\ _0\!\hom_R(I_B,B) +  \!
_{-j}\ext_B^{r-1}(S_2(I_{A/B}),K_B) -1 -\delta(I_{A/B})_0$$ where
$\delta(I_{A/B})_0 := \ _0\!\hom_R(I_B,I_{A/B})- \
_0\!\ext_B^1(I_B/I_B^2,I_{A/B}).$ Furthermore if $char(k)=0$ and $(B
\rightarrow A)$ is general with respect to $\ _{0}\!\hom_R(I_B,I_{A/B}) $,
then the codimension of the $H_B$-stratum of $A$ at $(B \rightarrow A)$ is $\
_0\!\ext_B^1(I_B/I_B^2,I_{A/B}).$
\end{corollary}

\begin{proof}
  By (i) we have that $B$ is unobstructed. Moreover (i), (iii) and the long
  exact sequence of $\Hom_B(I_B/I_B^2,-)$ applied to $ 0 \rightarrow I_{A/B}
  \rightarrow \ B \rightarrow A \rightarrow 0$ show the first line of
  assumptions of Theorem~\ref{varunobstr}A). Hence we conclude by
  Lemma~\ref{lemvarunobstr}.
\end{proof} 

\begin{remark} \label{remvarunobstr} 
We get   $_v\!\HH^2(R,A,A)=0$ by twisting the three $\Ext_B^i$-vanishing
assumptions of  Corollary~\ref{corvarunobstr} by $v$ because
Theorem~\ref{varunobstr}A) admits such a generalization (cf. \cite{K03},
Rem.\! 6$(a)$).
\end{remark}

\begin{corollary} \label{corvarunobstr2}
  Let $B = R/I_B$ be Cohen-Macaulay and let $A \simeq B/I_{A/B}$ be a
  graded Artinian Gorenstein quotient such that $K_A \simeq A(j)$ and $\dim B
  -\dim A = r$. If
  \begin{itemize}
  \item [(i)] $_0\!\Ext_B^1(I_B/I_B^2,B) =0$ and $(I_B)_{\wp}$ is syzygetic
    for any graded $\wp$ of $Ass(A) \cup Ass(B)$, 
   \item [(ii)]   $\dim S_2(I_{A/B})_j = \dim (I_{A/B}^2)_j \ $, and 
 % (e.g.  $(I_{A/B})_{j/2}=0$)
    \item [(iii)] $_0\!\Ext_B^2(I_B/I_B^2,I_{A/B})=0$, 
  \end{itemize} 
  then $A$ is unobstructed as a graded $R$-algebra (indeed
  $_0\!\HH^2(R,A,A)=0$), and
  $$
  \ \dim_{(A)} \PGor(H_A)=\ _0\!\hom_R(I_B,B) + \dim B_{j} -\dim
  S_2(I_{A/B})_j -1 -\delta(I_{A/B})_0$$ where $\delta(I_{A/B})_0 := \
  _0\!\hom_R(I_B,I_{A/B})- \ _0\!\ext_B^1(I_B/I_B^2,I_{A/B}).$ Furthermore if
  $char(k)=0$ and $(B \rightarrow A)$ is general with respect to $\
  _{0}\!\hom_R(I_B,I_{A/B}) $, then the codimension of the $H_B$-stratum of
  $A$ at $(B \rightarrow A)$ is $\ _0\!\ext_B^1(I_B/I_B^2,I_{A/B}).$
\end{corollary}

The proof is similar to the preceding Corollary. Note that if $B$ is a c.i.
and $(I_{A/B})_{v}=0$ for $v \leq j/2$, then all assumptions of
Corollary~\ref{corvarunobstr2} are satisfied. Hence \cite{IK}, Thm.\!\ 4.17 is
generalized by Corollary~\ref{corvarunobstr2} and we get in addition that
$\PGor(H_A)$ is generically smooth along the component of \cite{IK}, Thm.\!\
4.17.

\begin{remark} \label{remcorvarunobstr}
  If the Gorenstein algebra $A$ is given by \eqref{LESTHM}, then condition
  $(iii)$ of Corollary~\ref{corvarunobstr} and \ref{corvarunobstr2} is
  satisfied provided
\begin{equation} \label{varcorvar}
 _0\!\Ext_B^{i+1}(I_B/I_B^2,M_i(-is))= 0 \ \ {\rm for} \ 1 \leq i
  \leq r
\end{equation}
 letting $M_1=M$. Indeed splitting the exact sequence \eqref{LESTHM} into short
  exact sequences, we get 
  $$0 \rightarrow Z_i \rightarrow M_i(-is) \rightarrow Z_{i-1} \rightarrow 0$$
  where $Z_i = \ker (M_i(-is) \rightarrow M_{i-1}((1-i)s))$ for $i \geq 1$ and
  $Z_0=I_{A/B}$. Successively applying $_0\!\Ext^{i+1}_B(I_B/I_B^2,-)$ onto
  this sequences and using \eqref{varcorvar} we get
  $_0\!\Ext_B^{i+1}(I_B/I_B^2,Z_{i-1})= 0$ from
  $_0\!\Ext_B^{i+1}(I_B/I_B^2,M_i(-is))= 0$ and
  $_0\!\Ext_B^{i+2}(I_B/I_B^2,Z_i)= 0$, i.e. we get $(iii)$ because $Z_{r-1}=
  M_r(-rs)$.
\end{remark}

It is desirable to weaken the general assumption $_0\!\HH^2(B,A,A) =0$ of
Theorem~\ref{varunobstr}. In the appendix we prove a nice result,
Theorem~\ref{unobstr}, on deformations of the degeneracy locus of a regular
section of a maximal CM sheaf, which leads to a variation of
Theorem~\ref{varunobstr}B) (Theorem~\ref{gorgenth} below) in which the
mentioned assumption is replaced by assumptions closer to set-up of
Theorem~\ref{mainth}. The main idea of the proof of Theorem~\ref{unobstr} and
hence of Theorem~\ref{gorgenth} is to use the complete intersection property
(in large enough open subset) of a regular section to control all
deformations of the degeneracy locus. Since the proofs are probably of most
interest only for specialists in deformation theories, we delay them to the
appendix. For the case $r=1$ and $\depth_{I(Z)}B =1$, which requires special
attention, we refer to \cite{K03}.

\begin{theorem} \label{gorgenth} Let $R$ be a finitely generated polynomial
  $k$-algebra, let $B = R/I_B$ be a graded Cohen-Macaulay quotient of $R$
  of codimension $c$, and let $M$ be a finitely generated graded maximal
  Cohen-Macaulay $B$-module. Let $Y=\Proj(B)$, let $Z$ be a closed scheme such
  that $\dim B - \dim {B}/{I}(Z) \geq {\rm max}(r,2)$ and let $U=Y-Z$. Let $A$
  be defined by a regular section $\sigma$ of $\widetilde{M}^*(s)$ on $U$,
  i.e. given by \eqref{LESTHM}. Let ${ M}_i = \HH^0_*(U,\wedge^i \widetilde {
    M})$ for $i\geq 0$, and suppose ${\widetilde { M}}\arrowvert_U$ is locally
  free (of rank $r >0$) and
  $$\wedge^r{\widetilde { M}}\arrowvert_U \simeq \widetilde {{ K}}_{
    B}(t)\arrowvert_U \ {\rm (for \ some \ integer \ } t).$$
  Moreover suppose
  ${M}_i$ are maximal Cohen-Macaulay ${B}$-modules for $2 \leq i \leq r/2$.
   If
\begin{itemize}
  \item [(i)] $_0\!\Ext^i_B(M,M_i(-(i-1)s))= 0$ \ \ {\rm for} \ $2 \leq i \leq
    r-1$, and
\item [(ii)] $(M(-s),\sigma)$ is unobstructed along any graded deformation
  of $B$, and 
\item [(iii)] {\rm either} $_0\!\Ext^i_B(I_B/I_B^2,M_i(-is))= 0$ \ {\rm for}
  $1 \leq i \leq r$ and $(I_B)_{\wp}$ is syzygetic for any $\wp$ of $Ass(B)$,
  {\rm or} $_0\!\Ext^i_R(I_B,M_i(-is))= 0$ {\rm for} $1 \leq i \leq r$,
  \end{itemize}
  then $\ \ \ _0\!\hom_R(I_B,B)-\ \dim_{(B)}\GradAlg^{H_B}(R) =\ _0\!
  \hom_R(I_A,A)- \dim_{(A)} \GradAlg^{H_A}(R)  $,  \ $A$ is $H_B$-generic
  and Gorenstein, and
 $$
 \dim_{(A)} \GradAlg^{H_A}(R)= \dim_{(B)}\GradAlg^{H_B}(R)+\ 
 _0\! \hom_B(I_{A/B},A)-\ _0\!\hom_R(I_B,I_{A/B}).$$
 Moreover $A$ is unobstructed
 as a graded $R$-algebra if and only if $B$ is unobstructed as a graded
 $R$-algebra.
\end{theorem}

Fortunately, by results of the next sections we will see that many of the
$\Ext$-groups of Corollary~\ref{corvarunobstr} and Theorem~\ref{gorgenth}
vanish, and that the dimension formula turns out to be computable, provided we
apply it to a licci quotient $B$ of $R$ of small enough codimension.

\section{Families of Gorenstein quotients of low codimension on Cohen-Macaulay
  algebras}

An important issue of this paper is to study families of graded Gorenstein
quotients $A$ obtained by taking regular sections of the dual of a maximal
Cohen-Macaulay sheaf of rank r, i.e. quotients given by Theorem~\ref{mainth}.
In the following we will see how Corollary~\ref{corvarunobstr} and
Theorem~\ref{gorgenth} enable us to treat the low rank cases $r \leq 3$ on a
licci scheme satisfactorially. The different values of $r$, $1 \leq r \leq 3$,
require special attention and lead to the three theorems of this section.

The case $r=1$ was considered in \cite{K03}. There we proved the following
result in which $N_B:=\Hom_B(I_B/I_B^2,B)$ and $K_B^*:=\Hom_B(K_B,B) \simeq
\Hom_B(S_2(K_B),K_B)$. The conclusions below about when ``the codimension of
the stratum is zero'' overlap results proved by others (\cite{B992}, Thm.\!
3.2, cf. \cite{Hacol}, Thm.\! 3.5). Note that we say ``the stratum of
quotients given by \eqref{K(-s)} around or at $(A)$'' for ``the $H_B$-stratum
of $A$ at $(B \rightarrow A)$'' provided any $(B' \rightarrow A')$ in a small
enough {\it open} neighbourhood of $(B \rightarrow A)$ in $\GradAlg(H_B,H_A)$
is given by \eqref{K(-s)}.
\begin{theorem} \label{Kmainth}
  Let $B = R/I_B$ be a generically Gorenstein, graded Cohen-Macaulay quotient
  of a polynomial ring $R$, and let $A$ be a graded codimension one quotient
  of $B$, given by an exact sequence
 \begin{equation} \label{K(-s)}
 0 \rightarrow K_B(-s) \rightarrow B \rightarrow A \rightarrow 0 \ , \ \ \
 \  s \ \  {\rm  an \ \ integer} \ .
\end{equation}
\indent {\rm A)} If $B$ is licci, then $A$ is unobstructed as a graded
$R$-algebra (indeed $\ \HH^2(R,A,A)=0$), $A$ is Gorenstein and,
  $$
  \ \dim_{(A)} \GradAlg^{H_A}(R)= \ \dim_{(B)} \GradAlg^{H_B}(R) + \ \dim
  (K_B^*)_s -1 -\delta(K_B)_{-s}$$ where $ \delta(K_B)_{-s} =\
  _{-s}\hom_B(I_B/I_B^2,K_{B})-\ _{-s}\ext_B^1(I_B/I_B^2,K_{B})$. Moreover if
  char$(k)=0$ and $(B)$ is general with respect to $\
  _{-s}\!\hom_B(I_B/I_B^2,K_{B}) $, then the codimension
  of %the $H_B$-stratum of $A$ at
 %  $(B \rightarrow A)$, or equivalently 
  the stratum of quotients given by \eqref{K(-s)} around $(A)$ is $ \
  _{-s}\ext_B^1(I_B/I_B^2,K_{B})$. 

 \indent {\rm B)} If $s >> 0$ and $\
  \Proj(B)$ is locally licci, then $A$ is Gorenstein and the codimension of
  the stratum of quotients given by \eqref{K(-s)} around $(A)$ is zero (so $A$
  is $H_B$-generic), and
  $$
  \ \dim_{(A)} \GradAlg^{H_A}(R)= \ \dim_{(B)} \GradAlg^{H_B}(R) + \ \dim
  (K_B^*)_s -1.$$ Moreover $A$ is unobstructed as a graded $R$-algebra iff $B$
  is unobstructed as a graded $R$-algebra.
\end{theorem}

\begin{example} \label{exKmainth} (Arithmetically Gorenstein curves $\Proj(A)$
  in $ \pp^5$, obtained by \eqref{K(-s)}.) Let $R$ be a polynomial ring in 6
  variables, and let $Y= \Proj(B)$ be a generically Gorenstein ACM-surface
  with resolution
  $$0 \rightarrow R(-7)^3 \rightarrow R(-6)^3\oplus R(-5)^3 \rightarrow
  R(-3)^4 \rightarrow I_B \rightarrow 0$$ and Hilbert function $H_B(v) = 19
  \binom{v+1}{2} -26v+16$, $v \geq 2$. We may obtain $B$ by taking a c.i. $B'$
  of type $(2,2,2)$ and linking it to $B$ via some c.i. of type $(3,3,3)$.
  Since the linkage result, Prop.\!\ 33 of \cite{K03}, connects all invariants
  of $B$, appearing in Theorem~\ref{Kmainth}, to the corresponding invariants
  of $B'$ and since $N_{B'} \simeq B'(2)^{ \oplus 3}$, $\
  \Hom_{B'}(I_{B'}/I_{B'}^2,K_{B'}) \simeq K_{B'}(2)^{ \oplus 3}$ and $ K_{B'}
  \simeq B'$, we get by \cite{K03} Prop.\!\ 33 that $\ \dim (N_B)_0 = \ \dim
  (N_{B'})_0 + 3H_B(3)-3H_{B'}(3) =96 \ , $ and that
  $$\ \dim (K_B^*)_v = \dim B_{v-3} + 3 \dim (I_{B/D})_{v}-
  \delta(K_{B'})_{v-6} =67 + (19v^2-99v)/2 \ $$
  because $ \delta(K_{B'})_{v}=
  \ _{v}\hom_{B'}(I_{B'}/I_{B'}^2,K_{B'})$.  Similarly, $\ \delta(K_B)_{v-6} =
  \dim B'_{v-3} + 3 \dim (I_{B'/D})_{v}- \dim (K_{B'}^*)_v$. Hence we get $
  \delta(K_B)_{v}=0$ for $v \leq -7$ and $ (\delta(K_B)_{-5},
  \delta(K_B)_{-6}) = (-6,-1)$. Let $A$ be defined by
  \eqref{K(-s)} for $s \geq 5$, so $X = \Proj(A)$ is an AG curve in $ \pp^5$
  with Hilbert polynomial $p_X(v)= (19s-33)v - s (19s-33)/2$.  Then $A$ is
  unobstructed by Theorem~\ref{Kmainth} and (cf. \eqref{Grad}),
  $$
  \ \dim_{(A)} \GradAlg^{H_A}(R)= \ \dim_{(X)} \Hilb^{p_X}(\pp^{5}) =
  162+(19s^2-99s)/2- \delta(K_B)_{-s}\ , $$ where $
  _{-s}\ext_B^1(I_B/I_B^2,K_{B})= - \delta(K_B)_{-s}$ represents the
  codimension of quotients given by \eqref{K(-s)} at $(A)$.
\end{example}

In the following we concentrate on the case $r>1$, notably $r=2$. If $r=2$,
then Theorem~\ref{mainth} provides us with an exact sequence
\begin{equation} \label{M(-s)} 0 \rightarrow K_B(t-2s)\rightarrow M(-s)
  \stackrel{\sigma}{\rightarrow} B \rightarrow A \rightarrow 0 \ ,
\end{equation}
and we have $M \simeq \Hom_B(M,K_B)(t)$.  We define $\gamma(S_2M)_v$ by
$$\gamma(S_2M)_v =\ _{v}\!\hom_B(S_2(M),K_B(t))-\
_{v}\!\ext_B^1(S_2(M),K_B(t)) \  .
$$
Note that if $char(k) \neq 2$, then we get an exact sequence $\ 0 \rightarrow
\ \Hom_B(S_2(M),K_B(t)) \rightarrow \Hom_B(M,M) \rightarrow B \rightarrow 0$
from the split exact sequence $ \ 0 \rightarrow \ \wedge^2 M \rightarrow M
\otimes M \rightarrow S_2(M) \rightarrow 0$ and \eqref{twoo}. The split exact
sequence and \eqref{twoo} also lead to isomorphisms
\begin{equation} \label{M,M}
\ \Ext^i_B(M,M) \simeq \Ext_B^i(S_2(M),K_B(t)) \ \ {\rm for} \ i=1 \
({\rm resp.} \  i=2) 
\end{equation}
provided $\depth_{I(Z)}B \geq 3$ (resp. $\depth_{I(Z)}B \geq 4$). Indeed if
$\depth_{I(Z)}B \geq 3$ say, we have $$\Ext_B^{1}(M,M) \simeq
\Ext_{\sO_U}^1(\widetilde{M}\arrowvert_U , \widetilde {M}\arrowvert_U) \simeq
\Ext_{\sO_U}^1(\widetilde{M} \otimes \widetilde {M^*} \otimes \widetilde {
  K_B}\arrowvert_U, \widetilde {K_B} \arrowvert_U) \simeq
\Ext_B^1(S_2(M),K_B(t)) $$ where we have used $\widetilde{(M^* \otimes
  K_B)} \arrowvert_U \simeq \widetilde {M(-t)} \arrowvert_U$ which follows
from $M \simeq \Hom_B(M,K_B)(t)$. Hence if $\depth_{I(Z)}B \geq 3$ and
$char(k) \neq 2$, then $\gamma(S_2M)_v$ is given by $$\gamma(S_2M)_v =\
_{v}\!\hom_B(M,M)-\ _{v}\!\ext_B^1(M,M)-\ \dim B_v.$$ If $Q$ is a finitely
generated $B$-module, we define $\delta(Q)_v$ as previously by
$$
\delta(Q)_{v} =\ _{v}\!\hom_B(I_B/I_B^2,Q)-\ _{v}\!\ext_B^1(I_B/I_B^2,Q) \ 
.$$

If $A$ is a quotient of $B$ given by Theorem~\ref{mainth} and hence defined by
\eqref{M(-s)}, then a {\it stratum of quotients given by \eqref{M(-s)} around
  $(A)$} is of the form $p(W)$, $p:\GradAlg(H_B,H_A) \rightarrow
\GradAlg^{H_A}(R)$ the second projection, where $W$ is a maximal closed subset
of a small enough open affine neighbourhood of $(B \rightarrow A)$ in
$\GradAlg(H_B,H_A)$, consisting of quotients $(B' \rightarrow A')$ for which
there exists an extension $\xi': 0 \rightarrow K_{B'}(t-2s) \rightarrow M'(-s)
\rightarrow I_{A'/B'} \rightarrow 0$ which is induced by a regular section of
some $\widetilde {M'}\arrowvert_{U'}$ as in Theorem~\ref{mainth} ($\widetilde
{M'}\arrowvert_{U'}$ locally free and $\Proj(B')-U'$ of codimension at least
two). Moreover if $(B_W \rightarrow A_W)$ is the pullback to $W$ of the
universal element of $ \GradAlg(H_B,H_A)$ and $I_{A_W/B_W}= \ker(B_W
\rightarrow A_W)$, then there is an element $ \xi_W \in
\Ext_{B_W}^1(I_{A_W/B_W},K_{B_W}(t-2s))$ whose obvious pullbacks are $\xi'$
above and the extension given by \eqref{M(-s)}. Note that if $W$ is {\it open}
in $\GradAlg(H_B,H_A)$, then the stratum defined above coincides with the
$H_B$-stratum of $A$ at $(B \rightarrow A)$.

\begin{theorem} \label{Mmainth} Let $B = R/I_B$ be a graded licci
  quotient of $R$, %of codimension $c\geq 1$,
  let $M$ be a graded maximal Cohen-Macaulay B-module, and suppose
  $\widetilde{M}$ is locally free of rank $2$ in $U:=\Proj(B) -Z$, that $\dim
  { B}-\dim {B}/{I}(Z)\geq 2$ and $\wedge^2{\widetilde { M}}\arrowvert_U
  \simeq \widetilde {K_B}(t)\arrowvert_U$. Let $A$ be defined by a regular
  section $\sigma$ of $\widetilde{M}^*(s)$ on $U$, i.e. given by \eqref{M(-s)}
  and suppose $\ _s\!\Ext^1_B(M,B)=0$. \\[2mm]
  \indent {\rm A)} \ \ If $\ \ \ \ _t\!\Ext_B^2(S_2(M),K_B)=0 \ $ and $\
  _{-s}\!\Ext^2_B(I_B/I_B^2,M)= 0 \ $,   \\[2mm]
  then $A$ is unobstructed as a graded $R$-algebra (indeed
  $_0\!\HH^2(R,A,A)=0$), $A$ is Gorenstein, and {\small
$$  \ \dim_{(A)}
\GradAlg^{H_A}(R) \ = \ \dim (N_B)_0 + \ \dim (M^*)_s - 1 - \gamma(S_2M)_0 + \
\dim(K_B)_{t-2s} +\delta(K_B)_{t-2s} - \delta(M)_{-s}. $$ } Moreover if \
$char(k)=0$ and $(B \rightarrow A)$ is general with respect to $\
_{0}\!\hom_R(I_B,I_{A/B}) $, then the codimension of the $H_B$-stratum of $A$
at $(B \rightarrow A)$ is $$ _0\!\ext_B^1(I_B/I_B^2,I_{A/B})= \
_{-s}\!\ext^1_B(I_B/I_B^2,M) - \dim (\im \beta)$$ where $\beta$ is the
homomorphism $\ _{-2s}\!\Ext_B^1(I_B/I_B^2,K_B(t)) \rightarrow \
_{-s}\!\Ext_B^1(I_B/I_B^2,M)$ induced by \eqref{M(-s)}. This codimension also
equals the codimension of the stratum of quotients given by \eqref{M(-s)}
around $(A)$. \\[2mm]
\indent {\rm B)} If \ $(M,B)$ is unobstructed along any graded deformation of
$B$ and $ _{-s}\!\Ext^1_B(I_B/I_B^2,M)= 0$, then $A$ is Gorenstein and the
stratum of quotients given by \eqref{M(-s)} around $(A)$ is open in $
\GradAlg^{H_A}(R)$ (so $A$ is $H_B$-generic). Moreover $A$ is unobstructed as
a graded $R$-algebra and the dimension formula for $\ \dim_{(A)}
\GradAlg^{H_A}(R)$ of part {\rm A)} holds.
\end{theorem}

\begin{remark} \label{remMmainth}
  i) Looking to the proof below we can weaken the assumption ``$B$ is a licci
  quotient'' of Theorem~\ref{Mmainth}B) to ``$B$ is a generically syzygetic
  unobstructed CM quotient satisfying $
  _0\!\Ext^2_B(I_B/I_B^2,K_B(t-2s))=0$'', and conclude as in
  Theorem~\ref{Mmainth}B). \\
  \indent ii) The most natural vanishing condition for
  $_{-s}\!\Ext^1_B(I_B/I_B^2,M)$ in Theorem~\ref{Mmainth}B) seems to be $s \geq
  \max n_{2,j} - a$ where $ \max n_{2,j}$ is the largest degree of a minimal
  relation of $I_B$ (cf.  \eqref{complex2K}) and $a$ is an integer which
  satisfies $M_v = 0$ for $v \leq a$.  In this case we get $
  _{-s}\!\Hom_R(I_B,M)=\ _{-s}\!\Ext_R^1(I_B,M)=0 \ $ and hence $
  \delta(M)_{-s} = 0$. Since $K_B(t)_v=0$ for $v \leq a-s$ by \eqref{M(-s)} we
  get $\ _{-2s}\!\Ext_R^i(I_B,K_B(t))=0 \ $ for $i \leq 1$ and
  $\delta(K_B)_{t-2s}=0$, as well as $\dim(K_B)_{t-2s}=0$. It follows that
  $$
  \ \dim_{(A)} \GradAlg^{H_A}(R)= \ \dim (N_B)_0 + \ \dim (M^*)_s - 1 -
  \gamma(S_2M)_0.$$
  \indent iii) Arguing as in ii) one shows
  $_{-v}\!\Ext_B^i(I_B/I_B^2,K_B) = 0$ for $i \leq 1$ and hence
  $\delta(K_B)_{-v} = 0$ provided $v > \max n_{c,i} + \max n_{2,i} - n - c $,
  e.g. provided $ v > 2\reg(I_B) - n$.
\end{remark}

\begin{remark} \label{remMmainth2} If  $\dim {
    B}-\dim {B}/{I}(Z)\geq 4$ we can replace $ \GradAlg^{H_A}(R)$ by
  $\Hilb^p(\pp) $ in all conclusions, obviously modified ( $A$ by $\Proj(A)$
  etc.), of Theorem~\ref{Mmainth} by  \eqref{Grad}, cf.
  Remark~\ref{remunobstr}. 
\end{remark}

\begin{proof}
  A) We need to verify the assumptions of Corollary~\ref{corvarunobstr}. Since
  $B$ is licci and generically Gorenstein ($\widetilde {K_B}(t)\arrowvert_U$
  is locally free), one knows that $\Ext_B^1(I_B/I_B^2,B) =0$ (\cite{BU83})
  and $\Ext_B^i(I_B/I_B^2,K_B) =0$ for $i \geq 2$ (\cite{VAS}, Thm.4.2.6), and
  that $I_B$ is syzygetic (see \cite{K03}, proof of Thm.\!\ 16). It follows
  from Lemma~\ref{lemMmainth} below and Remark~\ref{remcorvarunobstr} that all
  assumptions of Corollary~\ref{corvarunobstr} hold. To see that the dimension
  formulas of $ \GradAlg^{H_A}(R)$ in Theorem~\ref{Mmainth} and
  Corollary~\ref{corvarunobstr} coincide, we use \eqref{M(-s)} and we get that
  $\delta(I_{A/B})_0 = \delta(M)_{-s} - \delta(K_B)_{t-2s}$ as well as $$\
  _0\!\ext_B^1(I_B/I_B^2,I_{A/B})= \ _{-s}\!\ext^1_B(I_B/I_B^2,M) - \dim (\im
  \beta).$$ Then we conclude by Lemma~\ref{lemMmainth}. For the final
  codimension statement, we refer to Lemma~\ref{lemMmainth2}.
  
  B) A generically Gorenstein licci quotient is generically syzygetic (cf.
  \cite{K03}, proof of Thm.\!\ 16) and satisfies $
  _0\!\Ext^2_B(I_B/I_B^2,K_B(t-2s))=0$. Hence we get Theorem~\ref{Mmainth}B)
  from Proposition~\ref{unobstrp1}$(i)$ and Theorem~\ref{gorgenth}, since the
  dimension formulas of Theorem~\ref{Mmainth} coincide with that in
  Theorem~\ref{gorgenth} by Lemma~\ref{lemMmainth} and Lemma~\ref{lemMmainth2}
  (see part A) of the proof).
\end{proof}

\begin{lemma} \label{lemMmainth}
  Let $B$ be Cohen-Macaulay, let $M$ and $Z$ be as in (the three first lines
  of) Theorem~\ref{Mmainth} and let $A \simeq B/I_{A/B}$ be given by
  \eqref{M(-s)} (so $r = \dim B -\dim A = 2$ and $K_A \simeq A(j)$ where $j =
  2s-t$). If $\ _{s}\!\Ext^1_B(M,B)=0$, then
  $$
  _0\!\hom_B(I_{A/B},A) = \ \dim (M^*)_s - 1 - \gamma(S_2M)_0 + \ 
  \dim(K_B)_{t-2s}. $$
  If in addition $ _t\!\Ext_B^2(S_2(M),K_B)=0$, then
  $_{-j}\!\Ext_B^2(S_2(I_{A/B}),K_B)=0$, i.e.  $_0\!\HH^2(B,A,A)=0$.
\end{lemma}

\begin{proof}
  To show the dimension formula we remark that $ _0\!\hom_B(I_{A/B},A) = \! 
  _{-j}\!\ext_B^{1}(S_2(I_{A/B}),K_B) -1 $ by Lemma~\ref{lemvarunobstr}.
  Moreover note that there exists an exact sequence
\begin{equation} \label{S2(M)}
     M \otimes K_B(t-3s) \rightarrow S_2(M(-s)) \rightarrow
 S_2(I_{A/B})   \rightarrow 0
 \end{equation} 
 whose leftmost map becomes injective if we sheafify and restrict to $U$.
 (Indeed $\wedge^2 \widetilde {K_B}\arrowvert_U = 0$, or if we prefer, the
 morphism $U \cap \Proj(A) \hookrightarrow \Proj (B)$ is an l.c.i.). Let $MK:=
 \ker( S_2(M(-s)) \rightarrow S_2(I_{A/B}) )$ and $\wedge:= \ker( M \otimes
 K_B(t-3s) \rightarrow MK)$. Since $\wedge$ is supported at $Z$, we get $
 \Hom_B( \wedge ,K_B)=0$. It follows that $ \Hom_B( MK,K_B) \simeq \Hom_B( M
 \otimes K_B(t-3s),K_B)$ and that $$ _{t-2s}\!\Ext_B^1(MK,K_B) \subseteq \!
 _{t-2s}\!\Ext_B^1(M \otimes K_B(t-3s),K_B) \simeq \ _0\!\Ext^1_B(M(-s),B) =
 0.$$ Hence \eqref{S2(M)} induces a long exact sequence (i.e. exact in degree
 $t-2s$)
 \begin{equation} \label{extS2(M)}
 ... \rightarrow  \Hom_B( M(t-3s),B) \rightarrow \Ext_B^1(S_2(I_{A/B}),K_B)
 \rightarrow \Ext_B^1(S_2(M(-s)),K_B) \rightarrow 0 .
\end{equation}
Since it is easy to see that the socle degree $j$ of $A$ is $j=2s-t$ from
\eqref{M(-s)}, we get the conclusion by counting dimensions in
\eqref{extS2(M)} in degree $t-2s$ and by observing that $
\Hom_B(S_2(I_{A/B}),K_B) \simeq \Hom_B(I_{A/B}^2,K_B) \simeq
\Hom_B(I_{A/B},K_B) \simeq K_B$. Indeed these isomorphisms follow from the
fact that $\ker(S_2(I_{A/B}) \rightarrow I_{A/B}^2)$ and $\coker(I_{A/B}^2
\hookrightarrow I_{A/B})$ are $A$-modules (see the proof of
Lemma~\ref{lemvarunobstr}).

By Lemma~\ref{lemvarunobstr}, we get  $_0\!\HH^2(B,A,A) \simeq \
_{-j}\!\Ext_B^2(S_2(I_{A/B}),K_B)$. Now continuing the exact sequence
\eqref{extS2(M)} we get an injection $ _{-j}\!\Ext_B^2(S_2(I_{A/B}),K_B)
 \hookrightarrow \ _{t-2s}\!\Ext_B^2(S_2(M(-s)),K_B)$ and we conclude easily.
\end{proof}

\begin{lemma} \label{lemMmainth2}
 Let $B$ be Cohen-Macaulay, let $M$ and $Z$ be as in (the three first lines
  of) Theorem~\ref{Mmainth} and let $A \simeq B/I_{A/B}$ be given by
  \eqref{M(-s)} (so $r = \dim B -\dim A = 2$ and $K_A \simeq A(2s-t)$). Then
   %the codimension of 
  the $H_B$-stratum of $A$ at $(B \rightarrow A)$ and the stratum of quotients
  given by \eqref{M(-s)} around $(A)$ coincide. Moreover if $\ 
  _{0}\!\Ext^1_B(I_B/I_B^2,I_{A/B}) = 0$, then the mentioned strata are open
  in $ \GradAlg^{H_A}(R)$.
\end{lemma}
\begin{proof}
  To see that the $H_B$-stratum and the stratum of quotients given by
  \eqref{M(-s)} coincide, it suffices to show that any $(B' \rightarrow A')$
  in a small enough open neighbourhood of $t:=(B \rightarrow A)$ in $G:=
  \GradAlg(H_B,H_A)$ is given by \eqref{M(-s)}. Let $(S,m_S)$ be the local
  ring of $G$ at $t$ and let $S_i=S/m_S^i$. Since $G$ is a scheme of finite
  type which represents the corresponding functor of graded flat quotients,
  there exists a universal quotient whose pullback to $\Spec(S)$ is denoted by
  $B_S \rightarrow A_S$ (with kernel $I_{A_S/B_S}$, flat over $S$). Since $
  K_{B_S}:= \Ext^{c}_{R_S}(B_S,R_S(-n-c))$ where $R_S:=R \otimes_k S$, is
  $S$-flat (cf. \cite{JS}, Prop.\!\ A1), it suffices to show that the natural
  map $ \ _t\!\Ext^1_{B_S}(I_{A_S/B_S},K_{B_S}(-2s)) \rightarrow \
  _t\!\Ext^1_{B}(I_{A/B},K_{B}(-2s)) $ is surjective because then, there is an
  extension of $I_{A_S/B_S}$ by $K_{B_S}(t-2s)$ over $S$, which reduces to the
  given extension of $I_{A/B}$ by $K_{B}(t-2s)$, and which extends to
  extensions in an open neighbourhood of $t$ in $G$ (the assumptions of
  Theorem~\ref{mainth} needed to get \eqref{M(-s)}, i.e. $\widetilde
  {M'}\arrowvert_{U'}$ is locally free and $\Proj(B')-U'$ is of codimension at
  least two is quite easy to get. Note that the maximal Cohen-Macaulayness of
  $M'$ may be deduced from the extension). To show the surjectivity, let $F_S$
  be a graded $B_S$-free module which surjects into $I_{A_S/B_S}$, and let
  $Q_S:= \ker (F_S \rightarrow I_{A_S/B_S})$. Let $Q_{S_i}:=Q_{S} \otimes_{S}
  S_i$ and $Q:=Q_{S} \otimes_S k$. A simple diagram chasing shows that it
  suffices to prove the surjectivity of $\ _t\!\Hom(Q_S,K_{B_{S}}(-2s))
  \rightarrow \ _t\!\Hom(Q, K_{B}(-2s)) $. As in the proof of \cite{K03},
  Prop.\!\ 13, it suffices to show that $$ \eta_i: \
  _t\!\Hom(Q_{S_i},K_{B_{S_i}}(-2s)) \rightarrow \ _t\!\Hom(Q_{S_{i-1}},
  K_{B_{S_{i-1}}}(-2s))$$ is surjective. Since $\ \Ext^1_{B}(Q,K_{B}) \simeq \
  \Ext^2_{B}(I_{A/B},K_{B})=0$, this follows by applying
  $_{t-2s}\!\Hom(Q_{S_i},-)$ onto $0 \rightarrow B \otimes_k {\mathfrak a}
  \rightarrow B_{S_{i}} \rightarrow B_{S_{i-1}} \rightarrow 0$ where
  ${\mathfrak a} := m_S^{i-1}/m_S^{i}$, i.e. we get the surjectivity of
  $\eta_i$ from $ _t\!\Ext^1_{B}(Q,K_{B}(-2s)) =0$ and we get what we want.
  Finally note that by Prop.\!\ 4$(ii)$ and (5) of \cite{K03} the assumption
  $\ _{0}\!\Ext^1_B(I_B/I_B^2,I_{A/B}) = 0$ implies that the second projection
  $ \GradAlg(H_B,H_A) \rightarrow \GradAlg^{H_A}(R)$ maps small enough open
  sets of $(B \rightarrow A)$ to open sets of $(A)$ and we are done.
  \end{proof}

  Since we need some vanishing results on $\Ext$-groups to use
  Theorem~\ref{Mmainth} effectively, we will delay specific examples to the
  next section. Instead we consider the rank 3 case which will allow us to
  treat certain Gorenstein families of codimension 5 or more. Here there are
  even more $\Ext$-groups involved.  Fortunately they all vanish for ``good''
  modules on licci schemes. Now if $r=3$, then Theorem~\ref{mainth} provides
  us with an exact sequence
\begin{equation} \label{M(-s)3}
0 \rightarrow K_B(t-3s) \rightarrow  M^{\vee}(t-2s) \rightarrow M(-s)
% \stackrel{\sigma}{\rightarrow} 
\rightarrow B \rightarrow A \rightarrow 0 \ ,
\end{equation}
where $ M^{\vee} = \Hom_B(M,K_B)$. Let $$
\ \gamma(S_2M)_v = \ 
_{v}\!\hom_B(S_2(M),K_B(t))-\ _{v}\!\ext_B^1(S_2(M),K_B(t)) +\ 
_{v}\!\ext_B^2(S_2(M),K_B(t)) \ 
$$
(one more term than in the rank 2 case !), let $\gamma(M,M)_v =\ 
_{v}\!\hom_B(M,M)-\ _{v}\!\ext_B^1(M,M)$ while let $ \delta(Q)_{v} =\ 
_{v}\!\hom_B(I_B/I_B^2,Q)-\ _{v}\!\ext_B^1(I_B/I_B^2,Q) \ $ be as previously.

\begin{theorem} \label{Mmainth3} Let $B = R/I_B$ be a graded licci
  quotient of $R$, let $M$ be a graded maximal Cohen-Macaulay B-module, and
  suppose $\widetilde{M}$ is locally free of rank $3$ in $U:=\Proj(B) -Z$,
  that $\dim { B}-\dim {B}/{I}(Z)\geq 3$ and $\wedge^3{\widetilde {
      M}}\arrowvert_U \simeq \widetilde {K_B}(t)\arrowvert_U$. Let $A$ be
  defined by a regular section $\sigma$ of $\widetilde{M}^*(s)$ on $U$, i.e.
  given by \eqref{M(-s)3}
  and suppose $\ _s\!\Ext^1_B(M,B)=0$. \\[2mm]
  \indent {\rm A)} \ \ \ If $ \ _0\!\Ext_B^2(M^{\vee} \otimes M,K_B)= \!
  _{-s}\!\Ext_B^3(S_2(M),K_B(t))=0 \ $ and $\ _{-s}\!\Ext^2_B(I_B/I_B^2,M)= \
  \ _{-2s}\!\Ext^3_B(I_B/I_B^2,M^{\vee}(t))= 0 \ $, then $A$ is unobstructed
  as a graded $R$-algebra (indeed $_0\!\HH^2(R,A,A)=0$), $A$ is Gorenstein,
  and
  $$
  \ \dim_{(A)} \GradAlg^{H_A}(R)= \ \dim (N_B)_0 + \ \dim (M^*)_s +
  \gamma(S_2M)_{-s} - \gamma(M,M)_0 - \ \dim(K_B)_{t-3s} - \delta \ ,$$ where
  $$ \delta:= \delta(M)_{-s} + \delta(K_B)_{t-3s} - \delta(M^{\vee})_{t-2s} -
  _{-2s}\ext^2_B(I_B/I_B^2,M^{\vee}(t)).$$ Moreover if $char(k)=0$ and $(B
  \rightarrow A)$ is general with respect to $\ _{0}\!\hom_R(I_B,I_{A/B}) $,
  then the codimension of the $H_B$-stratum of $A$ at $(B \rightarrow A)$ is
  $$\ _0\!\ext_B^1(I_B/I_B^2,I_{A/B})= \ _{-s}\!\ext^1_B(I_B/I_B^2,M) + \
  _{-2s}\ext^2_B(I_B/I_B^2,M^{\vee}(t)) -\ \dim (\im \beta)\ $$ where $\beta$
  is the homomorphism $\ _{-2s}\!\Ext_B^1(I_B/I_B^2,M^{\vee}(t)) \rightarrow \
  _{-s}\!\Ext_B^1(I_B/I_B^2,M)$ induced by \eqref{M(-s)3}.
%  This codimension also equals the codimension of the stratum of quotients
%  given by \eqref{M(-s)} around $(A)$.
  \\[2mm]
  \indent {\rm B)} If \ $(M,B)$ is unobstructed along any graded deformation
  of $B$ and $ _{-s}\!\Ext^2_B(M,M^{\vee}(t))= \ \ \
  _{-s}\!\Ext_B^2(S_2(M),K_B(t))=0$ and $ _{-s}\!\Ext^1_B(I_B/I_B^2,M)= \
  _{-2s}\!\Ext^2_B(I_B/I_B^2,M^{\vee}(t))= 0$, then $A$ is Gorenstein and
%  stratum of quotients given by \eqref{M(-s)} around $(A)$ is open in $
%  \GradAlg(H_A)$ (so $A$ is 
  $H_B$-generic. Moreover $A$ is unobstructed as a graded $R$-algebra and the
  dimension formula for $\ \dim_{(A)} \GradAlg^{H_A}(R)$ of part {\rm A)} holds
  (this formula simplifies a little due to the assumed vanishing of the 
  $\Ext$-groups).
\end{theorem}

\begin{remark} \label{remMmainth3}

  i) Using \eqref{twoo} one may see that $ _{-s}\!\Ext^2_B(M,M^{\vee}(t))  
  \simeq \ _{-s}\!\Ext_B^2(S_2(M),K_B(t))$ provided  $\depth_{I(Z)}B \geq 4$
  and $char(k) \neq 2$, cf. the assumptions of Theorem~\ref{Mmainth3}B).   
  
  ii) Moreover one may replace ``$ _{-s}\!\Ext^1_B(I_B/I_B^2,M)= \
  _{-2s}\!\Ext^2_B(I_B/I_B^2,M^{\vee}(t))= 0$'' in Theorem~\ref{Mmainth3}B) by
  $$ _{-s}\!\Ext^1_R(I_B,M)= \ _{-2s}\!\Ext^2_R(I_B,M^{\vee}(t))= 0 $$ 
  and still conclude as in Theorem~\ref{Mmainth3}B). This follows from
  Theorem~\ref{gorgenth}. This variation is particularly useful if the
  codimension of $B$ in $R$ is 2, in which case $\ \Ext^2_R(I_B,M^{\vee}(t))$
  vanishes. A natural vanishing condition for $ _{-s}\!\Ext^1_R(I_B,M)$ is
  again $s \geq \max n_{2,j} - a$ where $ \max n_{2,j}$ is the largest degree
  of a minimal relation of $I_B$ and $a$ is an integer which satisfies $M_v =
  0$ for $v \leq a$. In this case we get $\ _{-s}\!\Ext_R^i(I_B,M)=0 \ $ for
  $i=0,1$ and hence $ \delta(M)_{-s} = 0$.
\end{remark}

\begin{lemma} \label{lemMmainth3}
  Let $B$ be Cohen-Macaulay, let $M$ and $Z$ be as in (the three first lines
  of) Theorem~\ref{Mmainth3} and let $A \simeq B/I_{A/B}$ be given by
  \eqref{M(-s)3} (so $r = \dim B -\dim A = 3$ and $K_A \simeq A(j)$ where $j =
  3s-t$) and suppose $\ _{s}\!\Ext^1_B(M,B)=0$.
   
  {\rm A)} \ \ \ If $ \ _0\!\Ext_B^2(M^{\vee} \otimes M,K_B)= 0$, then 
  $$
  _0\!\hom_B(I_{A/B},A) = \ \dim (M^*)_s + \gamma(S_2M)_{-s} - \gamma(M,M)_0 -
  \ \dim(K_B)_{t-3s}$$ If in addition $ _{-s}\!\Ext_B^3(S_2(M),K_B(t))=0$,
  then
  $_{-j}\!\Ext_B^3(S_2(I_{A/B}),K_B) \simeq \ _0\! \HH^2(B,A,A)=0$. \\[2mm]
  \indent {\rm B)} \ \ If $\ _{-s}\!\Ext_B^2(S_2(M),K_B(t))=0$ then $
  _0\!\hom_B(I_{A/B},A)$ is given as in part {\rm A)}.
%  (where, obviously, one term of $ \gamma(S_2M)_{-s}$  vanishes).
\end{lemma}

\begin{proof}
  A) Thanks to  Lemma~\ref{lemvarunobstr} we have  $ _0\!\hom_B(I_{A/B},A) =
  \!  _{-j}\!\ext_B^{2}(S_2(I_{A/B}),K_B) -1 $. To compute $\  \!
  _{-j}\!\ext_B^{2}(S_2(I_{A/B}),K_B)$, we look to the exact sequence
\begin{equation} \label{S2(M)3}
     M^{\vee} \otimes M (t-3s) \rightarrow S_2(M(-s)) \rightarrow
 S_2(I_{A/B})   \rightarrow 0
 \end{equation} 
 where one knows that the kernel of $ \widetilde {M^{\vee} \otimes
   M}(t-3s)\arrowvert_U \rightarrow \widetilde {S_2(M)}(-2s) \arrowvert_U$ is
 $(\wedge^2 \widetilde {M^{\vee}}(2t-4s) \oplus \widetilde { K_B}(t-3s))
 \arrowvert_U \ $ ($U \cap \Proj(A) \hookrightarrow \Proj (B)$ is an l.c.i.).
 Let $M^{\vee}M:= \ker( S_2(M(-s)) \rightarrow S_2(I_{A/B}) )$.  Since
 $\wedge^2 \widetilde {M^{\vee}}\arrowvert_U \simeq \widetilde M \otimes
 \widetilde { K_B}(-t) \arrowvert_U $, we have an exact sequence
 $$
 0 \rightarrow ( \widetilde {K_B}(t-3s) \oplus ( \widetilde M \otimes
 \widetilde { K_B})(t-4s)) \arrowvert_U \rightarrow \widetilde {M^{\vee}
   \otimes M}(t-3s)\arrowvert_U \rightarrow \widetilde {M^{\vee}M}
 \arrowvert_U \rightarrow 0 \ $$
 to which we apply $ \Hom( - , \widetilde
 {K_B}\arrowvert_U)$. Since $\Ext_B^i(-, N) \simeq \Ext_{\sO_U}^i(\widetilde -
 \arrowvert_U,\widetilde N\arrowvert_U) $ for $i \leq 1$ for every maximal
 CM module $N$ by \eqref{twoo}, we get a long exact sequence starting as
 $$
 0 \rightarrow \Hom_B( M^{\vee}M,K_B(-j)) \rightarrow \Hom_B( M^{\vee}
 \otimes M(t-3s),K_B(-j)) \rightarrow \Hom_B(K_B,K_B) \oplus H \rightarrow ...
 $$
 where $ H:= \Hom_B( M \otimes K_B(t-4s),K_B(-j)) \simeq \Hom_B( M(-s),B) $
 and stopping at $ _0\! \Ext_B^1( M \otimes K_B(t-4s),K_B(-j)) \simeq \ _{s}\!
 \Ext_B^1( M,B) = 0$. Noting that $$_0\!\Ext_B^i(M^{\vee} \otimes M,K_B)
 \simeq \Ext_{\sO_U}^i(\widetilde M^* \otimes \widetilde {K_B} \otimes
 \widetilde M \arrowvert_U, \widetilde {K_B}\arrowvert_U) \simeq \ _0\!
 \Ext_B^i(M, M) \ \ {\rm for} \ \ i \leq 1,$$ we get
 $$
 \sum_{i=0}^{1} (-1)^i \! _{-j}\!\ext_B^i(M^{\vee}M,K_B)= \sum_{i=0}^{1}
 (-1)^i \! _{0}\!\ext_B^i(M, M) -1 -\! _s\! \hom_B(M,B) \ . $$
 Note that the
 long exact sequence above and the assumption $_0\!\Ext_B^2(M^{\vee} \otimes
 M,K_B)= 0$ also show $ _0\!\Ext_B^2(M^{\vee} M,K_B)= 0$. Now combining with
 the long exact sequence which we get by applying $ _0\! \Hom( - ,K_B(-j))=0$
 to $ 0 \rightarrow M^{\vee}M \rightarrow S_2(M)(-2s) \rightarrow S_2(I_{A/B})
 \rightarrow 0$;
 \begin{equation} \label{extS2(M)3}
 \small
 ... \rightarrow \  _0\! \Ext_B^1(M^{\vee} M,K_B(-j)) \rightarrow \
  _{0}\! \Ext_B^2(S_2(I_{A/B}),K_B(-j))  \rightarrow \  _{0}\!
 \Ext_B^2(S_2(M)(-2s),K_B(-j)) \rightarrow 0  
\end{equation}
which implies
$$
\sum_{i=0}^{2} (-1)^i \! _{-j}\!\ext_B^i(S_2(I_{A/B}),K_B)=
\gamma(S_2M)_{-s} - \sum_{i=0}^{1} (-1)^i \! _{-j}\!\ext_B^i(M^{\vee}M,K_B)\ ,
$$
we get the dimension formula because $$ \ \Hom_B(S_2(I_{A/B}),K_B) \ \simeq \
\Hom_B(I_{A/B}^2,K_B)\ \simeq \ \Hom_B(I_{A/B},K_B) \simeq K_B$$ and
$\Ext_B^1(S_2(I_{A/B}),K_B) \simeq \Ext_B^1(I_{A/B},K_B)=0 \ $
($\ker(S_2(I_{A/B}) \rightarrow I_{A/B}^2)$ and $\coker(I_{A/B}^2
\hookrightarrow I_{A/B})$ are $A$-modules). Part B) is proven in exactly the
same way (the sequence \eqref{extS2(M)3} stops by ``one module earlier'').
Moreover continuing \eqref{extS2(M)3}, we get the ``if in addition'' statement
of part A) from the assumption $ _{-s}\!\Ext_B^3(S_2(M),K_B(t))=0$.
\end{proof}

\begin{proof}[Proof of Theorem~\ref{Mmainth3}]
  A) follows from Corollary~\ref{corvarunobstr}, Remark~\ref{remcorvarunobstr}
  and Lemma~\ref{lemMmainth3} because $B$ is licci and generically Gorenstein
  (cf. the proof of Theorem~\ref{Mmainth} for details). Note that we need to
  split \eqref{M(-s)3} into two short exact sequences, and apply $
  _0\!\Hom_B(I_B/I_B^2,-)$ to them, to see the formula of $\delta$ in
  Theorem~\ref{Mmainth3} ($\delta$ obviously equals $\delta(I_{A/B})_0$ by
  definition of the latter). The same splitting, together with $\
  _{-s}\!\Ext^2_B(I_B/I_B^2,M)= 0$ and $ _{-3s}\!
  \Ext^i_B(I_B/I_B^2,K_B(t))=0$ for $i > 1$, shows that $$
  _0\!\ext_B^1(I_B/I_B^2,I_{A/B})= \ _{-s}\!\ext^1_B(I_B/I_B^2,M)+ \
  _{-2s}\!\ext^2_B(I_B/I_B^2,M^{\vee}(t)) - \dim (\im \beta),$$ and we get the
  codimension statement.
  
  B) We get Theorem~\ref{Mmainth3}B) from Proposition~\ref{unobstrp1}$(i)$,
  Theorem~\ref{gorgenth} and Lemma~\ref{lemMmainth3} (see part A) of the proof
  for the dimension formula).
\end{proof}

\begin{remark} \label{rem2Mmainth3}
If we in Theorem~\ref{Mmainth}A) and  Theorem~\ref{Mmainth3}A) replace
the vanishing of all $ _v\! \Ext_B^i$-groups by the vanishing of the
corresponding  $\ \Ext_B^i$-groups (skipping the index $v$), we get that $A$
is strongly unobstructed in the sense  $\ \HH^2(R,A,A)=0$ (cf.
Remark~\ref{remvarunobstr}).
\end{remark}

\section{Rank two and three sheaves on codimension two quotients}
In this section we will see that Theorem~\ref{Mmainth} applies to the normal
module, $M=N_B$ and Theorem~\ref{Mmainth} and \ref{Mmainth3} to the 1. Koszul
homology module, $M=H_1$. Of course Theorem~\ref{Mmainth} also applies to
modules of the form $M=S_2(K_B) \oplus K_B^*(t)$, where $t$ is an integer. For
such decomposable modules Theorem~\ref{Mmainth} does not really lead to new
results since they may be treated by applying Theorem~\ref{Kmainth}.

Before giving our applications we need two propositions to handle the
vanishing of the $\Ext$-groups involved. In what follows, $B=R/I$ is an
$n$-dimensional codimension two CM quotient of $R$,
\begin{equation} \label{Icodim2}
  0 \rightarrow  G_2:= \oplus_{j=1}^{\mu -1} R(-n_{2,j}) \rightarrow
  G_1:= \oplus_{i=1}^{\mu} R(-n_{1,i}) \rightarrow I \rightarrow 0  
\end{equation}
is a minimal resolution and $Y=\Proj(B)$ is an l.c.i. in an open set $U=Y-Z$
where $\depth_{I(Z)}B \geq 1$. Taking $R$-duals, we get a minimal resolution
 \begin{equation} \label{Kcodim2}
  0 \rightarrow R \rightarrow  \oplus R(n_{1,i}) \rightarrow
  \oplus R(n_{2,j}) \rightarrow K_B(n+2) \rightarrow 0  
\end{equation}
to whom we apply $\Hom(-,B)$ to see the exactness to the left in the exact
sequence 
\begin{equation}  \label{K*codim2}
 0 \rightarrow K_B(n+2)^* \rightarrow \oplus
  B(-n_{2,j}) \rightarrow \oplus B(-n_{1,i})\rightarrow I/I^2 \rightarrow 0
\end{equation}
Note that \eqref{K*codim2} splits into two short exact sequences ``via $
\oplus B(-n_{2,j}) \twoheadrightarrow H_1 \hookrightarrow  \oplus
B(-n_{1,i})$'', one of which is \eqref{es1} with $ \HH_2(R,B,B) =0$.

\begin{proposition} \label{extI}
  Let $l$ be a natural number, let $B=R/I$ be a codimension two CM quotient of
  $R$, and suppose $\depth_{I(Z)}B \geq l$.  Then
  $$\ \Ext_B^i(I/I^2,I/I^2)=0 \ \ {\rm for} \ 1 \leq i \leq l. $$
\end{proposition}

\begin{proof} Since $pd_R I = 1$, we have $\ \Ext_R^i(I,I/I^2) = 0 $ and $\
  \Tor_i^R(I,B)=0$ for $i \geq 2$. Moreover, using the spectral sequence $\
  \Ext_B^p( \Tor_q^R(I,B),I/I^2)$ which converges to $\
  \Ext_R^{p+q}(I,I/I^2)$, we get isomorphisms $\ \Ext_B^{i-2}(
  \Tor_1^R(I,B),I/I^2) \simeq \Ext_B^i(I/I^2,I/I^2) $ for $i > 2$ and an exact
  sequence $$\ 0 \rightarrow \Ext_B^1(I/I^2,I/I^2) \rightarrow
  \Ext_R^1(I,I/I^2) \rightarrow \Hom_B( \Tor_1^R(I,B),I/I^2) \rightarrow
  \Ext_B^2(I/I^2,I/I^2) \rightarrow 0.$$ Note that $ \Ext_R^1(I,-)$ is right
  exact by $pd_R I = 1$. In particular $ \Ext_R^1(I,R) \simeq \Ext_R^1(I,A)$
  and it follows that $ N_B \simeq \Ext_R^1(I,I) \simeq \Ext_R^1(I,R) \otimes
  I \simeq K_B(n+2) \otimes I/I^2$. Similarly we get $$ \Ext_R^1(I,I/I^2)
  \simeq \Ext_R^1(I,R) \otimes I/I^2 \simeq N_B \ .$$ Hence there is an
  injection $\HH_{I(Z)}^0( \Ext_B^1(I/I^2,I/I^2)) \hookrightarrow
  \HH_{I(Z)}^0(N_B)$. Since $ \widetilde {I/I^2} \arrowvert_U $ is locally
  free, we get $\HH_{I(Z)}^0( \Ext_B^1(I/I^2,I/I^2)) \simeq
  \Ext_B^1(I/I^2,I/I^2)$. Moreover one knows that $N_B$ is a maximal CM module
  (cf. the text after \eqref{es1}), and we conclude that $
  \Ext_B^1(I/I^2,I/I^2)=0$.
  
  Let $p \leq l$ be a natural number and suppose we have proved $
  \Ext_B^i(I/I^2,I/I^2)=0$ for $1 \leq i < p$. By the spectral sequence (resp.
  the two displayed formulas) above, it suffices to show $ \Ext_B^{p-2}(
  \Tor_1^R(I,B),I/I^2) = 0$ for $p > 2$ (resp. $ \Ext_B^{p-2}(
  \Tor_1^R(I,B),I/I^2) \hookrightarrow N_B$ a graded injection for $p = 2$).
  We have $\HH_{I(Z)}^i(I/I^2)=0$ for $i \leq p-2$ by \eqref{es1}. Hence we
  get an injective graded map 
 \begin{equation} \label{Tor}
 \small
 \Ext_B^{p-2}( \Tor_1^R(I,B),I/I^2) \hookrightarrow
  \Ext_{\sO_U}^{p-2}(\widetilde {\Tor_1^R(I,B)} \arrowvert_U, \widetilde
  {I/I^2} \arrowvert_U) \simeq \Ext_{\sO_U}^{p-2}(\widetilde
  {\Tor_1^R(I,K_B(n+2))} \arrowvert_U, \widetilde {N_B} \arrowvert_U),
 \end{equation}
 (cf. \eqref{twoo}), noting that $ \widetilde {\Tor_1^R(I,K_B)} \simeq
 \widetilde {\Tor_1^R(I,B) \otimes K_B}$ and $ \widetilde {I/I^2 \otimes K_B}
 \simeq \widetilde {N_B}$ are isomorphic on $U$. Since $\HH_{I(Z)}^i(N_B)=0$
 for $i \leq p-1$, the rightmost $\Ext$-group in \eqref{Tor} is, by
 \eqref{twoo} mainly, further isomorphic to $
 \Ext_{B}^{p-2}(\Tor_1^R(I,K_B(n+2)), N_B)$. We {\it claim} that $
 \Tor_1^R(I,K_B(n+2)) \simeq R/I$. Indeed we have by \eqref{Kcodim2} that $
 \Tor_1^R(I,K_B(n+2))$ is the homology group in the middle of the complex $$ 0
 \rightarrow I \rightarrow \oplus I(n_{1,i}) \rightarrow \oplus I(n_{2,j})
 \rightarrow 0.$$ Applying $\Hom_R(-,I)$ to \eqref{Icodim2}, we see that
 $\ker[ \oplus I(n_{1,i}) \rightarrow \oplus I(n_{2,j})] \simeq \Hom(I,I)
 \simeq R $ and putting things together we get the claim. Now using the claim
 we get $ \Ext_B^{p-2}( \Tor_1^R(I,B),I/I^2) \simeq \Ext_{B}^{p-2}(R/I, N_B) $
 which vanishes for $p>2$ and equals $N_B$ for $p = 2$, and we are done.
\end{proof}

\begin{remark}  \label{remextI}
 Applying  $\Hom_R(-,I/I^2)$ to \eqref{Icodim2}, noting that
  $\Hom_R(I,I/I^2) \simeq \Hom_B(I/I^2,I/I^2)$ and $ \Ext_R^1(I,I/I^2) \simeq
  N_B$ (see the proof above) we get the following exact sequence
  $$
  0 \rightarrow \Hom_B(I/I^2,I/I^2) \rightarrow \oplus I/I^2(n_{1,i})
  \rightarrow \oplus I/I^2 (n_{2,j}) \rightarrow N_B \rightarrow 0.$$
  This
  sequence of graded $B$-modules can be used to find $_v\!\hom_B(I/I^2,I/I^2)
  $ because $ \dim (I^2)_v$ and $ \dim (N_B)_v$ can be computed from
 \begin{equation} \label{I2codim2}
  0 \rightarrow \wedge^2( \oplus R(-n_{2,j})) \rightarrow  (\oplus R(-n_{1,i}))
  \otimes ( \oplus R(-n_{2,j})) \rightarrow 
  S_2( \oplus R(-n_{1,i})) \rightarrow I^2 \rightarrow 0 \ , 
\end{equation}
 \begin{equation} \label{NBcodim2}
  0 \rightarrow G_1^* \otimes_R G_2 \rightarrow  (( G_1^* \otimes_R G_1)
  \oplus (G_2^* \otimes_R G_2))/R \rightarrow 
  G_2^* \otimes_R G_1 \rightarrow N_B \rightarrow 0 \  
\end{equation}
Indeed the latter sequence is deduced from the exact sequence $ 0 \rightarrow
R \rightarrow \oplus I(n_{1,i}) \rightarrow \oplus I(n_{2,j}) \rightarrow N_B
\rightarrow 0$ which we get by applying $\Hom_R(-,I)$ to \eqref{Icodim2},
(cf. \cite{K03},\! (26)). Note that if $s > \max n_{2,j} -2 \min n_{1,i}$,
then $ I^2(n_{i,j})_{-s}=0$ for all $i,j$ and we get $
_{-s}\!\Hom_B(I/I^2,I/I^2) \simeq R_{-s}$.
\end{remark}

Now we consider the 1.Koszul homology module $H_1$. In this case all
assumptions of  Theorem~\ref{Mmainth}A) are satisfied, due to  

\begin{proposition} \label{propH}
  Let $B=R/I$ be a codimension two CM quotient of $R$, and suppose
  $\depth_{I(Z)}B \geq 2$.  Then $S_2(H_1)$ is a maximal CM $B$-module.
  Moreover $$\ \ \Ext_B^1(H_1,H_1) \simeq \Ext_B^2(I/I^2,H_1)=0 \ \ {\rm and}
  \ \ \Ext_B^1(H_1,B)=0 . $$
\end{proposition}

\begin{proof} In the sequence \eqref{es1}, $ \HH_2(R,B,B) =0$. Applying
  $\Hom(I/I^2,-)$ to \eqref{es1} we get isomorphisms $\ \Ext_B^1(I/I^2,I/I^2)
  \simeq \Ext_B^2(I/I^2,H_1)$ because $\Ext_B^i(I/I^2,B)=0$ for $1 \leq i \leq
  2$ (\cite{KP1}, (7)). Hence $ \Ext_B^2(I/I^2,H_1)=0$ by
  Proposition~\ref{extI}. Moreover applying $\Hom(-,B)$ (resp.  $\Hom(-,H_1)$)
  to \eqref{es1} we get $ \Ext_B^1(H_1,B) \simeq \Ext_B^2(I/I^2,B)=0 $ 
  (resp. $ \Ext_B^1(H_1,H_1) \simeq \Ext_B^2(I/I^2,H_1)$). 
  Finally to see that $S_2(H_1)$ is a maximal CM module, we consider the short
  exact sequence deduced from \eqref{K*codim2} ``to the left''. It induces an
  exact sequence
\begin{equation} \label{S2(H)}
0 \rightarrow  K_B(n+2)^* \otimes ( \oplus B(-n_{2,j}))
\stackrel{\psi}{\rightarrow}  S_2( 
\oplus B(-n_{2,j})) \rightarrow S_2(H_1)   \rightarrow 0
 \end{equation} 
 because $\ker \psi$, which is supported at $Z$ and contained in a maximal
CM module, must vanish. It follows from \eqref{S2(H)} that $ S_2(H_1)$ has
codepth at most one. Moreover dualizing \eqref{S2(H)} we get
\begin{equation} \label{S2(H*)}
0 \rightarrow \Hom_B( S_2(H_1), K_B) \rightarrow K_B \otimes S_2( \oplus
B(n_{2,j})) \rightarrow  \oplus S_2(K_B)(n_{2,j}+n+2)    
 \end{equation} 
 because $\ \Hom_B( S_2(\oplus B(-n_{2,j})),K_B) \ \simeq \ K_B \otimes S_2(
 \oplus B(n_{2,j}))$ and $$ \Hom_B( K_B(n+2)^* \otimes ( \oplus B(-n_{2,j}))
 ,K_B) \simeq \oplus S_2(K_B)(n_{2,j}+n+2).$$ Note that the cokernel of $ K_B
 \otimes S_2( \oplus B(n_{2,j})) \rightarrow \oplus S_2(K_B)(n_{2,j }+n+2)$ is
 $ \Ext_B^1( S_2(H_1), K_B)$. One may, however, easily see that this map is
 surjective because it is, via symmetrization and tensorization, obtained from
 the surjective map $ \oplus B(n_{2,j}) \rightarrow K_B(n+2)$. Hence $
 \Ext_B^1( S_2(H_1), K_B)=0$, and by Gorenstein duality and the fact that the

 codepth of $S_2(H_1)$ is at most one, we get that $S_2(H_1)$ is maximally
 Cohen-Macaulay.
\end{proof}

If the number of minimal generators of $I$ is $\mu(I)=4$, then the rank of
$H_1$ is $r=2$ by \eqref{es1} and Theorem~\ref{Mmainth}A) applies to
$B$-module $M=H_1$. Since all assumptions of Theorem~\ref{mainth} are
satisfied (\cite{KP1} and Proposition~\ref{propnormal}) we have an exact
sequence
\begin{equation} \label{H1(-s)}
0 \rightarrow K_B(t-2s)\rightarrow H_1(-s) \rightarrow  I_{A/B} \rightarrow 0
\ 
\end{equation}
with $t = n+2- \sum n_{1,i}$.  Hence
\begin{corollary} \label{corH} Let $B=R/I$ be a graded codimension two CM
  quotient of $R$, let $U= \Proj(B)-Z \hookrightarrow \pp^{n+1}$ be an l.c.i.
  and suppose $\mu(I)=4$ and $\depth_{I(Z)}B \geq 2$. If $A$ is defined by a
  regular section of $\widetilde{H_1}^*(s)$ on $U$, i.e. given by
  \eqref{H1(-s)}, then $A$ is unobstructed as a graded $R$-algebra (indeed
  $_0\!\HH^2(R,A,A)=0$), $A$ is Gorenstein of codimension 4 in $R$, and
  {\small
  $$
  \ \dim_{(A)} \GradAlg(H_A)= \ \dim (N_B)_0 + \ \dim (H_1^*)_s - 1 -
  \gamma(S_2H_1)_0 + \ \dim(K_B)_{t-2s} +\delta(K_B)_{t-2s} -
  \delta(H_1)_{-s}.$$}Moreover if $char(k)=0$ and $(B \rightarrow A)$ is
general with respect to $\ _{0}\!\hom_R(I_B,I_{A/B}) $, then the stratum of
quotients given by \eqref{H1(-s)} around $(A)$ is irreducible and its
codimension is $\ _{-s}\!\ext^1_B(I/I^2,H_1) - \dim (\im \beta)$ where $\beta$
is the homomorphism $\ _{-2s}\!\Ext_B^1(I/I^2,K_B(t)) \rightarrow \
_{-s}\!\Ext_B^1(I/I^2,H_1)$ induced by \eqref{H1(-s)}. Furthermore if $(B')
\in \GradAlg^{H_B}(R)$ satisfies the same assumptions as $B$ above and defines
$A'$ as $B$ defined $A$, then the closures in $ \GradAlg^{H_A}(R)$ of the
stratum of quotients given by \eqref{H1(-s)} around $(A)$ and the
corresponding stratum around $(A')$ coincide. If in addition $\ s > \max
n_{2,j}-\min n_{2,j} $, then the stratum of quotients given by \eqref{H1(-s)}
around $(A)$ is open (so $A$ is $H_B$-generic), and $$\ \ \dim_{(A)}
\GradAlg^{H_A}(R)= \ \dim (N_B)_0 + \ \dim (H_1^*)_s - 1 - \ _0\!\hom_B(
S_2(H_1), K_B(t)) \ . $$
\end{corollary}

Note that $ \dim (H_1^*)_s $ is easily computed from the exact sequence
\begin{equation} \label{Ncodim2}
  0 \rightarrow  N_B \rightarrow
  \oplus_{i=1}^{\mu} B(n_{1,i}) \rightarrow H_1^* \rightarrow 0  
\end{equation}
(the $B$-dual sequence of \eqref{es1} is short exact due to $
\Ext_B^1(I/I^2,B)=0$). Moreover $\gamma(S_2H_1)_0=\ _0\!\hom_B( S_2(H_1),
K_B(t))$ is computed from \eqref{S2(H*)} and the minimal resolution
(\cite{EIS}, p.\! 595), 
\begin{equation} \label{S2Kcodim2}
  0 \rightarrow \wedge^2( \oplus R(n_{1,i})) \rightarrow  (\oplus R(n_{1,i}))
  \otimes ( \oplus R(n_{2,j})) \rightarrow 
  S_2( \oplus R(n_{2,j})) \rightarrow S_2(K_B)(2n+4) \rightarrow 0 \  
\end{equation}

\begin{proof} Since all assumptions of Theorem~\ref{Mmainth}A) are taken care
  of by Proposition~\ref{propH}, we get Corollary~\ref{corH} from
  Theorem~\ref{Mmainth}A) except the irreducibility of the stratum, its
  uniqueness (i.e. that the strata above coincide, up to closure) and the
  final statement of Corollary~\ref{corH}. The irreducibility is trivial
  because the stratum is the image of an irreducible set. Indeed
  $\GradAlg(H_B,H_A)$ is smooth at $(B \rightarrow A)$ because $B$ is
  unobstructed and $_0\!\HH^2(B,A,A)=0$ by Lemma~\ref{lemMmainth}. The proof
  of the uniqueness is essentially the same as for \cite{K03}, Prop.\!\
  23$(i)$ and Thm.\!\ 24 (see the two first lines of the proof of \cite{K03},
  Thm.\!\ 24) because the open subscheme of $ \GradAlg^{H_B}(R)$
  consisting of CM quotients is irreducible. Moreover, since $
  \oplus B(-n_{2,i}) \twoheadrightarrow H_1$ is surjective (cf.
  \eqref{K*codim2} ) we have $(H_1)_v=0$ for $v < \min n_{2,j} $. By
  Remark~\ref{remMmainth}$(ii)$ and the assumption $s > \max n_{2,j}-\min
  n_{2,j}$, it follows that $ _{-s}\!\Ext_B^1(I/I^2,H_1)=0$. Hence $A$ is
  $H_B$-generic and the formula for $ \ \dim_{(A)} \GradAlg^{H_A}(R)$
  simplifies as in Corollary~\ref{corH} according to
  Remark~\ref{remMmainth}$(ii)$.
\end{proof}

We illustrate Corollary~\ref{corH} by a example in which we compute all
numbers in the dimension formula of $ \ \dim_{(A)} \GradAlg^{H_A}(R)$ as well
as the number leading to the codimension of the $H_B$-stratum. Even though the
example is relatively simple, our methods of computation may be used quite
generally (to treat sections of $H_1^*(s)$ on a codimension 2 CM quotient of
$R$ generated by $ \mu(I)=4$ generators).

\begin{example}
  It is well known how to find a minimal resolution of $H_1$. Indeed by
  \cite{AH}, 
 \begin{equation} \label{wedgeH1}
0 \rightarrow \wedge^2 G_2 \rightarrow \wedge^2 G_1 \rightarrow G_2
\rightarrow  H_1 \rightarrow 0
\end{equation} 
where the free modules $G_i$ belong to \eqref{Icodim2}. Applying
the mapping cone construction to \eqref{H1(-s)} we get the following
resolution of the Gorenstein algebra $A=R/I_A$ of Corollary~\ref{corH},
\begin{equation} \label{wedgeH1A}
0 \rightarrow R(b) \rightarrow  G_1^*(b) \oplus \wedge^2 G_2(-s) \rightarrow
G_2^*(b) \oplus \wedge^2 G_1(-s) \oplus G_2 \rightarrow G_2(-s) \oplus G_1
\rightarrow I_A \rightarrow 0 
\end{equation} 
where $b:= -\sum n_{1,i}-2s$ and $ \mu(I)=4$. As a special case we suppose the
resolution \eqref{Icodim2} is linear (i.e. all $n_{1,i}=3$). We get the
following resolution (obviously minimal for $s > 0$)  of $A$,
\begin{equation} \label{H1A}
 \begin {aligned}
  & 0 \rightarrow R(-12-2s) \rightarrow  R(-9-2s)^4 \oplus R(-8-s)^8 \\ 
 & \rightarrow R(-8-2s)^3 \oplus R(-6-s)^6 \oplus R(-4)^3 \rightarrow
 R(-4-s)^3 \oplus R(-3)^4  \rightarrow R \rightarrow A 
 \rightarrow 0 \ \ . 
 \end{aligned}
\end{equation} 
Let $R$ be a polynomial ring in four variables. Then $Y= \Proj(B)$ is a curve
of degree $d=6$ and with Hilbert polynomial $p_Y(v)= 6v-2$. Moreover $A$ is
Artinian of socle degree $2s+8$ and with $h$-vector
$(1,4,10,..,6s+16,6s+19,6s+16,..,10,4,1)$. We suppose $s \geq -1$ to avoid
discussing very special cases. If $Y$ is an l.c.i., then $A$ is unobstructed
by Corollary~\ref{corH} of a dimension which we now calculate. Indeed $\ \dim
(N_B)_0 = h^0(\widetilde {N_B}) = 4d=24$ while, by \eqref{Ncodim2}, $$ \ \dim
(H_1^*)_s = 4 \cdot p_Y(s+3)- h^0(\widetilde {N_B}(s)) = 24(s+3)-8 - 4d -2ds =
12s+40$$ for $s \geq -1$ (one may see that $\ (N_B)_s = 12 = 4d+2ds$ also for
$ s = -1$ by \eqref{NBcodim2}). Moreover $$ _0\!\hom_B( S_2(H_1), K_B(-8)) = 6
\dim (K_B)_0 - 3 \dim (S_2(K_B)_0) = 0$$ by \eqref{S2(H*)}, \eqref{Kcodim2}
and \eqref{S2Kcodim2}. Hence by Corollary~\ref{corH}, if $s>0$, then $ \
\dim_{(A)} \PGor(H_A)=\ $
$$
\ \dim_{(A)} \GradAlg^{H_A}(R)= \ \dim (N_B)_0 + \ \dim (H_1^*)_s - 1 - \
_0\!\hom_B( S_2(H_1), K_B(-8)) = 12s+63$$ and $A$ is $H_B$-generic, i.e. the
closure of the family given by \eqref{H1(-s)} forms a $(12s+63)$-dimensional
generically smooth, irreducible component of $ \GradAlg^{H_A}(R)$. It is
interesting to observe that we can also use Lemma~\ref{lemvarunobstr}$(iii)$
to compute $_0\!\hom_B(I_{A/B},A)$ and hence $\ \dim (H_1^*)_s - 1 - 
\ _0\!\hom_B( S_2(H_1), K_B(-8))$ (cf. Lemma~\ref{lemMmainth}). Indeed since
we get $\dim S_2(I_{A/B})_{2s+8}=6$ from \eqref{H1A}, it follows that $$\
_0\!\hom_B(I_{A/B},A) = \dim B_{2s+8} -\dim S_2(I_{A/B})_{2s+8}=12s+40$$ which
again leads to $ \ \dim_{(A)} \GradAlg^{H_A}(R)= 12s+63$.

Finally suppose $s \leq 0$. We still have $\delta(K_B)_{-8-2s}= 0 $ and $
_{-2s}\!\Ext_B^1(I/I^2,K_B(-8))=0$ for $s > -2$ by
Remark~\ref{remMmainth}$(iii)$ while $(K_B)_{-8-2s}= 0$ by \eqref{Kcodim2}. It
remains to compute $ \delta(H_1)_{-s} $ and $\ _{-s}\!\ext^1_B(I/I^2,H_1)$. To
do so we apply $\Hom(I/I^2,-)$ to \eqref{es1}, and we get
 \begin{equation} \label{deltacodim2}
 \small
  0 \rightarrow  \Hom(I/I^2,H_1) \rightarrow
 \Hom(I/I^2,\oplus B(-n_{1,i})) \rightarrow \Hom(I/I^2,I/I^2) \rightarrow
 \Ext^1_B(I/I^2,H_1) \rightarrow 0  
\end{equation}
Note that $\ (N_B)_v = 0$ for $ v < -1$ by \eqref{NBcodim2}. Hence $$\
_{-s}\!\ext^1_B(I/I^2,H_1) = \! _{-s}\!\hom(I/I^2,I/I^2) = \dim R_{-s}$$ and $
\delta(H_1)_{-s} = \ - _{-s}\!\ext^1_B(I/I^2,H_1)$ for $-1 \leq s \leq 0$ by
Remark~\ref{remextI}. By Corollary~\ref{corH}, $$ \ \dim_{(A)}
\GradAlg^{H_A}(R)= 12s+63+ _{-s}\!\ext^1_B(I/I^2,H_1) \ \ {\rm for} \ \ -1
\leq s \leq 0,$$ and the codimension, $_{-s}\!\ext^1_B(I/I^2,H_1)$, of the
$H_B$-stratum of $A$ is $1,4$ for $s=0,-1$ respectively.
\end{example}

\begin{remark} \label{remcodim} In view of Remark~\ref{remextI}, we can always
  use \eqref{deltacodim2} to find $\delta(H_1)_{-s}$. Moreover
  $\delta(K_B)_{t-2s}$ is always computable from the exact sequence we get by
  applying $\Hom(-,K_B(t))$ to \eqref{es1}, cf. \cite{K03}, Remark 14$(d)$.
  Since $\dim (K_B)_{t-2s}$ is given by \eqref{Kcodim2}, we see that all
  members of the dimension formula for $ \GradAlg^{H_A}(R)$ in
  Corollary~\ref{corH} are easily computed, even by hand.  On the other side,
  the formula for the codimension is not always straightforward. However, in
  the range $\ \max n_{2,j}- 2 \min n_{1,i} < s \leq \max n_{2,j}-\min n_{2,j}
  \ $, one may show
  $$_{-s}\!\ext^1_B(I/I^2,H_1) = -\delta(H_1)_{-s}= \dim R_{-s}.$$
  Indeed by
  \eqref{deltacodim2} it suffices to show $N_B(- n_{1,i})_{-s}=0$ for every
  $i$. Since $\oplus R( n_{2,j}- n_{1,i'}) \rightarrow N_B$ is surjective and
  $\oplus R( n_{2,j}- n_{1,i'}- n_{1,i})_{-s} = 0$ we conclude by
  Remark~\ref{remextI}.  Finally note that in the range $\ \max n_{2,j}- 2
  \min n_{1,i} < s \leq \max n_{2,j}-\min n_{2,j} \ $,
  $_{-s}\!\ext^1_B(I/I^2,H_1)$ is the codimension of the $H_B$-stratum of
  Corollary~\ref{corH} because we have $ _{-2s}\!\Ext_B^1(I/I^2,K_B(t))=0$.
  Indeed by Remark~\ref{remMmainth}$(iii)$ we get the last mentioned vanishing
  provided $2s-t > 2 \max n_{2,j}-n-2$, i.e.  provided $s > \max n_{2,j}- \sum
  n_{1,i}/2$ (by $t = n+2- \sum_{i=1}^{i=4} n_{1,i}$) which holds since $
  \sum_{i=1}^{i=4} n_{1,i}/2 \geq 2 \min n_{1,i} $.
\end{remark}

\begin{example} (Arithmetically Gorenstein curves $\Proj(A)$ in $ \pp^5$,
  obtained by \eqref{H1(-s)}.) Here we reconsider the preceding example, in
  dimension two higher. Let $R$ be a polynomial ring in 6 variables, and let
  $Y= \Proj(B)$ be a threefold satisfying $\depth_{I(Z)}B \geq 2$ whose
  resolution \eqref{Icodim2} is linear with $\mu = 4$. Let $A$ be defined by a
  regular section of $\widetilde{H_1}^*(s) \arrowvert_{Y-Z}$, so $X =
  \Proj(A)$ is an AG curve in $ \pp^5$. The Hilbert polynomial of $Y$ is given
  by ``integrating'' $6v-2$ two times, or more precisely $p_Y(v)= 6
  \binom{v+2}{3} - 2 \binom{v+2}{2} + 3 \binom{v+1}{1}$. The $h$-vector of $A$
  is of course still $(1,4,10,..,6s+16,6s+19,6s+16,..,10,4,1)$ and the Hilbert
  polynomial of $X$ is $p_X(v)= (6s^2+44s+81)(v-s-3)$. This time we suppose $s
  \geq -2$, only avoiding the degenerate case. $A$ is unobstructed by
  Corollary~\ref{corH}. To find the dimension of $ \GradAlg^{H_A}(R)$, let
  $\eta(v):= \dim (I/I^2)_v$ and note $\eta(v)= \dim (I)_v$ for $v <6$. Then $
  \ \dim (N_B)_0 = 3 \eta(4)-4 \eta(3) =48$ by Remark~\ref{remextI}. Moreover,
  by \eqref{Ncodim2}, $ \ \dim (H_1^*)_s = 4 \cdot \dim B_{(3+s)} -\ \dim
  (N_B)_s$ where $\ \dim (N_B)_s = 3 \dim I_{(4+s)} - 4 \dim I_{(3+s)} + \dim
  R_s $ by Remark~\ref{remextI}. Hence {\small $$ \ \dim (H_1^*)_s = 4 \cdot
    \dim R_{(3+s)} - 3 \dim I_{(4+s)} - \dim R_s = 4 \binom{s+8}{5} -12
    \binom{s+6}{5} + 8 \binom{s+5}{5} = 2(s+4)^2(s+5).$$ }Moreover $_0\!\hom_B(
    S_2(H_1), K_B(-6)) = 6 \dim (K_B)_2 - 3 \dim S_2(K_B)_4 = 0$ by
    \eqref{S2(H*)}, \eqref{Kcodim2} and \eqref{S2Kcodim2}. Hence by
    Corollary~\ref{corH} and \eqref{Grad}, if $s>0$, then
  $$
  \ \dim_{(A)} \GradAlg^{H_A}(R)= \ \dim_{(X)} \Hilb^{p_X}(\pp^5)=
  2(s+4)^2(s+5)+47.
  $$
  and the quotients given by \eqref{H1(-s)} generate a generically smooth,
  irreducible component of $ \GradAlg^{H_A}(R)$. Now $\delta(K_B)_{-6-2s}= 0 $
  and $ _{-2s}\!\Ext_B^1(I/I^2,K_B(-6))=0$ for $s > -2$ by
  Remark~\ref{remMmainth}$(iii)$ while $(K_B)_{-6-2s}= 0$ for $s \geq -2$ by
  \eqref{Kcodim2}. To compute $ \delta(H_1)_{-s} $ and $\
  _{-s}\!\ext^1_B(I/I^2,H_1)$ we use \eqref{deltacodim2} and we get $\
  _{-s}\!\ext^1_B(I/I^2,H_1) = \dim R_{-s}$ and $ \delta(H_1)_{-s} = \ -
  _{-s}\!\ext^1_B(I/I^2,H_1)$ for $-1 \leq s \leq 0$ by Remark~\ref{remcodim}.
  By Corollary~\ref{corH}, $$ \ \dim_{(A)} \GradAlg(H_A)= 2(s+4)^2(s+5)+47 + \
  _{-s}\!\ext^1_B(I/I^2,H_1) \ \ {\rm for} \ \ -1 \leq s \leq 0,$$ and the
  codimension, $_{-s}\!\ext^1_B(I/I^2,H_1)$, of the $H_B$-stratum of $A$ is
  $1,6$ for $s=0,-1$ respectively. Furthermore, for $s=-2$, we compute
  $\delta(H_1)_2$ and $\delta(K_B)_{-2}$ exactly as described in
  Remark~\ref{remcodim} above and we get its values to be $-3$ and $3$
  respectively. Hence $ \dim_{(A)} \GradAlg^{H_A}(R)= 2(s+4)^2(s+5)+47 = 71$
  in this case.
  
  Finally we pay some extra attention to the case $s=0$ of AG curves of degree
  $d=81$ and genus $g=244$. In this case the stratum given by \eqref{H1(-s)}
  forms a 207 dimensional irreducible family contained in an irreducible
  component of $\ \GradAlg^{H_A}(R)$ or of $\ \Hilb^{p_X}(\pp^5)$ of dimension
  208. It is interesting to observe that we have exactly the same degree and
  genus as for the AG curve with s = 6 of Example~\ref{exKmainth} where a 207
  dimensional stratum in a 208 dimensional component was constructed by means
  of Theorem~\ref{Kmainth}!  Their resolutions, however, differ (see
  \eqref{H1A} for the AG curve of the first family where $ R(-4)^3$ appears as
  a repeated factor, and use linkage to see that $ R(-5)^3$ appears as a
  repeated factor in the resolution of the second family). We would like to
  pose the following two questions. Are these two strata subschemes of
  codimension one of the {\it same} irreducible component of $
  \GradAlg^{H_A}(R)$? Is the minimal resolution of the generic Gorenstein
  algebra of this component obtained from the resolution of \eqref{H1A} by
  deleting all repeated factors?
\end{example}

Now we apply Theorem~\ref{Mmainth} to $M=N_B$. By Theorem~\ref{mainth} and
Proposition~\ref{propnormal} we have the exact sequence
\begin{equation} \label{NB(-s)}
0 \rightarrow K_B(t-2s)\rightarrow N_B(-s) \rightarrow  I_{A/B} \rightarrow 0
\ 
\end{equation}
where $t=n+2$. Since we have not been able to verify the assumption $\ 
_t\!\Ext_B^2(S_2(N_B),K_B)=0$ of Theorem~\ref{Mmainth}A), even by increasing
$\depth_{I(Z)}B$ to its largest possible value, we must use the B)-part of the
Theorem. Here we get advantage of developing the concept ``unobstructed along
any graded deformation of $B$'', and we must suppose $s$ ``large enough''.
Indeed we get

\begin{corollary} \label{corNB} Let $B=R/I$ be a graded codimension two CM
  quotient of $R$, let $U= \Proj(B)-Z \hookrightarrow \pp^{n+1}$ be an l.c.i.
  and suppose $\depth_{I(Z)} B \geq 4$. If $A$ is given by a regular section
  of $\widetilde{N_B}^*(s)$ on $U$ % i.e gi \eqref{NB(-s)},
  and if $s > 2 \max n_{2,j}- \min n_{1,i}$ and $char(k) \neq 2$, then $A$ is
  unobstructed as a graded $R$-algebra, and the stratum of quotients given by
  \eqref{NB(-s)} around $(A)$ is open and irreducible (so $A$ is
  $H_B$-generic). Moreover $A$ is Gorenstein of codimension 4 in $R$, and
  $$
  \ \dim_{(A)} \GradAlg^{H_A}(R)= \ \dim_{(X)} \Hilb^{p}(\pp^{n+1})= \ \dim
  (N_B)_0 + \ \dim (I/I^2)_s - \ _0\!\hom(I/I^2,I/I^2) \ . $$ Letting
  $\eta(v):= \dim (I/I^2)_v$, we also have
  $$
  \ \dim_{(A)} \GradAlg^{H_A}(R)= \ \dim_{(X)} \Hilb^{p}(\pp^{n+1})= \eta(s)+
  \sum_{j=1}^{\mu -1} \eta( n_{2,j})-\sum_{i=1}^{\mu} \eta( n_{1,i}) \ . $$
  Furthermore if $(B') \in \GradAlg^{H_B}(R)$ satisfies the same assumptions
  as $B$ above and defines $A'$ as $B$ defined $A$, then the closures in $
  \GradAlg^{H_A}(R)$ of the stratum of quotients given by \eqref{NB(-s)}
  around $(A)$ and the corresponding stratum around $(A')$ coincide, i.e. they
  form the same irreducible component of $ \GradAlg^{H_A}(R)$.
\end{corollary}

\begin{proof} Firstly to see that $(N_B,B)$ is unobstructed along any graded
  deformation of $B$, it suffices by Proposition~\ref{unobstrp1}$(iii)$ to
  show that $\ _0\!\Ext_B^1(N_B,N_B)=0$ since $N_{B_T}:= \Hom(I_T/I_T^2,B_T)$
  ($B_T :=R_T/I_T$ a graded deformation of $B$ to an Artinian $T$) is a graded
  deformation of $N_B$ to $B_T$ by \cite{JS},\! Prop.\! A1 and $\
  _0\!\Ext_B^1(I/I^2,B)=0$. Since $\ \depth_{I(Z)} B \geq 4 \ $ implies
  $\depth_{I(Z)} N_B \geq 4 \ $ and $ \depth_{I(Z)} I/I^2 \geq 3 \ $ we get by
  \eqref{twoo} and Proposition~\ref{extI}, $$\ _0\!\Ext_B^1(N_B,N_B)\ \simeq \
  \HH_*^1(U,\sH om_{\sO_U}(\widetilde {I/I^2}^* \arrowvert_U, \widetilde
  {I/I^2}^* \arrowvert_U)) \ \simeq \ \Ext_B^1(I/I^2,I/I^2)=0$$ and we get
  what we want. Similarly $\ \Ext_B^1(N_B,B)\ \simeq \HH_*^1(U,\sH
  om_{\sO_U}(\widetilde {I/I^2}^*, \widetilde {B})) \simeq
  \HH_*^1(U,\widetilde {I/I^2}) =0$.
  
  Secondly to show the remaining assumption of Theorem~\ref{Mmainth}B), we use
  Remark~\ref{remMmainth}$(ii)$. Since $\oplus R( n_{2,j}- n_{1,i})
  \rightarrow N_B$ is surjective by Remark~\ref{remextI}, it follows that
  $(N_B)_v=0$ for $v < \min n_{1,i} - \max n_{2,j}$, and we conclude by
  Remark~\ref{remMmainth}$(ii)$. Note that also the dimension formula above
  follows from Remark~\ref{remMmainth} and \eqref{M,M} (see the text before
  \eqref{M,M} for $ _0\!\hom(I/I^2,I/I^2) =\ _0\!\hom(N_B,N_B)= \
  _0\!\hom(S_2(N_B),K_B(t))+1 $). Finally the proof of the irreducibility is
  trivial because the stratum is the image of an irreducible set. Indeed $
  \GradAlg(H_B,H_A)$ is smooth at $(B \rightarrow A)$ by the first conclusion
  of Theorem~\ref{unobstr}, cf. proofs of Theorem~\ref{Mmainth}B)
  Theorem~\ref{gorgenth}. The proof of the uniqueness (i.e. that the strata
  above coincide, up to closure) is essentially the same as for \cite{K03},
  Prop.\!\ 23$(i)$ and Thm.\!\ 24 (see the two first lines of the proof of
  \cite{K03}, Thm.\!\ 24) because the open subscheme of $ \GradAlg^{H_B}(R)$
  consisting of CM quotients is irreducible. For the connection with the
  Hilbert scheme, see \eqref{Grad}.
\end{proof}

\begin{remark} \label{remcorNB}
  Let $\depth_{I(Z)} B \geq 4$ and $char(k) \neq 2$. Using \eqref{twoo} (see
  the proof above), we get
  $$\ _0\!\Ext_B^2(N_B,N_B)\ \simeq \ \Ext_{\sO_U}^2(\widetilde {I/I^2}
  \arrowvert_U,\widetilde {I/I^2} \arrowvert_U) \simeq \ _0\!\Hom_B(I/I^2,
  \HH^3_{I(Z)}(I/I^2))\ ,$$
  and $ \HH^3_{I(Z)}(I/I^2) \simeq
  \HH^4_{I(Z)}(I^2)$. Sometimes we can use this connection to prove the
  vanishing of $\ _0\!\Ext_B^2(N_B,N_B)$ in which case
  Theorem~\ref{Mmainth}A) applies to $M=N_B$ provided we can show
  $_{-s}\!\Ext_B^2(I/I^2,N_B)=0$. To show this vanishing, we play on
  $N_B \simeq I/I^2 \otimes K_B(n+2)$ and we get $
 \  \Ext_B^2(I/I^2,N_B) % \simeq \ \Ext_B^2((I/I^2)^* \otimes I/I^2, K_B(n+2))
 \ \simeq \ \Ext_B^2( \Hom(I/I^2,I/I^2), K_B(n+2))=0$ from the fact that the
 codepth of $\Hom(I/I^2,I/I^2) $ is at most one by Remark~\ref{remextI}.  The
 benefit of using Theorem~\ref{Mmainth}A) (instead of
 Theorem~\ref{Mmainth}B) which requires $_{-s}\!\Ext_B^1(I/I^2,N_B)=0$) is
 that we don't need to assume $s > 2 \max n_{2,j}- \min n_{1,i}$. Hence if we
 in the case $s \leq 2 \max n_{2,j}- \min n_{1,i}$ suppose $ _0\!\Hom_B(I/I^2,
 \HH^4_{I(Z)}(I^2)) =0$, we get that $A$ is an unobstructed graded Gorenstein
 quotient of codimension 4 in $R$, and $$ \ \dim_{(A)} \GradAlg^{H_A}(R) =
 \epsilon + \delta + \dim(K_B)_{t-2s}$$ where $\delta := \ \delta(K_B)_{t-2s} -
 \delta(N_B)_{-s}$ and $\epsilon:= \ \dim (N_B)_0 + \ \dim (I/I^2)_s - \ 
 _0\!\hom(I/I^2,I/I^2)$, i.e. $\epsilon$ is equal to the expression of
 $\eta's$ in Corollary~\ref{corNB}.  Continuing this argument (for this final
 statement we omit the details of the proof) we may even show that we can skip
 $s > 2 \max n_{2,j}- \min n_{1,i}$ in Corollary~\ref{corNB} and at least get
 inequalities
 $$
 \epsilon + \delta + \dim(K_B)_{t-2s} - \ _0\!\ext_B^2(N_B,N_B) \leq \ 
 \dim_{(A)} \GradAlg^{H_A}(R) \leq \epsilon + \delta + \dim(K_B)_{t-2s} \ .$$
 Moreover, if the inequality to the right is an equality, then $A$ is
 unobstructed.
 \end{remark}

 \begin{example} (Arithmetically Gorenstein curves $\Proj(A)$ in $ \pp^5$,
   obtained by \eqref{NB(-s)}.) Let $R$ be a polynomial $k$-algebra in 6
   variables, let $B=R/I$ be a codimension two quotient with minimal
   resolution
\begin{equation} \label{I}
 0 \rightarrow R(-3)^2 \rightarrow  R(-2)^3  \rightarrow R
 \rightarrow B \rightarrow 0  \ . 
\end{equation}
and suppose $Y= \Proj(B)$ is an l.c.i in $\pp^5$. Let $A$ be given by a regular
section of $\widetilde {I/I^2}(s)$.  The Hilbert polynomial/function of $Y$ is
 $p_Y(v)= H_B(v)= 3 \binom{v+2}{3} + \binom{v+2}{2}$ for
$v \geq 0$  by \eqref{I} (e.g. by ``integrating'' $3v+1$ two times). Thanks to
\eqref{NBcodim2} and \eqref{Kcodim2} and the mapping cone construction applied
to \eqref{NB(-s)} we get the following resolution (obviously minimal for $s >
4$) of the Gorenstein algebra $A=R/I_A$ of Corollary~\ref{corNB},
\begin{equation} \label{NBA}
 \begin {aligned}
  & 0 \rightarrow R(-2s) \rightarrow  R(2-2s)^3 \oplus R(-1-s)^6 \\ 
 & \rightarrow R(3-2s)^2 \oplus R(-s)^{12} \oplus R(-3)^2 \rightarrow
 R(1-s)^6 \oplus R(-2)^3  \rightarrow R \rightarrow A 
 \rightarrow 0 \ \ . 
 \end{aligned}
\end{equation} 
Let $h^i(\sO_X(v))= \dim \HH^i(X,\sO_X(v))$. Since $K_A \simeq A(2s-6)$, the
Hilbert polynomial of $X$ is of the form $p_X(v)= h^0(\sO_X(v))-h^1(\sO_X(v))=
d(v-s+3)$ because $h^1(\sO_X(s-3))= h^0(\sO_X(s-3))$. Moreover looking to the
resolution of $I_A$ we see that {\small $$p_X(s-2)=
  h^0(\sO_X(s-2))-h^0(\sO_X(s-4))=
  h^0(\sO_Y(s-2))-h^0(\sO_Y(s-4))=3s^2-10s+9.$$ } So $X = \Proj(A)$ is an AG
curve of degree $d = 3s^2-10s+9$ and arithmetic genus $g =1+d(s-3)$ in $
\pp^5$. If $s > 2 \max n_{2,j}- \min n_{1,i}=4$, then Corollary~\ref{corNB}
applies. Letting $\eta(v):= \dim (I/I^2)_v$, we get $ \sum
\eta(n_{2,j})-\sum\eta( n_{1,i}) = 2 \dim I_3- 3 \dim I_2 = 23$. Calculating
$\eta(s)$ and using \eqref{I2codim2} which implies
\begin{equation} \label{I2}
 0 \rightarrow R(-6) \rightarrow R(-5)^6 \rightarrow  R(-4)^6  \rightarrow I^2
\rightarrow 0  \, 
\end{equation}
we get $\eta(s)=(s+1)(s-1)^2$. Hence, if $s \geq 5$, then $A$ is unobstructed
and
$$
\ \dim_{(A)} \GradAlg^{H_A}(R)= \ \dim_{(X)} \Hilb^{p_X}(\pp^{5})= \ \eta(s)
+ 2\eta(3)- 3\eta(2) = (s+1)(s-1)^2+23 \ .$$
  
Finally we discuss the cases $3 \leq s \leq 4$. By \eqref{I2}, we get an
injection $ \HH_{\mathfrak m}^4(I^2) \rightarrow \HH_{\mathfrak m}^6(R(-6))$.
By Remark~\ref{remcorNB}, $\ _0\!\Ext_B^2(N_B,N_B)=0$ and it follows that $A$
is unobstructed and $ \ \dim_{(A)} \GradAlg^{H_A}(R)= (s+1)(s-1)^2+23 + \delta
$ where $\delta := \ \dim(K_B)_{6-2s} + \delta(K_B)_{6-2s} -
\delta(N_B)_{-s}$. We have $(K_B)_{6-2s} = 0$ by \eqref{Kcodim2}.  Moreover
applying $\Hom(-,K_B(6))$ to \eqref{es1}, we get $\ -\delta(K_B)_{6-2s}=\ 
_{-2s}\!\hom(H_1,K_B(6))$ for $s \geq 3$. Since $H_1$ has rank one by
\eqref{es1} we get $ H_1 \simeq \Hom(H_0,K_B(6- \sum n_{1,i})) \simeq K_B $ by
the last conclusion of Theorem~\ref{mainth}. Hence $\delta(K_B)_{6-2s}=\ -\dim
B_{(6-2s)} $, i.e.  $\delta(K_B)_{6-2s} = -1,0$ for $s=3,4$ respectively. It
remains to compute $\delta(N_B)_{-s} $ for which we use $H_1 \simeq K_B$ and
the following exact sequence
 \begin{equation} \label{deltaNBcodim2}
 \small
  0 \rightarrow  \Hom(I/I^2,N_B) \rightarrow
 \Hom(I/I^2,\oplus B(n_{1,i})) \rightarrow \Hom(I/I^2,K_B^*) \rightarrow
 \Ext^1_B(I/I^2,N_B) \rightarrow 0  
\end{equation}
which we get by applying $\Hom(I/I^2,-)$ to \eqref{Ncodim2}. Since $
\Hom(I/I^2,K_B^*) \simeq \Hom(I/I^2 \otimes K_B(6),B(6)) \simeq I/I^2(6)$, we
get $\delta(N_B)_{-s} = 2,-3$ for $ s = 3,4$ respectively. By
Remark~\ref{remcorNB}, $ \ \dim_{(A)} \GradAlg^{H_A}(R)= 36,71$ for $s=3,4$
respectively. Note that this result confirms the well known formula
$$6d+2(1-g)\leq \dim_{(A)}\GradAlg^{H_A}(R) \leq 6d+2(1-g)+h^1(\sN_X)$$
because for $s=3$ (resp. $s=4$) the curve has degree $d = 6$ and genus $g=0$
(resp. $d = 17$ and genus $g=18$, and one may verify $h^1(\sN_X)=3$ by other
methods).
\end{example}

We will finish the section by looking to the rank $r = 3$ case, i.e. we apply
Theorem~\ref{Mmainth3}A) to the $B$-module $M=H_1$ provided the number of
minimal generators of $I$ is $\mu(I)=5$. Note that \eqref{M(-s)3} translates
to $t = n+2- \sum n_{1,i}$, $ H_2 = H_1^{\vee}(t) \simeq \Hom(H_1,K_B(t)) $
and
\begin{equation} \label{H1(-s)3}
0 \rightarrow K_B(t-3s)\rightarrow H_2(-2s) \rightarrow H_1(-s) \rightarrow
I_{A/B} \rightarrow 0 \ .
\end{equation}
\begin{corollary} \label{corH3} Let $B=R/I$ be a graded codimension two CM
  quotient of $R$, let $U= \Proj(B)-Z \hookrightarrow \pp^{n+1}$ be an l.c.i.
  and suppose $\mu(I)=5$, $char(k) \neq 2$ and $\depth_{I(Z)}B \geq 3$. If $A$
  is defined by a regular section of $\widetilde{H_1}^*(s)$ on $U$, i.e. given
  by \eqref{H1(-s)3}, then $A$ is unobstructed as a graded $R$-algebra (indeed
  $_0\!\HH^2(R,A,A)=0$), $A$ is Gorenstein of codimension 5 in $R$, and $ \
  \dim_{(A)} \GradAlg^{H_A}(R)= $
  $$
  \ \dim (N_B)_0 + \ \dim (H_1^*)_s + _{-s}\!\hom_B(S_2(H_1), K_B(t)) - _0\!
  \hom_B(H_1, H_1)- \dim(K_B)_{t-3s} - \delta \ ,$$ where $ \delta:=
  \delta(H_1)_{-s} + \delta(K_B)_{t-3s} - \delta(H_2)_{-2s}$. Moreover if
  $char(k)=0$ and $(B \rightarrow A)$ is general with respect to $\
  _{0}\!\hom_R(I_B,I_{A/B}) $, then the codimension of the $H_B$-stratum of
  $A$ at $(B \rightarrow A)$ is $$ \ _{-s}\!\ext^1_B(I_B/I_B^2,H_1) - \ \dim
  (\im \beta)\
  $$ where $\beta$ is the homomorphism $\ _{-2s}\!\Ext_B^1(I/I^2,H_2)
  \rightarrow \ _{-s}\!\Ext_B^1(I/I^2,H_1)$ induced by \eqref{H1(-s)3}. If in
  addition $\ s > \max n_{2,j}-\min n_{2,j} \ $, then $A$ is $H_B$-generic, $\
  \delta = 0$ and $ \dim(K_B)_{t-3s}=0$.
\end{corollary}

\begin{proof} This is a corollary to Theorem~\ref{Mmainth3}A). Indeed, thanks
  to Proposition~\ref{propH} we only need to show $ \ _0\!\Ext_B^2(H_2(-t)
  \otimes H_1,K_B)= 0 \ $ and $ \ _{-2s}\!\Ext^i_B(I/I^2,H_2)= 0 \ $ for
  $2 \leq i \leq 3$ and the ``if in addition''-statement. Using \eqref{es1} we
  get $ \ _{-2s}\!\Ext^i_B(I/I^2,H_2) \simeq \ 
  _{0}\!\Ext^{i-1}_B(H_1,H_2(-2s))$.  Since $\depth_{I(Z)}B \geq 3$, we have
  $$ \ _{0}\!\Ext^{1}_B(H_1,H_2(-2s)) \simeq \Ext_{\sO_U}^{1}(\widetilde {H_1}
  \arrowvert_U,\widetilde {H_1^* \otimes K_B}(t-2s) \arrowvert_U) \simeq
  \Ext_{\sO_U}^{1}(\widetilde {H_1 \otimes H_1} \arrowvert_U,\widetilde
  {K_B}(t-2s) \arrowvert_U) $$ by \eqref{twoo} which vanishes because $0
  \rightarrow \ \widetilde H_2 \arrowvert_U \simeq \wedge^2 \widetilde H_1
  \arrowvert_U \rightarrow \widetilde { H_1 \otimes H_1} \arrowvert_U
  \rightarrow S_2(H_1)\arrowvert_U \rightarrow 0 $ is exact and $H_2$ and
  $S_2(H_1)$ are maximal CM modules (this argument also shows $ \
  _{0}\!\Ext^{2}_B(H_1,H_2(-2s)) = 0$ provided $\depth_{I(Z)}B \geq 4$).
    
  To show $ \ _{0}\!\Ext^{2}_B(H_1,H_2(-2s)) = 0$ under the assumption
  $\depth_{I(Z)}B \geq 3$, we use \eqref{K*codim2}. Hence it suffices to show
  the $ \ _v\! \Ext^{1}_B(K_B^*,H_2) = 0$ for every $v$. By \eqref{twoo} this
  $\Ext$-group is isomorphic to %$ \ _v\! \Ext^{1}_B(K_B^*,H_2) \simeq $
  {\small $$ \Ext_{\sO_U}^{1}(\widetilde {K_B}^* \arrowvert_U, \widetilde
    {H_1^* \otimes K_B}(t+v) \arrowvert_U) \simeq \Ext_{\sO_U}^{1}(\widetilde
    {K_B^* \otimes H_1}\arrowvert_U ,\widetilde {K_B}(t+v) \arrowvert_U)
    \simeq \ _v\! \Ext^{1}_B(\Hom_B(K_B, H_1),K_B(t))$$ }which vanishes since
    $\Hom_B(K_B, H_1)$ is a maximal CM module. Indeed by using
    \eqref{K*codim2}, we showed in \cite{KR02} the exactness of
\begin{equation} \label{H2codim2}
  0 \rightarrow \Hom_B(K_B(n+2), H_1) \rightarrow \wedge^2(\oplus
  B(-n_{2,j})) \rightarrow H_2 \rightarrow 0 \ 
\end{equation}
  (for that sequence it suffices to have   $\depth_{I(Z)}B \geq 2$), from
  which we get that $\Hom_B(K_B, H_1)$ is a maximal CM $B$-module.
  
  To prove $ \ _t\!\Ext_B^2(H_2 \otimes H_1,K_B)= 0 \ $, we use again
  \eqref{K*codim2} to get $0 \rightarrow \widetilde {K_B(n+2)}^* \otimes
  \widetilde { H_2} \arrowvert_U \rightarrow (\oplus \widetilde B(-n_{2,j}))
  \otimes \widetilde H_2 \arrowvert_U \rightarrow \widetilde H_2 \otimes
  \widetilde H_1 \arrowvert_U \rightarrow 0 \ $. Applying
  $\HH^0_*(U,{-})$ to it, we get the exact sequence
  $$0 \rightarrow \Hom_B(K_B(n+2), H_2) \rightarrow (\oplus B(-n_{2,j}))
  \otimes H_2 \rightarrow H_2 \otimes H_1 / \tau \rightarrow 0 \ $$
  where
  $\tau:= \HH^0_{I(Z)}( H_2 \otimes H_1)$. Invoking also \eqref{twoo} it
  follows that $ \ _{v-n-2}\! \Ext_B^2(H_2 \otimes H_1/ \tau,K_B) \simeq \ 
  _v\! \Ext^{1}_B(\Hom_B(K_B, H_2),K_B) \simeq \Ext_{\sO_U}^{1}(\widetilde
  {K_B}^*\otimes \widetilde {H_1^* \otimes K_B}(t) \arrowvert_U,\widetilde {
    K_B}(v) \arrowvert_U) \simeq \ _v\! \Ext^{1}_B(H_1^*(t),K_B)$, and that the
  latter vanishes since $H_1^*$ is a maximal CM $B$-module. Indeed due to
  \eqref{Ncodim2}, $H_1^*$ has codepth at most one, and applying $\Hom(-,K_B)$
  to \eqref{Ncodim2}, we will see $ \Ext^{1}_B(H_1^*,K_B)=0$ as well because
  $\Hom(N_B,K_B) \simeq N_B(-n-2) \simeq I/I^2 \otimes K_B$. Then we conclude
  by the exact sequence $$ \Ext_B^{2}( H_2 \otimes H_1/ \tau,K_B) \rightarrow
  \Ext_B^{2}( H_2 \otimes H_1 ,K_B ) \rightarrow \Ext_B^{2}( \tau, K_B) $$
  because $ \Ext_B^{2}( \tau, K_B) \hookrightarrow \Ext_{\sO_U}^{2}(\widetilde
  {\tau} \arrowvert_U, \widetilde { K_B}\arrowvert_U) $ is injective by
  \eqref{twoo} and $\widetilde {\tau} \arrowvert_U=0$.
  
  Finally if $s > \max n_{2,j} - \min n_{2,j}$ we get $
  _{-s}\!\Hom_R(I,H_1)=\ _{-s}\!\Ext_R^1(I,H_1)=0 \ $ and hence $
  \delta(H_1)_{-s} = 0$ by Remark~\ref{remMmainth3}$(ii)$ and the ``left''
  short exact sequence deduced from \eqref{K*codim2}.  Using \eqref{H2codim2}
  instead of \eqref{K*codim2} we similarly get $ \delta(H_2)_{-2s} = 0$ as
  well as $ \delta(K_B)_{t-3s} = 0$ and we conclude easily.
 \end{proof}

\begin{remark} \label{remcorH3}
  Due to Remark~\ref {rem2Mmainth3} the quotients $A$ of Corollary~\ref{corH}
  and Corollary~\ref{corH3} are strongly unobstructed in the sense $\ 
  \HH^2(R,A,A)=0$. This follows from the proofs of the corollaries because the
  proofs of the vanishing of the $\ _v\! \Ext_B^i$-groups involved are easily
  extended to the vanishing of the corresponding $\ \Ext_B^i$-groups. As a
  consequence of this, look to the scheme $ZGor(H)$ parametrizing {\it not
    necessarily graded} Gorenstein quotients $R \rightarrow A$ with Hilbert
  function $H$ (\cite{IK}, p.126).  Since $_v\!\HH^2(R,A,A)=0$ for $v \geq 0$,
  we can use \cite{K98}, Thm.\!\ 1.10 to prove the smoothness of $ZGor(H)$ at
  $(A)$, and one may also find the dimension of the scheme $ZGor(H)$ at $(A)$
  by a formula analogous to the dimension formulas of Corollary~\ref{corH} and
  Corollary~\ref{corH3}, cf. \cite{K03}, Prop.\!\ 29.
\end{remark}

\begin{example}
  Let $R$ be a polynomial ring in five variables, and let $B=R/I$ be
  a CM-quotient with $\mu(I) = 5$ and with linear resolution
(i.e. all $n_{1,i}$ of \eqref{Icodim2} are  $n_{1,i}=4$). Then $Y= \Proj(B)$
  is a surface with Hilbert polynomial $p_Y(v)= 5v^2-5v+5$. By \cite{AH}, 
 \begin{equation} \label{wedgeH2}
0 \rightarrow \wedge^3 G_2 \rightarrow \wedge^3 G_1 \rightarrow  \wedge^2 G_2
\rightarrow  H_2 \rightarrow 0
\end{equation} 
is exact, and the mapping cone construction applied to \eqref{H1(-s)3} lead to
the following resolution of the Gorenstein algebra $A=R/I_A$ of
Corollary~\ref{corH3}A),
\begin{equation} \label{H2A}
 \begin {aligned}
  & 0 \rightarrow R(-20-3s) \rightarrow  R(-16-3s)^5 \oplus R(-15-2s)^8
 \rightarrow  R(-15-3s)^4 \oplus R(-12-2s)^{10}  \oplus \\ & R(-10-s)^6  
  \rightarrow R(-10-2s)^6 \oplus R(-8-s)^6 \oplus R(-5)^4 \rightarrow
 R(-5-s)^4 \oplus R(-4)^5  \rightarrow I_A 
 \rightarrow 0 .
 \end{aligned}
\end{equation} 
By \eqref{H1(-s)3} or \eqref{H2A} $A$ is Artinian of socle degree $3s+15$ and
with symmetric $h$-vector $(1,...,1)$ given by $H_A(v)=p_Y(v)$ for $1 \leq v
\leq s+4$, $H_A(v)=p_Y(v)-4$ for $v = s+5$ and
$$H_A(v)=p_Y(v)-15(v-s-4)^2+35(v-s-4)-30 \ \ {\rm for} \ \ s+5 < v < 2(s+5).$$
So if $s=-3,-2,...$, then the $h$-vector of $A$ is $(1,5,11,15,11,5,1)$,
$(1,5,15,31,45,45,31,15,5,1),...$ respectively. Suppose $s \geq -3$ and $Y$ an
l.c.i.. Then $A$ is unobstructed by Corollary~\ref{corH3}. By
Remark~\ref{remextI}, $\ \dim (N_B)_s = 4\dim I_{s+5}-5\dim I_{s+4}+\dim R_s$,
so $\ \dim (N_B)_0=60$. Moreover {\small $$ \ \dim (H_1^*)_s = 5 \dim B_{s+4}-
  \dim (N_B)_s = 5 \binom{s+8}{4} -4 \dim I_{s+5} - \binom{s+4}{4} =
  15s^2+125s+265$$ }by \eqref{Ncodim2}. Applying $\Hom(-,H_1)$ to \eqref{es1},
  and we get
 \begin{equation} \label{deltacodim3}
 \small
  0 \rightarrow  \Hom(I/I^2,H_1) \rightarrow
 \Hom(\oplus B(-n_{1,i}),H_1) \rightarrow \Hom(H_1,H_1) \rightarrow
 \Ext^1_B(I/I^2,H_1) \rightarrow 0  
\end{equation}
Hence $\ _0\! \hom(H_1,H_1) = 1$  by Remark~\ref{remcodim} and $
_0\!\hom_B( S_2(H_1), K_B(-15-s)) = 0$ by \eqref{S2(H*)} and \eqref{Kcodim2}.
Finally $ \dim(K_B)_{-15-3s}=0$ and the sum $
 \ \dim (H_1^*)_s + _{-s}\!\hom_B(S_2(H_1), K_B(t)) - \
  _0\! \hom_B(H_1, H_1)- \dim(K_B)_{t-3s}$ of  Corollary~\ref{corH3} is equal
  to $15s^2+125s+264$. We get
  $$
  \ \dim_{(A)} \GradAlg^{H_A}(R)= \ \dim_{(A)} \PGor(H_A)= 15s^2+125s+324-
  \delta.$$ Moreover if $s > 0$, then $ \delta =0$ and $A$ is $H_B$-generic.
  Looking to $s \leq 0$, we note that if $\ 2s > \max n_{2,j} - 2\min
  n_{2,j}$, then $ _{-2s}\!\Hom_R(I,H_2)=\ _{-2s}\!\Ext_R^1(I,H_2)=0 \ $
  and hence $ \delta(H_2)_{-2s} = 0$ by Remark~\ref{remMmainth3} and
  \eqref{H2codim2}. Combining with Remark~\ref{remcodim}, we get the
  codimension of the $H_B$-stratum of Corollary~\ref{corH3} to be $\ -\delta =
  - \delta(H_1)_{-s} = _{-s}\!\ext^1_B(I/I^2,H_1) = \dim R_{-s}$ for $-2 \leq
  s \leq 0$. In particular $$ \ \dim_{(A)} \GradAlg^{H_A}(R)= 15s^2+125s+324 +
  \binom{-s+4}{4} \ \ {\rm for} \ \ -2 \leq s \leq 0.$$ In the final case
  $s=-3$ (where we skip a few details), one may see $ \delta(H_2)_{-2s} = 5
  \dim H_2(4)_6- \ _6\!\hom_B( S_2(H_1), K_B(-15)) = 0 $ (from the sequence we
  get by applying $\Hom(-,H_2)$ to \eqref{es1}), and $$ \delta(H_1)_{-s} = 5
  \dim N_B(-4)_3- \ _3\!\hom_B(I/I^2,I/I^2)= 5 \cdot 20 - (4 \dim I_8^2+\dim
  R_3) = 5$$ by Remark~\ref{remcodim}, and we get $ \ \dim_{(A)}
  \GradAlg^{H_A}(R)= 15s^2+125s+324 -5 = 79$.
\end{example}

\section{Appendix: Deformations of quotients of zerosections}
Let $X= \Proj(A)$ be a subscheme of $Y= \Proj(B) \subset \pp =
\Proj(R)$ defined as the degeneracy locus of a regular section of some sheaf
$\sM$ supported on $Y$. In this appendix we prove a quite general result,
Theorem~\ref{unobstr}, concerning the unobstructedness and the
``family-dimension'' of a quotient $A$ obtained from $B$ as the homogeneous
coordinate ring of a zerosection as above, in which we neither assume $B$ to
be Cohen-Macaulay, nor $A$ to be Gorenstein. This leads to a main result of
this paper (Theorem~\ref{gorgenth}). We have included a version of
Theorem~\ref{unobstr} for the Hilbert scheme (Corollary~\ref{corunobstr}).
This result essentially generalizes Thm.\!\ 9.4 of \cite{KMMNP}, which treats
the case where $\sM$ is locally free of rank $r=1$, to higher ranks.
Unfortunately our methods often lead to assumptions on $\depth_{I(Z)}B$ which
imply that $A$ is non-Artinian. We have, however, succeeded in proving
Theorem~\ref{gorgenth} also for an Artinian Gorenstein algebra $A$. Below $B
= R/I_B \rightarrow A = R/I_A$ is a graded surjection with kernel
$I_{A/B}$ and $M^* = \Hom_B(M,B)$.

\begin{theorem}\label{unobstr}
  Let $r \geq 1$ and $s$ be integers. Let $B$ be a graded quotient of a
  finitely generated polynomial $k$-algebra $R$.  Let $M$ be a finitely
  generated graded $B$-module, let $Y:= \Proj(B)$ and $U=Y-Z$ be an open subset
  of $\ \Proj(B)$ such that $\depth_{I(Z)}B \geq 2$ and such that $\widetilde
  {M}\arrowvert_U$ is locally free of rank $r$.  Let $
  M_i=\HH^0_*(U,\wedge^i\widetilde M)$ for $0 \leq i \leq r$, and let $\sigma
  \in \HH^0(U,\widetilde M^*(s))$ be a regular section on $U$. Let $X=\Proj(A)$
  be the zero locus of $\sigma$ defined by $A:= \coker (\HH^0_*(U,\widetilde
  M(-s)) \stackrel{\sigma}{\rightarrow}B)$. Let $K_1=\ker
  {\sigma}(s)$ and suppose
 \begin{itemize}
  \item [(i)] $\HH^i_*(U,\widetilde M_{i+1})=0 \ \ {\rm for}\ \
    1\leq i \leq r-1$ 
    % (or  $\HH^1_*(U,\widetilde M_2)=0$ and
    %$\HH^i_*(U,\widetilde M_{i+2})=0 \ \ {\rm for}\ \ 
    %1\leq i \leq r-2$)
 % then $\dim C = \dim S -r$
 % and $\aaa$ is equidimensional and $S_1$. Hence $\aaa$ is the
 % homogeneous coordinate ring of $C$.
 %  \item [(ii)] $B$ is unobstructed as a graded $R$-algebra
   \item [(ii)] $_0\!\Ext^2_B(M_1,K_1)=0$, and
 \item [(iii)] $(M_1(-s),\sigma)$ is unobstructed %(resp. unobstructed
   along any graded deformation of $B$.
 \end{itemize}
 Then the first projection $\GradAlg(H_B,H_A) \rightarrow \GradAlg^{H_B}(R)$ is
 smooth at $(B \rightarrow A)$ and, \\[2mm]
\indent $\ \ \ \dim_{(B \rightarrow A)}
 \GradAlg(H_B,H_A)= \ \dim_{(B)}\GradAlg^{H_B}(R)+\ _0hom_B(I_{A/B},A).$ \\[2mm]
 Moreover, $\depth_{\mathfrak{m}}A \geq 1 $. If, in addition,
 \begin{itemize}
 \item [(iv)] {\rm either} $_0\!\Ext^1_B(I_B/I_B^2,I_{A/B})=0$ and $(I_B)_{\wp}$ is
   syzygetic for any graded prime $\wp$ of $Ass(I_{A/B})$, {\rm or}
   $_0\!\Ext^1_R(I_B,I_{A/B})=0$, 
 \end{itemize}
 then $\ \ \ _0\hom_R(I_B,B)-\ \dim_{(B)}\GradAlg^{H_B}(R) =\ _0\hom_R(I_A,A)-
 \dim_{(A)} \GradAlg^{H_A}(R)  $,  \ $A$ is $H_B$-generic, and
 $$
 \dim_{(A)} \GradAlg^{H_A}(R)= \dim_{(B)}\GradAlg^{H_B}(R)+\ 
 _0\hom_B(I_{A/B},A)-\ _0\hom_R(I_B,I_{A/B}).$$
 Moreover $A$ is unobstructed
 as a graded $R$-algebra if and only if $B$ is unobstructed as a graded
 $R$-algebra.
\end{theorem}

Using the arguments appearing in \eqref{incl} below, we will see that the
condition (i) of Theorem~\ref{unobstr} implies $\HH^0_{I(Z)}(A) = 0$. In
particular $A$ is the homogeneous coordinate ring of $X$.  If we in addition
suppose

\ \ \ $(v)$ \ \ \ \ \ \ \ \ \ \ \ \ \ \ \ \ \ \ \ \ \ \ \
$\HH^i_*(U,\widetilde
M_{i})=0$ \ {\rm  for} $ \ 1\leq i \leq r$ \ , \\[2mm]
we can mainly argue as in \eqref{incl} (or as in the proof of \cite{KP2},
Lemma 12) to see that $\HH^1_{I(Z)}(A) = 0$, i.e. that $\depth_{I(Z)}A \geq
2$. Hence $\GradAlg^H(R) \simeq \Hilb^{ p}(\pp)$ at $(X \subset \pp)$ by
\eqref{Grad}. 

Let $X \subset Y$ be closed subschemes of $\pp$ of Hilbert polynomials $p_X$
and $p_Y$ respectively, let $\DD(p_X,p_Y)$ be the Hilbert-flag scheme
parametrizing all such ``pairs'' of closed subschemes and let $p_1:\DD(p_X,p_Y)
\rightarrow \Hilb^{p_X}(\pp)$ be the projection induced by $p_1((X' \subset
Y'))= (X')$. $X$ is called {\it $p_Y$-generic} if there is an {\em open}
subset $U_X$ of $ \Hilb^{ p_X}(\pp)$ such that $(X) \in U_X \subset
p_1(\DD(p_X,p_Y))$.

\begin{corollary}\label{corunobstr}
  In addition to the notations and assumptions of Theorem~\ref{unobstr}
  $(i)$-$(iv)$, suppose $(v)$. Then $X$ is $p_Y$-generic, and
  $$
  \ \dim_{(X)} \Hilb^{p_X}(\pp)= \dim_{(Y)}
  \Hilb^{p_Y}(\pp)+\hom_{\sO_Y}(\sI_{X/Y},\sO_X)
  -\hom_{\sO_{\pp}}(\sI_{Y},\sI_{X/Y}).$$
  Moreover $X$ is unobstructed if and only
  if $\ Y$ is unobstructed.
\end{corollary}

\begin{remark} \label{remunobstr}
  Since $\depth_{I(Z)}B \geq 2$, the modules $M_i$ of Theorem~\ref{unobstr}
  satisfy $M_i= (\wedge^i M)^{**}$ (cf. \cite{KP2}, Remark 8). In particular
  $\widetilde M_i= (\wedge^i\widetilde M)^{**}$ and clearly $\widetilde M_i
  \arrowvert_U = (\wedge^i\widetilde M) \arrowvert_U$.  Moreover note that
  Theorem~\ref{unobstr}$(i)$ holds if $$\depth_{I(Z)}M_i \geq i+1 \ \ {\rm
    for} \ \ 2 \leq i \leq r.$$
  %cf.  \eqref{ONE}. 
  This holds in particular if each $M_i$ is a maximal CM
  $B$-module and $\depth_{I(Z)}B \geq r+1$. In this case it follows from
  \cite{KP2}, Prop.\!\ 6 that $A$ is equidimensional and satisfies Serre's
  condition $S_1$ provided $B$ is Cohen-Macaulay. Indeed since a regular
  section by definition leads to an exact Koszul resolution of $\widetilde A$
  on $U$, we get that $U \cap X$ is equidimensional and without embedded
  components, and \cite{KP2}, Prop.\!\ 6 applies.  In the same way $(v)$
  above holds if each $M_i$ is maximally Cohen-Macaulay and $\depth_{I(Z)}B
  \geq r+2$, in which case $A$ satisfies Serre's condition $S_2$ if $B$ is
  Cohen-Macaulay.
\end{remark}

\begin{proof}
  Once we have proved the smoothness of the first projection
  $$q:\GradAlg(H_B,H_A) \rightarrow \GradAlg^{H_B}(R) \ \ {\rm at} \ \ (B
  \rightarrow A),$$ we 
  shall see that we get all the conclusions of the Theorem rather quickly. To
  prove the smoothness of $q$, let $(T,{m}_T) \rightarrow (S,
  {m}_S)$ be a small Artin surjection with kernel ${\mathfrak a}$.
  Let $B_S \rightarrow A_S$ be a (flat and graded) deformation of $B
  \rightarrow A$ to $S$ and let $B_T$ be a deformation of $B_S$ to $T$. It
  suffices to find a deformation $A_T$ of $A_S$ to $T$ and a map $B_T
  \rightarrow A_T$ over $B_S \rightarrow A_S$. Let $I_{A_S/B_S} = \ker(B_S
  \rightarrow A_S)$.

Firstly we {\it claim} that there exists a graded deformation
$\phi_S:M_{1S}(-s) \rightarrow I_{A_S/B_S}$ (of $S$-flat $B_S$-modules) of
$\phi:M_{1}(-s) \rightarrow I_{A/B}$ to $S$ where $\phi$ composed with
$I_{A/B} \hookrightarrow B$ is the map $\sigma$ in the definition of $A$.
Indeed we can by induction suppose there is a graded deformation
$\phi_{S_1}:M_{1S_1}(-s) \rightarrow I_{A_{S_1}/B_{S_1}}$ of $\phi$ to $S_1$
where $S \rightarrow S_1$ is small Artin surjection (with kernel ${\mathfrak
  a_1}$). Composing $\phi_{S_1}$ with $I_{A_{S_1}/B_{S_1}} \rightarrow
B_{S_1}$ we
deduce by the %resp. 
assumption of $(iii)$ the existence of a deformation $M'_{1S}$ of $M_{1S_1}$ to
$S$. 
    
By Remark~\ref{unobstrrem} the existence of a homogeneous map
$\phi'_{S}:M'_{1S}(-s) \rightarrow I_{A_{S}/B_{S}}$ such that $\phi'_{S}
\otimes_S id_{S_1}=\phi_{S_1}$ is equivalent to the vanishing of a well
defined element (obstruction) $o_0(M'_{1S}, I_{A_S/B_S}) \in \ 
_0\!\Ext_B^1(M_1,I_{A/B}(s)) \otimes_k {\mathfrak a_1}$. Since
$$_0\!\Ext_B^1(M_1,M_1) \otimes_k {\mathfrak a_1} \stackrel{\psi}{\rightarrow}
\ _0\!\Ext_B^1(M_1,I_{A/B}(s)) \otimes_k {\mathfrak a_1} \rightarrow \ 
_0\!\Ext_B^2(M_1,K_1) \otimes_k {\mathfrak a_1}$$
is exact, there is by $(ii)$
an element $\lambda \in \ _0\!\Ext_B^1(M_1,M_1) \otimes {\mathfrak a_1}$ such
that $\psi(\lambda)=o_0(M'_{1S}, I_{A_S/B_S})$. On the other hand, by
deformation theory, one knows more generally that $M'_{1S}- \lambda$ defines a
graded deformation $M_{1S}$ of $M_{1S_1}$ such that
$\psi(\lambda)=o_0(M'_{1S}, I_{A_S/B_S})-o_0(M_{1S}, I_{A_S/B_S})$ (analogous
to \cite{K03}, last part of Remark 3).  Hence $o_0(M_{1S}, I_{A_S/B_S})=0$ for
some deformation $M_{1S}$ and the {\rm claim} is proved.

The composition of $\phi_{S}:M_{1S}(-s) \rightarrow I_{A_{S}/B_{S}}$ with
$I_{A_{S}/B_{S}} \rightarrow B_{S}$ yields a homogeneous map $\sigma_S$ such
that $\sigma_S \otimes_S id_k= \sigma$ where $id_k:k \rightarrow k$ is the
identity. By $(iii)$ there is a deformation $\sigma_T:M_{1T}(-s) \rightarrow
B_T$ over $ \sigma_S$.

Secondly we {\it claim} that $A_T:= \coker \sigma_T $ is a graded
deformation of $A_S$ to $T$, i.e. that the surjection
$I_{A_{T}/B_{T}}\otimes_T k \twoheadrightarrow I_{A/B}$ is an isomorphism
(where $I_{A_{T}/B_{T}}:= \im \sigma_T$). Indeed since $A_T \otimes_TS
\simeq A_S$, it follows that $A_T$ is a deformation of $A_S$ if it is T-flat,
i.e. if $\Tor_1^T(A_T,k) =0$ or equivalently, if $I_{A_{T}/B_{T}}\otimes_T k
\simeq I_{A/B}$. The rough idea for proving this isomorphism of ideals is just
to see that the following part of the Koszul complex
\begin{equation} \label{kosz}
 \wedge^2 M_{1}(-2s) \longrightarrow M_{1}(-s) \longrightarrow B
 \longrightarrow A 
\end{equation}
(induced by the section $\sigma$), which obviously commutes with the
corresponding complex, $ \wedge ^2 M_{1T}(-2s) \rightarrow M_{1T}(-s)
\rightarrow B_T \rightarrow A_T$ over $T$ (because $\sigma_T$ commutes with
$\sigma$), is exact (or ``exact enough'', cf. below for a precise
formulation).

More precisely since $\sigma$ is a regular section on $U$, the Koszul complex
induced by $\sigma$ is {\it exact on} $U$ (\cite{KP2}, Thm.\!\ $7(4)$).
Applying $\HH^0_*(U,-)$ to this Koszul complex (which is really what we did in
the proof of Theorem~\ref{mainth} in \cite{KP2}), we get in particular the
following part of a complex
\begin{equation} \label{kos}
 M_{2}(-2s) \longrightarrow M_{1}(-s) \longrightarrow B \longrightarrow A
 \longrightarrow 0 
\end{equation}
where $ M_{1}(-s) \longrightarrow B \longrightarrow A
 \longrightarrow 0$ is exact. % because $\depth_{I(Z)}B \geq 2$. 
 Using correspondingly $\sigma_T$ we get a complex $ M_{2T}(-2s) \rightarrow
 M_{1T}(-s) \rightarrow B_T \rightarrow A_T$ over $T$ which commutes with
 \eqref{kos}, where $M_{2T}:=\HH^0_*(U, \wedge^2 \widetilde M_{1T})$ (slightly
 abusing the notation of $U$ by letting $U$ be the set in $\Proj(B_T)$
 which corresponds to $U \subseteq \Proj(B)$).  Moreover note that $ M_{1T}
 \stackrel{\simeq}{\longrightarrow} \HH^0_*(U,\widetilde M_{1T})$ follows from
 $\HH^0_*(U,\widetilde M_{1}) \simeq M_{1}$ and the fact that $M_{1T}$ is a
 deformation of $M_1$. Indeed since we by induction may suppose
 $\HH^0_*(U,\widetilde M_{1S}) \simeq M_{1S}$, we conclude easily by applying
 $\HH^0_*(U,-)$ to the exact sequence
\begin{equation} \label{smalldef}
 0 \rightarrow  \widetilde M_{1} \arrowvert_U  \otimes_k
 {\mathfrak a} \rightarrow  
 \widetilde M_{1T} \arrowvert_U  \rightarrow \widetilde M_{1S} \arrowvert_U 
 \rightarrow 0,
\end{equation}
and by comparing with $ 0 \rightarrow M_{1} \otimes_k {\mathfrak a}
\rightarrow M_{1T} \rightarrow M_{1S} \rightarrow 0.$ Let $Z_i = \ker
(M_i(-is) \rightarrow M_{i-1}((1-i)s))$ for $i \geq 1$ and let $Z_{1T}:=
\ker(M_{1T}(-s) \rightarrow I_{A_T/B_T})$. Applying $(-) \otimes_Tk$ to $$
Z_{1T} \longrightarrow M_{1T}(-s) \longrightarrow I_{A_T/B_T} \longrightarrow
0,$$
we see that the exactness of $ Z_1 \rightarrow M_1(-s) \rightarrow
I_{A/B} \rightarrow 0$ and the isomorphism $ M_{1T} \otimes_T k \simeq M_{1}$
imply the claim provided we can prove that $ Z_{1T} \rightarrow Z_{1}$ is
surjective.
 
Now we will show that $ \HH^1_*(U,\widetilde Z_1)=0$ implies the surjectivity
of $ Z_{1T} \rightarrow Z_{1}$. To prove this we remark that the $T$-flatness
of $\widetilde Z_{1T} \arrowvert_U$ (which is true because the claim is true
locally in $U$) yields an exact sequence \eqref{smalldef} in which we have
replaced every $\widetilde M_1$ by $\widetilde Z_1$. Applying $\HH^0_*(U,-)$ to
such an exact sequence, we get the exact sequence
$$Z_{1T} \rightarrow Z_{1S} \rightarrow \ \HH^1_*(U,\widetilde Z_1) \otimes_k
{\mathfrak a}.$$
because  we have $\HH^0_*(U,\widetilde Z_{1T}) \simeq Z_{1T}$ 
from $\HH^0_*(U,\widetilde M_{1T}) \simeq M_{1T}$. If $
\HH^1_*(U,\widetilde Z_1)=0$, we get the surjectivity of $ Z_{1T} \rightarrow
Z_{1S}$ and hence the surjectivity of $ Z_{1T} \rightarrow Z_{1}$ by
induction.

Hence it suffices to prove $\HH^1_*(U,\widetilde Z_1)=0$. Taking cohomology of
the sequence $0 \rightarrow Z_i \rightarrow M_i(-is) \rightarrow Z_{i-1}
\rightarrow 0$, the assumption $(i)$ leads to injections
\begin{equation} \label{incl}
 \HH^1_*(U, \widetilde Z_1)
\hookrightarrow \HH^2_*(U, \widetilde Z_2)\hookrightarrow...\hookrightarrow
\HH^{r-1}_*(U, \widetilde Z_{r-1})
\end{equation}
and we get the vanishing $ \HH^1_*(U, \widetilde Z_1)=0$ because $ \widetilde
Z_{r-1} \arrowvert_U \simeq \widetilde M_{r}(-rs) \arrowvert_U $ and
$\HH^{r-1}_*(U, \widetilde M_{r})=0$ and the second {\it claim} is proved.

Combining the two claims we get the smoothness of the projection 
$q:\GradAlg(H_B,H_A) \rightarrow \GradAlg^{H_B}(R)$ at $(B \rightarrow A)$.

From \eqref{incl} we have $\HH_{I(Z)}^2(Z_1) \simeq \HH^1_*(U, \widetilde
Z_1)=0$ and using the exact sequence
\begin{equation} \label{exaZ_1}   
 0 \rightarrow Z_1 \rightarrow M_{1}(-s) \rightarrow B
 \rightarrow A \rightarrow 0
\end{equation} 
we get $\HH_{I(Z)}^1(I_{A/B})=0$ and $\HH_{I(Z)}^0(A)=0$. It follows that
$\depth_{\mathfrak{m}}A \geq 1$.

Once we have the smoothness of $q$ and the assumption $(iv)$, all remaining
conclusions of Theorem~\ref{unobstr} follow by exactly the same proof as in
the proof of Thm.\!\ 5B) in \cite{K03}. Indeed $(iv)$ implies that the {\it
  second} projection $p:\GradAlg(H_B,H_A) \rightarrow \GradAlg^{H_A}(R)$ is
smooth at $(B \rightarrow A)$ (cf. \cite{K03}, Prop.\!\ 4$(ii)$ and
Remark~\ref{mre}$(i)$ to conclude that the assumption $\
_0\!\Ext_R^1(I_B,I_{A/B})=0$ also implies the smoothness of $p$). Since the
tangent space $T_{B \rightarrow A}$ of $\GradAlg(H_B,H_A)$ at $(B \rightarrow
A)$ is given by the {\it cartesian} square in the following diagram of exact
sequences
\begin{equation} \label{square}
\begin{array}{ccccccccc}
    & & & & & &  _0\! \Hom_R(I_{B},I_{A/B}) &  \\
    & & & & &  & \downarrow &  &  \\
      & & & & T_{B \rightarrow A}
    & \stackrel{T_q}{\rightarrow} &  _0\! \Hom_R(I_{B},B) &  \\
    & & & & \ \ \downarrow \  & \square  & \downarrow &  &  \\
   0 & \rightarrow & _0\! \Hom_R(I_{A/B},A) & \rightarrow &
    _0\! \Hom_R(I_{A},A) & \rightarrow &  _0\! \Hom_R(I_{B},A)
   &  &
\end{array}
\end{equation}
(\cite{K98}, (10)) where the top vertical map is injective and the tangent map
$T_q$ is surjective by the smoothness of $q$, we get
$$\dim_{(B \rightarrow A)} \GradAlg(H_B,H_A)=\dim T_{B \rightarrow A}=
\ \dim_{(B)} \GradAlg^{H_B}(R)+ \ _0\!\hom_R(I_{A/B},A) $$ 
as well as the other conclusions (by e.g. using that the
tangent map $T_p :T_{B \rightarrow A} \rightarrow \ _0\! \Hom_R(I_{A},A)$
of $p$ is surjective by $(iv)$ and \eqref{square}). 
\end{proof}

\begin{remark}  \label{remunobstr4}
  If we, instead of $ \ (i) \ \ \ \HH^i_*(U,\widetilde M_{i+1})=0 \ \ 
  {\rm for}\ \ 1\leq i \leq r-1$, assume
  $$(i') \ \ \HH^1_*(U,\widetilde M_2)=0 \ \ {\rm and} \ \ 
  \HH^i_*(U,\widetilde M_{i+2})=0 \ \ {\rm for} \ \ 1\leq i \leq r-2,$$
  and
  keep the other assumptions of Theorem~\ref{unobstr}, we still get all
  conclusions of Theorem~\ref{unobstr}. Indeed, looking to the proof of (in
  particular the second claim of) Theorem~\ref{unobstr}, it suffices to prove
  that $ M_{2T} \rightarrow M_{2}$ is surjective and that \eqref{kos} is
  exact.  Now applying $\HH^0_*(U, \wedge^2-)$ to the exact sequence
  \eqref{smalldef}, remarking that $\widetilde M_2 \arrowvert_U \simeq
  \wedge^2 \widetilde M_{1} \arrowvert_U$, we get the exact sequence $M_{2T}
  \rightarrow M_{2S} \rightarrow \ \HH^1_*(U,\widetilde M_2) \otimes_k
  {\mathfrak a}.$ Since we have $ \HH^1_*(U,\widetilde M_2)=0$ by the
  assumption $(i')$, we get the surjectivity of $ M_{2T} \rightarrow M_{2S}$
  and hence the surjectivity of $ M_{2T} \rightarrow M_{2}$ by induction.
  Finally to show that \eqref{kos} is exact, it suffices to show the
  surjectivity of $\HH^0_*(U,\widetilde M_2(-2s)) \rightarrow
  \HH^0_*(U,\widetilde Z_1)$ or the vanishing of $ \HH^1_*(U, \widetilde
  Z_2)$.  Since the conditions of $(i')$ lead to inclusions as in \eqref{incl}
  provided we have replaced $Z_i$ by $Z_{i+1}$, we get precisely $
  \HH^1_*(U,\widetilde Z_2)=0$ and we are done.
\\
Finally note that the assumption $ \ \HH^1_*(U,\widetilde M_2)=0 $ in $(i')$
is superfluous if we in proving the second claim above can show that there is a
deformation $M'_{2T}$ of $M_2$ which locally on $U$ is $M_{2T}$ for any $T$.
Indeed the argument in the proof of Theorem~\ref{unobstr} where we showed
$\HH^0_*(U,\widetilde M_{1T}) \simeq M_{1T}$ as consequence of
$\HH^0_*(U,\widetilde M_{1}) \simeq M_{1}$ and the fact that $M_{1T}$ is a
deformation of $M_1$, apply to get $M'_{2T} \simeq M_{2T}$ and hence the
surjectivity of $ M_{2T} \rightarrow M_{2}$.
\end{remark}

If we apply Theorem~\ref{unobstr} under the assumptions of
Theorem~\ref{mainth}, we get Theorem~\ref{gorgenth} stated in the background
section. Here we include a proof of Theorem~\ref{gorgenth}.

\begin{proof}
  It suffices to verify $(i)$, $(ii)$ and $(iv)$ of Theorem~\ref{unobstr}.
  Firstly we show $(ii)$. Note that \eqref{LESTHM} is given by applying $
  \HH^0_*(U,-)$ onto the Koszul resolution induced by the regular
  section $\sigma$. Now splitting the exact sequence \eqref{LESTHM} into short
  exact sequences, we get $0 \rightarrow Z_i \rightarrow M_i(-is) \rightarrow
  Z_{i-1} \rightarrow 0$ where $Z_i = \ker (M_i(-is) \rightarrow
  M_{i-1}((1-i)s))$ for $i \geq 2$ and $Z_1(s)=K_1$. Hence
  $_0\!\Ext^2_B(M,K_1)=0$ provided $_0\!\Ext^2_B(M,M_2(-s))=0$ and
  $_0\!\Ext^3_B(M,Z_2(s))=0$, while $_0\!\Ext^3_B(M,Z_2(s))=0$ provided
  $_0\!\Ext^3_B(M,M_3(-2s))=0$ and $_0\!\Ext^4_B(M,Z_3(s))=0$ etc. By the
  assumption $(i)$ of Theorem~\ref{gorgenth}, it suffices to show
  $_0\!\Ext^r_B(M,Z_{r-1}(s))= 0$. Since $Z_{r-1} = M_r(-rs) = K_B(t-rs)$, cf.
  Theorem~\ref{mainth}, we get this vanishing by Gorenstein duality. Hence we
  conclude by the assumption $(i)$ of Theorem~\ref{gorgenth}. Moreover $(iv)$
  of Theorem~\ref{unobstr} is proven by the same argument, using $(iii)$ of
  Theorem~\ref{gorgenth} (cf. Remark~\ref{remcorvarunobstr}), and remarking
  that $Ass(I_{A/B}) \subset Ass(B)$. Note that the argument shows that $ \
  _0\!\Ext^i_R(I_B,M_i(-is))= 0$ for $1 \leq i \leq r$ implies $
  _0\!\Ext_R^1(I_B,I_{A/B})=0$ and we conclude as required.
 
  It remains to show $(i)$ of Theorem~\ref{unobstr} or $(i')$ of
  Remark~\ref{remunobstr4}. If $\depth_{I(Z)}B = \dim B- \dim {B}/{I}(Z) \geq
  r+1$, we have $$ \HH^{i}_*(U,\widetilde M_{i+1}) \simeq
  \HH^{i+1}_{{I}(Z)}(M_{i+1}) = 0 \ \ {\rm for} \ \ 1 \leq i \leq r-1$$ since
  all ${M}_i$ are maximal CM ${B}$-modules by Theorem~\ref{mainth}. If
  $\depth_{I(Z)}B = r$ and $r \geq 3$ we similarly verify $(i')$ of
  Remark~\ref{remunobstr4}. In the final case $\depth_{I(Z)}B = r=2$ we can
  use the two final sentences in Remark~\ref{remunobstr4} to conclude. Indeed
  in this case $M_2 = K_B(t)$ and since $ \Ext^{j}_{R}(B,R(-n-c))=0$ for $j
  \neq c$ and $K_{B}= \Ext^{c}_{R}(B,R(-n-c))$ we may take $M'_{2T} =
  K_{B_T}:= \Ext^{c}_{R_T}(B_T,R_T(-n-c))$ and prove that $K_{B_T}$ is a
  deformation of $K_{B}$ to $B_T$ by e.g. \cite{JS}, Prop.\! (A1).
   %(as in \cite{KMMNP}, proof of Prop.\!\ 10.7). 
  This concludes the proof.  
\end{proof}

{\small

\bibliographystyle{amsalpha}

}

\bigskip
\bigskip

\end{document}